\catcode'32=9
\magnification=1200

\voffset=1cm

\font\tenpc=cmcsc10

\font\eightrm=cmr8
\font\eighti=cmmi8
\font\eightsy=cmsy8
\font\eightbf=cmbx8
\font\eighttt=cmtt8
\font\eightit=cmti8
\font\eightsl=cmsl8
\font\sixrm=cmr6
\font\sixi=cmmi6
\font\sixsy=cmsy6
\font\sixbf=cmbx6

\skewchar\eighti='177 \skewchar\sixi='177
\skewchar\eightsy='60 \skewchar\sixsy='60

\font\tengoth=eufm10
\font\tenbboard=msbm10
\font\eightgoth=eufm7 at 8pt
\font\eightbboard=msbm7 at 8pt
\font\sevengoth=eufm7
\font\sevenbboard=msbm7
\font\sixgoth=eufm5 at 6 pt
\font\fivegoth=eufm5

\font\tengoth=eufm10
\font\tenbboard=msbm10
\font\eightgoth=eufm7 at 8pt
\font\eightbboard=msbm7 at 8pt
\font\sevengoth=eufm7
\font\sevenbboard=msbm7
\font\sixgoth=eufm5 at 6 pt
\font\fivegoth=eufm5

\newfam\gothfam
\newfam\bboardfam

\catcode`\@=11

\def\raggedbottom{\topskip 10pt plus 36pt
\r@ggedbottomtrue}
\def\pc#1#2|{{\bigf@ntpc #1\penalty
\@MM\hskip\z@skip\smallf@ntpc #2}}

\def\tenpoint{%
  \textfont0=\tenrm \scriptfont0=\sevenrm \scriptscriptfont0=\fiverm
  \def\rm{\fam\z@\tenrm}%
  \textfont1=\teni \scriptfont1=\seveni \scriptscriptfont1=\fivei
  \def\oldstyle{\fam\@ne\teni}%
  \textfont2=\tensy \scriptfont2=\sevensy \scriptscriptfont2=\fivesy
  \textfont\gothfam=\tengoth \scriptfont\gothfam=\sevengoth
  \scriptscriptfont\gothfam=\fivegoth
  \def\goth{\fam\gothfam\tengoth}%
  \textfont\bboardfam=\tenbboard \scriptfont\bboardfam=\sevenbboard
  \scriptscriptfont\bboardfam=\sevenbboard
  \def\bboard{\fam\bboardfam}%
  \textfont\itfam=\tenit
  \def\it{\fam\itfam\tenit}%
  \textfont\slfam=\tensl
  \def\sl{\fam\slfam\tensl}%
  \textfont\bffam=\tenbf \scriptfont\bffam=\sevenbf
  \scriptscriptfont\bffam=\fivebf
  \def\bf{\fam\bffam\tenbf}%
  \textfont\ttfam=\tentt
  \def\tt{\fam\ttfam\tentt}%
  \abovedisplayskip=12pt plus 3pt minus 9pt
  \abovedisplayshortskip=0pt plus 3pt
  \belowdisplayskip=12pt plus 3pt minus 9pt
  \belowdisplayshortskip=7pt plus 3pt minus 4pt
  \smallskipamount=3pt plus 1pt minus 1pt
  \medskipamount=6pt plus 2pt minus 2pt
  \bigskipamount=12pt plus 4pt minus 4pt
  \normalbaselineskip=12pt
  \setbox\strutbox=\hbox{\vrule height8.5pt depth3.5pt width0pt}%
  \let\bigf@ntpc=\tenrm \let\smallf@ntpc=\sevenrm
  \let\petcap=\tenpc
  \normalbaselines\rm}
\def\eightpoint{%
  \textfont0=\eightrm \scriptfont0=\sixrm \scriptscriptfont0=\fiverm
  \def\rm{\fam\z@\eightrm}%
  \textfont1=\eighti \scriptfont1=\sixi \scriptscriptfont1=\fivei
  \def\oldstyle{\fam\@ne\eighti}%
  \textfont2=\eightsy \scriptfont2=\sixsy \scriptscriptfont2=\fivesy
  \textfont\gothfam=\eightgoth \scriptfont\gothfam=\sixgoth
  \scriptscriptfont\gothfam=\fivegoth
  \def\goth{\fam\gothfam\eightgoth}%
  \textfont\bboardfam=\eightbboard \scriptfont\bboardfam=\sevenbboard
  \scriptscriptfont\bboardfam=\sevenbboard
  \def\bboard{\fam\bboardfam}%
  \textfont\itfam=\eightit
  \def\it{\fam\itfam\eightit}%
  \textfont\slfam=\eightsl
  \def\sl{\fam\slfam\eightsl}%
  \textfont\bffam=\eightbf \scriptfont\bffam=\sixbf
  \scriptscriptfont\bffam=\fivebf
  \def\bf{\fam\bffam\eightbf}%
  \textfont\ttfam=\eighttt
  \def\tt{\fam\ttfam\eighttt}%
  \abovedisplayskip=9pt plus 2pt minus 6pt
  \abovedisplayshortskip=0pt plus 2pt
  \belowdisplayskip=9pt plus 2pt minus 6pt
  \belowdisplayshortskip=5pt plus 2pt minus 3pt
  \smallskipamount=2pt plus 1pt minus 1pt
  \medskipamount=4pt plus 2pt minus 1pt
  \bigskipamount=9pt plus 3pt minus 3pt
  \normalbaselineskip=9pt
  \setbox\strutbox=\hbox{\vrule height7pt depth2pt width0pt}%
  \let\bigf@ntpc=\eightrm \let\smallf@ntpc=\sixrm
  \normalbaselines\rm}

\tenpoint

\frenchspacing


\newif\ifpagetitre
\newtoks\auteurcourant \auteurcourant={\hfil}
\newtoks\titrecourant \titrecourant={\hfil}

\def\appeln@te{}
\def\vfootnote#1{\def\@parameter{#1}\insert\footins\bgroup\eightpoint
  \interlinepenalty\interfootnotelinepenalty
  \splittopskip\ht\strutbox 
  \splitmaxdepth\dp\strutbox \floatingpenalty\@MM
  \leftskip\z@skip \rightskip\z@skip
  \ifx\appeln@te\@parameter\indent \else{\noindent #1\ }\fi
  \footstrut\futurelet\next\fo@t}

\pretolerance=500 \tolerance=1000 \brokenpenalty=5000
\newdimen\hmargehaute \hmargehaute=0cm
\newdimen\lpage \lpage=13.3cm
\newdimen\hpage \hpage=20cm
\newdimen\lmargeext \lmargeext=1cm
\hsize=11.25cm
\vsize=18cm
\parskip 0pt
\parindent=12pt

\def\margehaute{\vbox to \hmargehaute{\vss}}%
\def\margebasse{\vss}

\output{\shipout\vbox to \hpage{\margehaute\nointerlineskip
  \corpsdepage\margebasse}
  \advancepageno \global\pagetitrefalse
  \ifnum\outputpenalty>-20000 \else\dosupereject\fi}

\def\corpsdepage{\hbox to \lpage{\hss\pagetexte\hskip\lmargeext}}
\def\pagetexte{\vbox{\makeheadline\pagebody\makefootline}}
\headline={\ifpagetitre\titleheadline \else
  \ifodd\pageno\rightheadline \else\leftheadline\fi\fi}
\def\leftheadline{\eightpoint\hfil\the\auteurcourant\hfil}
\def\rightheadline{\eightpoint\hfil\the\titrecourant\hfil}
\def\titleheadline{\hfill}
\pagetitretrue

\def\footnoterule{\kern-6\p@
  \hrule width 2truein \kern 5.6\p@} 

\def\pd#1#2 {\pc#1#2| }

\def\pointir{\discretionary{.}{}{.\kern.35em---\kern.7em}\nobreak
\hskip 0em plus .3em minus .4em }

\def\abstract#1{\vbox{\eightpoint \pc ABSTRACT|\pointir #1}}

\def\titre#1|{\message{#1}
              \par\vskip 30pt plus 24pt minus 3pt\penalty -1000
              \vskip 0pt plus -24pt minus 3pt\penalty -1000
              \centerline{\bf #1}
              \vskip 5pt
              \penalty 10000 }

\def\section#1|{\par\vskip .3cm
                {\bf #1}\pointir}

\def\ssection#1|{\par\vskip .2cm
                {\it #1}\pointir}

\long\def\th#1|#2\finth{\par\medskip
              {\petcap #1\pointir}{\it #2}\par\smallskip}

\long\def\tha#1|#2\fintha{\par\medskip
                    {\petcap #1.}\par\nobreak{\it #2}\par\smallskip}

\def\rem#1|{\par\medskip
            {{\it #1}.\quad}}

\def\rema#1|{\par\medskip
             {{\it #1.}\par\nobreak }}

\def\article#1|#2|#3|#4|#5|#6|#7|
    {{\leftskip=7mm\noindent
     \hangindent=2mm\hangafter=1
     \llap{[#1]\hskip.35em}{#2}.\quad
     #3, {\sl #4}, vol.\nobreak\ {\bf #5}, {\oldstyle #6},
     p.\nobreak\ #7.\par}}
\def\livre#1|#2|#3|#4|
    {{\leftskip=7mm\noindent
    \hangindent=2mm\hangafter=1
    \llap{[#1]\hskip.35em}{#2}.\quad
    {\sl #3}.\quad #4.\par}}
\def\divers#1|#2|#3|
    {{\leftskip=7mm\noindent
    \hangindent=2mm\hangafter=1
     \llap{[#1]\hskip.35em}{#2}.\quad
     #3.\par}}
\mathchardef\conj="0365
\def\proof{\par{\it Proof}.\quad}
\def\qed{\quad\raise -2pt\hbox{\vrule\vbox to 10pt{\hrule width 4pt
\vfill\hrule}\vrule}}

\def\decale#1|{\par\noindent\hskip 28pt\llap{#1}\kern 5pt}

\catcode`\@=12


\catcode`\@=11
\def\matrice#1{\null \,\vcenter {\normalbaselines \m@th
\ialign {\hfil $##$\hfil &&\  \hfil $##$\hfil\crcr
\mathstrut \crcr \noalign {\kern -\baselineskip } #1\crcr
\mathstrut \crcr \noalign {\kern -\baselineskip }}}\,}

\def\petitematrice#1{\left(\null\vcenter {\normalbaselines \m@th
\ialign {\hfil $##$\hfil 
&&\thinspace  \hfil $##$\hfil\crcr
\mathstrut \crcr \noalign {\kern -\baselineskip } #1\crcr
\mathstrut \crcr \noalign {\kern -\baselineskip }}}\right)}

\catcode`\@=12

\def\qed{\quad\raise -2pt\hbox{\vrule\vbox to 10pt{\hrule width 4pt
   \vfill\hrule}\vrule}}


%

\def\bfc{\mathop{\bf c}}


\def\il{\bigl]\kern-.25em\bigl]}
\def\ir{\bigr]\kern-.25em\bigr]}

\def\iil{\bigl>\kern-.25em\bigl>}
\def\iir{\bigr>\kern-.25em\bigr>}



\def\Grille{\setbox13=\vbox to 5mm{\hrule width 110mm\vfill}
\setbox13=\vbox{\offinterlineskip
   \copy13\copy13\copy13\copy13\copy13\copy13\copy13\copy13
   \copy13\copy13\copy13\copy13\box13\hrule width 110mm}
\setbox14=\hbox to 5mm{\vrule height 65mm\hfill}
\setbox14=\hbox{\copy14\copy14\copy14\copy14\copy14\copy14
   \copy14\copy14\copy14\copy14\copy14\copy14\copy14\copy14
   \copy14\copy14\copy14\copy14\copy14\copy14\copy14\copy14\box14}
\ht14=0pt\dp14=0pt\wd14=0pt
\setbox13=\vbox to 0pt{\vss\box13\offinterlineskip\box14}
\wd13=0pt\box13}


\def\fleche(#1,#2)\dir(#3,#4)\long#5{%
\noalign{\nointerlineskip\leftput(#1,#2){\vector(#3,#4){#5}}\nointerlineskip}}


\def\hfl#1#2#3{\smash{\mathop{\hbox to#3{\rightarrowfill}}\limits
^{\scriptstyle#1}_{\scriptstyle#2}}}

\def\gfl#1#2#3{\smash{\mathop{\hbox to#3{\leftarrowfill}}\limits
^{\scriptstyle#1}_{\scriptstyle#2}}}


 \message{`lline' & `vector' macros from LaTeX}
 \catcode`@=11
\def\{{\relax\ifmmode\lbrace\else$\lbrace$\fi}
\def\}{\relax\ifmmode\rbrace\else$\rbrace$\fi}
\def\newcount{\alloc@0\count\countdef\insc@unt}
\def\newdimen{\alloc@1\dimen\dimendef\insc@unt}
\def\newwrite{\alloc@7\write\chardef\sixt@@n}

\newwrite\@unused
\newcount\@tempcnta
\newcount\@tempcntb
\newdimen\@tempdima
\newdimen\@tempdimb
\newbox\@tempboxa

\def\@spaces{\space\space\space\space}
\def\@whilenoop#1{}
\def\@whiledim#1\do #2{\ifdim #1\relax#2\@iwhiledim{#1\relax#2}\fi}
\def\@iwhiledim#1{\ifdim #1\let\@nextwhile=\@iwhiledim
        \else\let\@nextwhile=\@whilenoop\fi\@nextwhile{#1}}
\def\@badlinearg{\@latexerr{Bad \string\line\space or \string\vector
   \space argument}}
\def\@latexerr#1#2{\begingroup
\edef\@tempc{#2}\expandafter\errhelp\expandafter{\@tempc}%
\def\@eha{Your command was ignored.
^^JType \space I <command> <return> \space to replace it
  with another command,^^Jor \space <return> \space to continue without it.}
\def\@ehb{You've lost some text. \space \@ehc}
\def\@ehc{Try typing \space <return>
  \space to proceed.^^JIf that doesn't work, type \space X <return> \space to
  quit.}
\def\@ehd{You're in trouble here.  \space\@ehc}

\typeout{LaTeX error. \space See LaTeX manual for explanation.^^J
 \space\@spaces\@spaces\@spaces Type \space H <return> \space for
 immediate help.}\errmessage{#1}\endgroup}
\def\typeout#1{{\let\protect\string\immediate\write\@unused{#1}}}

\font\tenln    = line10
\font\tenlnw   = linew10

\newdimen\@wholewidth
\newdimen\@halfwidth
\newdimen\unitlength 

\unitlength =1pt


\def\thinlines{\let\@linefnt\tenln \let\@circlefnt\tencirc
  \@wholewidth\fontdimen8\tenln \@halfwidth .5\@wholewidth}
\def\thicklines{\let\@linefnt\tenlnw \let\@circlefnt\tencircw
  \@wholewidth\fontdimen8\tenlnw \@halfwidth .5\@wholewidth}

\def\linethickness#1{\@wholewidth #1\relax \@halfwidth .5\@wholewidth}

\newif\if@negarg

\def\lline(#1,#2)#3{\@xarg #1\relax \@yarg #2\relax
\@linelen=#3\unitlength
\ifnum\@xarg =0 \@vline
  \else \ifnum\@yarg =0 \@hline \else \@sline\fi
\fi}

\def\@sline{\ifnum\@xarg< 0 \@negargtrue \@xarg -\@xarg \@yyarg -\@yarg
  \else \@negargfalse \@yyarg \@yarg \fi
\ifnum \@yyarg >0 \@tempcnta\@yyarg \else \@tempcnta -\@yyarg \fi
\ifnum\@tempcnta>6 \@badlinearg\@tempcnta0 \fi
\setbox\@linechar\hbox{\@linefnt\@getlinechar(\@xarg,\@yyarg)}%
\ifnum \@yarg >0 \let\@upordown\raise \@clnht\z@
   \else\let\@upordown\lower \@clnht \ht\@linechar\fi
\@clnwd=\wd\@linechar
\if@negarg \hskip -\wd\@linechar \def\@tempa{\hskip -2\wd\@linechar}\else
     \let\@tempa\relax \fi
\@whiledim \@clnwd <\@linelen \do
  {\@upordown\@clnht\copy\@linechar
   \@tempa
   \advance\@clnht \ht\@linechar
   \advance\@clnwd \wd\@linechar}%
\advance\@clnht -\ht\@linechar
\advance\@clnwd -\wd\@linechar
\@tempdima\@linelen\advance\@tempdima -\@clnwd
\@tempdimb\@tempdima\advance\@tempdimb -\wd\@linechar
\if@negarg \hskip -\@tempdimb \else \hskip \@tempdimb \fi
\multiply\@tempdima \@m
\@tempcnta \@tempdima \@tempdima \wd\@linechar \divide\@tempcnta \@tempdima
\@tempdima \ht\@linechar \multiply\@tempdima \@tempcnta
\divide\@tempdima \@m
\advance\@clnht \@tempdima
\ifdim \@linelen <\wd\@linechar
   \hskip \wd\@linechar
  \else\@upordown\@clnht\copy\@linechar\fi}

\def\@hline{\ifnum \@xarg <0 \hskip -\@linelen \fi
\vrule height \@halfwidth depth \@halfwidth width \@linelen
\ifnum \@xarg <0 \hskip -\@linelen \fi}

\def\@getlinechar(#1,#2){\@tempcnta#1\relax\multiply\@tempcnta 8
\advance\@tempcnta -9 \ifnum #2>0 \advance\@tempcnta #2\relax\else
\advance\@tempcnta -#2\relax\advance\@tempcnta 64 \fi
\char\@tempcnta}

\def\vector(#1,#2)#3{\@xarg #1\relax \@yarg #2\relax
\@linelen=#3\unitlength
\ifnum\@xarg =0 \@vvector
  \else \ifnum\@yarg =0 \@hvector \else \@svector\fi
\fi}

\def\@hvector{\@hline\hbox to 0pt{\@linefnt
\ifnum \@xarg <0 \@getlarrow(1,0)\hss\else
    \hss\@getrarrow(1,0)\fi}}

\def\@vvector{\ifnum \@yarg <0 \@downvector \else \@upvector \fi}

\def\@svector{\@sline
\@tempcnta\@yarg \ifnum\@tempcnta <0 \@tempcnta=-\@tempcnta\fi
\ifnum\@tempcnta <5
  \hskip -\wd\@linechar
  \@upordown\@clnht \hbox{\@linefnt  \if@negarg
  \@getlarrow(\@xarg,\@yyarg) \else \@getrarrow(\@xarg,\@yyarg) \fi}%
\else\@badlinearg\fi}

\def\@getlarrow(#1,#2){\ifnum #2 =\z@ \@tempcnta='33\else
\@tempcnta=#1\relax\multiply\@tempcnta \sixt@@n \advance\@tempcnta
-9 \@tempcntb=#2\relax\multiply\@tempcntb \tw@
\ifnum \@tempcntb >0 \advance\@tempcnta \@tempcntb\relax
\else\advance\@tempcnta -\@tempcntb\advance\@tempcnta 64
\fi\fi\char\@tempcnta}

\def\@getrarrow(#1,#2){\@tempcntb=#2\relax
\ifnum\@tempcntb < 0 \@tempcntb=-\@tempcntb\relax\fi
\ifcase \@tempcntb\relax \@tempcnta='55 \or
\ifnum #1<3 \@tempcnta=#1\relax\multiply\@tempcnta
24 \advance\@tempcnta -6 \else \ifnum #1=3 \@tempcnta=49
\else\@tempcnta=58 \fi\fi\or
\ifnum #1<3 \@tempcnta=#1\relax\multiply\@tempcnta
24 \advance\@tempcnta -3 \else \@tempcnta=51\fi\or
\@tempcnta=#1\relax\multiply\@tempcnta
\sixt@@n \advance\@tempcnta -\tw@ \else
\@tempcnta=#1\relax\multiply\@tempcnta
\sixt@@n \advance\@tempcnta 7 \fi\ifnum #2<0 \advance\@tempcnta 64 \fi
\char\@tempcnta}

\def\@vline{\ifnum \@yarg <0 \@downline \else \@upline\fi}

\def\@upline{\hbox to \z@{\hskip -\@halfwidth \vrule
  width \@wholewidth height \@linelen depth \z@\hss}}

\def\@downline{\hbox to \z@{\hskip -\@halfwidth \vrule
  width \@wholewidth height \z@ depth \@linelen \hss}}

\def\@upvector{\@upline\setbox\@tempboxa\hbox{\@linefnt\char'66}\raise
     \@linelen \hbox to\z@{\lower \ht\@tempboxa\box\@tempboxa\hss}}

\def\@downvector{\@downline\lower \@linelen
      \hbox to \z@{\@linefnt\char'77\hss}}

\thinlines

\newcount\@xarg
\newcount\@yarg
\newcount\@yyarg
\newcount\@multicnt
\newdimen\@xdim
\newdimen\@ydim
\newbox\@linechar
\newdimen\@linelen
\newdimen\@clnwd
\newdimen\@clnht
\newdimen\@dashdim
\newbox\@dashbox
\newcount\@dashcnt
 \catcode`@=12


\newbox\tbox
\newbox\tboxa

\def\leftzer#1{\setbox\tbox=\hbox to 0pt{#1\hss}%
     \ht\tbox=0pt \dp\tbox=0pt \box\tbox}

\def\rightzer#1{\setbox\tbox=\hbox to 0pt{\hss #1}%
     \ht\tbox=0pt \dp\tbox=0pt \box\tbox}

\def\centerzer#1{\setbox\tbox=\hbox to 0pt{\hss #1\hss}%
     \ht\tbox=0pt \dp\tbox=0pt \box\tbox}

%
\def\image(#1,#2)#3{\vbox to #1{\offinterlineskip
    \vss #3 \vskip #2}}


\def\leftput(#1,#2)#3{\setbox\tboxa=\hbox{%
    \kern #1\unitlength
    \raise #2\unitlength\hbox{\leftzer{#3}}}%
    \ht\tboxa=0pt \wd\tboxa=0pt \dp\tboxa=0pt\box\tboxa}

\def\rightput(#1,#2)#3{\setbox\tboxa=\hbox{%
    \kern #1\unitlength
    \raise #2\unitlength\hbox{\rightzer{#3}}}%
    \ht\tboxa=0pt \wd\tboxa=0pt \dp\tboxa=0pt\box\tboxa}

\def\centerput(#1,#2)#3{\setbox\tboxa=\hbox{%
    \kern #1\unitlength
    \raise #2\unitlength\hbox{\centerzer{#3}}}%
    \ht\tboxa=0pt \wd\tboxa=0pt \dp\tboxa=0pt\box\tboxa}

\unitlength=1mm

\def\put(#1,#2)#3{\noalign{\nointerlineskip
                               \centerput(#1,#2){$#3$}
                                \nointerlineskip}}
\def\segment(#1,#2)\dir(#3,#4)\long#5{%
\leftput(#1,#2){\lline(#3,#4){#5}}}
\auteurcourant={DOMINIQUE FOATA {\sevenrm AND} GUO-NIU HAN}
\titrecourant={ANDR\'E PERMUTATION CALCULUS}

\def\grn{\mathop{\bf grn}\nolimits}

\def\And{\mathop{\rm And}\nolimits}
\def\AndI{\mathop{\rm And}\nolimits^{I}}  
\def\AndII{\mathop{\rm And}\nolimits^{I\!I}}

\def\brullet{{\scriptscriptstyle\bullet}}
\def\Deltaa{\mathop{\hbox{$\Delta$}}\limits}
\def\bfF{{\bf F}\,}
\def\bfL{{\bf L}\,}

\def\NL{{\bf N\!L}}
\def\spi{\mathop{\bf spi}\nolimits}
\def\pit{\mathop{\bf pit}\nolimits}

\def\twodigit#1#2{#1\hskip -1pt#2}
\def\mapright#1{\smash{\mathop{\longrightarrow}\limits^{#1}}}
\def\esp{\hskip-7.8pt}

\rightline{January 14, 2016}
\bigskip\bigskip
\centerline{\bf Andr\'e Permutation Calculus; a Twin Seidel Matrix Sequence}
\bigskip
\centerline{\sl Dominique Foata and Guo-Niu Han}
\footnote{}{
{\it Key words and phrases.} 
Entringer numbers, 
tangent and secant numbers, alternating permutations, Andr\'e permutations, 
Seidel matrix sequence, increasing binary trees, 
greater neighbor of maximum, spike, pit,
tight permutations, hooked permutations, Seidel triangle sequence,
formal Laplace transform.
\par
{\it Mathematics Subject Classifications.} 
05A15, 05A30, 11B68, 33B10.}

\bigskip\bigskip

{\narrower\narrower
\eightpoint 

\noindent
{\bf Abstract}.\quad 
Entringer numbers occur in the Andr\'e permutation combinatorial set-up under several forms.
This leads to the construction of a matrix-analog refinement of the tangent (resp. secant) numbers. Furthermore, closed expressions for the
three-variate exponential generating functions for pairs of so-called Entringerian statistics are derived.

}
\bigskip

{\eightpoint

\halign{\indent\hfil #.&\ #\hfil\cr
1&Introduction\cr
\omit&1. Entringer Numbers\cr
\omit&2. Andr\'e Permutations\cr
\omit&3. Statistics on Andr\'e Permutations\cr
\omit&4. A further bijection between Andr\'e Permutations\cr
\omit&5. The Twin Seidel Matrix Sequence\cr
\omit&6. Tight and hooked permutations\cr
\omit&7. Trivariate Generating Functions\cr
2&From Alternating to Andr\'e Permutations of the first kind\cr
3&The bijection $\phi$ between Andr\'e I and Andr\'e II permutations\cr
4&The bijection $g$ of the set of Andr\'e I permutations onto itself\cr
5&The proof of Theorem~1.1 (iii) and (iv)\cr
6&Combinatorics of the Twin Seidel Matrix Sequence\cr
\omit&1. The first evaluations\cr
\omit&2. Tight and not tight Andr\'e I Permutations\cr
\omit&3. Hooked and unhooked Permutations\cr
\omit&4. A Bijection of $B_n(m,k)$ onto $N\!H_n(m+1,k)$\cr
7&The Making of Seidel Triangle Sequences\cr
\omit&1. The Seidel Tangent-Secant Matrix\cr
\omit&2. The generating function for the Entringer numbers\cr
\omit&3. Seidel Triangle Sequences\cr
8&Trivariate Generating Functions \cr
\omit&1. The Upper Triangles of ${\rm Twin}^{(1)}$\cr
\omit&2. The Upper Triangles of ${\rm Twin}^{(2)}$\cr
\omit&3. The Bottom rows of the Matrices $B_n$\cr
\omit&4. The Lower Triangles of ${\rm Twin}^{(1)}$\cr
\omit&5. The Lower Triangles of ${\rm Twin}^{(2)}$\cr
9&The Formal Laplace Transform\cr
\omit&References\cr
}

}

\bigskip

\vfill\eject

\centerline{\bf 1. Introduction}

\medskip
The notion of {\it alternating} or {\it zizag} permutation devised by D\'esir\'e Andr\'e, back in 1881 [An1881], for interpreting the coefficients $E(n)$ $(n\ge 0)$ of the Taylor expansion of $\tan u+\sec u$, the so-called {\it tangent} and {\it secant numbers}, has remained some sort of a curiosity for a long time, until it was realized that the geometry of those alternating permutations could be exploited to obtain further arithmetic refinements of those numbers. Classifying alternating permutations according to the number of inversions directly leads to the constructions of their $q$-analogs (see [AF80, AG78, St76]). Sorting them according to their first letters led Entringer [En66]  to obtain a fruitful refinement $E_n=\sum_mE_n(m)$ that has been described under several forms [OEIS, GHZ11, KPP94, MSY96, St10], the entries $E_n(m)$ satisfying a simple finite difference equation (see (1.1) below).

In fact, those numbers $E_n(m)$, further called {\it Entringer numbers}, appear in other contexts, in particular when dealing with analytical properties of the {\it Andr\'e permutations}, of the two kinds~I and~II, introduced by Sch\"utzenberger and the first author ([FSch73, FSch71]). For each $n\ge 1$ let $\And_n^I$ (resp. $\And_n^{I\!I}$) be the set of all Andr\'e permutations of $12\cdots n$ (see \S1.2). It was shown that $\#\And_n^I=\#\And_n^{I\!I}=E_n$. The first purpose of this paper is to show
that there are several natural statistics ``stat,''  defined on $\And_n^I$ (resp. $\And_n^{I\!I}$), whose distributions are {\it Entringerian}, that is,
 integer-valued mappings ``stat,''  satisfying $\#\{w\in \And_n^I\
({\rm resp.}\  \And_n^{I\!I}):{\rm stat}(w)=m\}=E_n(m)$.

The second purpose is to work out a matrix-analog refinement $E_n=\sum_{m,k}a_n(m,k)$ of the tangent and secant numbers, whose 
row and column sums $\sum_ka_n(m,k)$ and $\sum_ma_n(m,k)$  are themselves refinements
 of the Entringer numbers.
This will be achieved, first by inductively defining the so-called {\it twin Seidel matrix sequence} $(A_n,B_n)$ $(n\ge 2)$ (see \S1.5), then by proving that the entries of those matrices
provide the joint distributions of pairs of Entringerian statistics defined on Andr\'e permutations of each kind (Theorem~1.2). See \S1.6 for the plan of action.

The third purpose is to obtain analytical expressions for the joint exponential generating functions for pairs of those Entringerian statistics. See \S1.7 and the contents of Section~7 and~8. Let us give more details on the notions introduced so far.

\goodbreak
\medskip

1.1. {\it Entringer numbers}.\quad
According to D\'esir\'e Andr\'e [An1879, An1881] each permutation $w=x_{1}x_{2}\cdots x_{n}$ of $12\cdots n$ is said to be {\it (increasing) alternating} if $x_{1}<x_{2}$, $x_{2}>x_{3}$, $x_{3}<x_{4}$, etc. in an alternating way. Let ${\rm Alt}_{n}$ be the set of all alternating permutations of $12\cdots n$. He then proved that
$\#\,{\rm Alt}_{n}=E_{n}$, where $E_{n}$ is  the {\it tangent number} (resp. {\it secant number}) when $n$ is odd (resp. even), those numbers
 appearing in the Taylor expansions of $\sec u$ and $\tan u$:
$$
\leqalignno{ \noalign{\vskip-5pt}
\tan u&=\sum_{n\ge 1} {u^{2n-1}\over
(2n-1)!}
E_{2n-1}
={u\over 1!}1\!+\!{u^3\over 3!}2\!+\!{u^5\over 5!}16\!+\!{u^7\over
7!}272\!+\! {u^9\over 9!}7936\!+\!\cdots\cr
\sec u&=\sum_{n\ge 0}
{u^{2n}\over (2n)!}E_{2n}
=1\!+\!{u^2\over 2!}1\!+\!{u^4\over 4!}5\!+\!{u^6\over 6!}61\!+\!{u^8\over
8!}1385\!+\! {u^{10}\over 10!}50521\!+\!\cdots\cr 
}
$$
(See, e.g., [Ni23, p.~177-178], [Co74, p.~258-259]).

Let $\bfF w:=x_{1}$ be the {\it first} letter of a permutation $w=x_{1}x_{2}\cdots x_{n}$ of $12\cdots n$. 
For each $m=1,\ldots, n$,
the {\it Entringer numbers} are defined by
$E_n(m):=\#\{w\in {\rm Alt}_{n}:\bfF w=m\}$,
as was introduced by Entringer [En66].
In particular, $E_n(n)=0$ for $n\geq 2$.
He showed that those numbers satisfied the recurrence:
$$\leqalignno{\noalign{\vskip-5pt}
E_1(1):=1;\quad E_n(n)&:=0\ {\rm for\ all\ }n\ge 2;&(1.1)\cr
\Delta E_n(m)+E_{n-1}(n-m)&=0\quad
(n\geq 2; m=n-1, \ldots, 2,1);\cr
\noalign{\vskip-5pt}}
$$
where $\Delta$ stands for the classical finite difference operator (see, e.g. [Jo39])
$\Delta E_n(m):=E_n(m+1)-E_n(m)$. See Fig. 1.1 for the table of their first values.
Those numbers are registered as the $A008282$ sequence 
in  Sloane's On-Line Encyclopedia of Integer Sequences, 
together with an abundant bibliography [OEIS].
They naturally constitute a refinement of the tangent and secant numbers:
$$
\sum_mE_n(m)=E_n=\cases{\hbox{tangent number},&if $n$ is odd;\cr
\hbox{secant number},&if $n$ is even.\cr}\leqno(1.2)
$$

$$
\vbox{\halign{\vrule\ \hfil$#$\ \vrule
&\strut\ \hfil$#$\hfil
&\ \hfil$#$\hfil
&\ \hfil$#$\hfil
&\ \hfil$#$\hfil
&\ \hfil$#$\hfil
&\ \hfil$#$\hfil
&\ \hfil$#$\hfil
&\ \hfil$#$\hfil
&\ \hfil$#$\hfil\ \vrule
&\ \hfil$#$\ \vrule\cr
\noalign{\hrule}
m={}&1&2&3&4&5&6&7&8&9&{\rm Sum}\cr
\noalign{\hrule}
n=1&1&&&&&&&&&1\cr
2&1&0&&&&&&&&1\cr
3&1&1&0&&&&&&&2\cr
4&2&2&1&0&&&&&&5\cr
5&5&5&4&2&0&&&&&16\cr
6&16&16&14&10&5&0&&&&61\cr
7&61&61&56&46&32&16&0&&&272\cr
8&272&272&256&224&178&122&61&0&&1385\cr
9&1385&1385&1324&1202&1024&800&544&272&0&7936\cr
\noalign{\hrule}
}}
$$
\centerline{Fig. 1.1. The Entringer Numbers $E_{n}(m)$}

\bigskip
Now, let  $\bfL w:=x_{n}$ denote the {\it last} letter of a permutation $w=x_{1}x_{2}\cdots x_{n}$ of $12\cdots n$.  In our previous paper [FH14] we made a full study of the so-called {\it Bi-Entringer numbers} defined by
$$E_n(m,k):=\#\{w\in {\rm Alt}_n: \bfF w=m,\,\bfL w=k\},
$$
and showed that the sequence of the matrices $(E_n(m,k)_{1\le m,k\le n})\;(n\!\ge\! 1)$ was fully determined by a partial difference equation system and the three-variable exponential generating function for those matrices could be calculated.
As the latter analytical derivation essentially depends on the {\it geometry} of alternating permutations, it is natural to ask whether other combinatorial models, counted by tangent and secant numbers, are likely to have a parallel development.

Let $E(u):=\tan u+\sec u=\sum_{n\ge 0}(u^n/n!)\,E_n$. Then, the first and second derivatives of $E(u)$ are equal to: $E'(u)=E(u)\,\sec u$ and $E''(u)=E(u)\,E'(u)$, two identities equivalent to the two recurrence relations:
$$\displaylines{\rlap{(*)}\hfill
E_{n+1}=\sum_{0\le 2j\le n}{n\choose 2j}E_{n-2j}E_{2j}\quad (n\ge 0),\qquad E_0=1;\hfill\cr
\rlap{(**)}\hfill
E_{n+2}=\sum_{0\le j\le n}{n\choose j}E_j\,E_{n+1-j}\quad (n\ge 0),\qquad
E_0=E_1=1.\hfill\cr}
$$
The first of those relations can be readily interpreted in terms of {\it alternating} permutations, or in terms of the so-called {\it Jacobi} permutations introduced by Viennot [Vi80]. The second one leads naturally to the model of {\it Andr\'e} permutations, whose geometry will appear to be rich and involves several analytic developments.

\medskip
1.2. {\it Andr\'e permutations}.\quad Those permutations were introduced in [FSch73, FSch71], and further studied in [Str74, FSt74, FSt76]. Other properties have been developed in the works by Purtill [Pu93], Hetyei [He96], Hetyei and Reiner [HR98], the present authors [FH01], Stanley [St94], in particular in the study of the   $cd$-index in a Boolean algebra. More recently, Disanto [Di14] has been able to calculate the joint distribution of the right-to-left minima and left-to-right minima in those permutations.

In the sequel, permutations of a finite set $Y=\{ y_1<y_2<\cdots<y_n\}$ of positive integers will be written as {\it words} $w=x_1x_2\cdots x_n$, where the {\it letters}~$x_i$ are the elements of~$Y$
in some order. The {\it minimum} (resp. {\it maximum}) letter of~$w$, in fact, $y_1$ (resp. $y_n$), will be denoted by $\min(w)$ (resp. $\max (w)$). When writing $w=v\,\min(w)\,v'$ it is meant that the word~$w$ is the juxtaposition product of the left factor~$v$, followed by the letter $\min(w)$, then by the right factor~$v'$.
\medskip
{\it Definition}.\quad
Say that the empty word~$e$ and each one-letter word are both {\it Andr\'e~I} and {\it Andr\'e~II permutations}.
Next, if $w=x_1x_2\cdots x_n$ ($n\geq 2$) is a permutation of a set of positive integers $Y=\{ y_1<y_2<\cdots<y_n\}$, write $w=v\,\min(w)\,v'$. Then, $w$ is said to be an {\it Andr\'e~I} (resp. {\it Andr\'e~II}) {\it permutation} if 
both $v$ and $v'$ are themselves  Andr\'e~I (resp. Andr\'e~II) permutations, and furthermore if~$\max(vv')$ (resp.~$\min(vv')$) is a letter of~$v'$.

\medskip
The set of all Andr\'e I (resp. Andr\'e II) permutations of $Y$ is denoted by $\AndI_{Y}$ (resp. $\AndII_{Y}$), and simply by
$\AndI_{n}$ (resp. $\AndII_{n}$) when $Y=\{1,2,\ldots, n\}$.
In the sequel, an Andr\'e I (resp. Andr\'e II) permutation, with no reference to a set $Y$,
is meant to be an element of 
$\AndI_{n}$ (resp. $\AndII_{n}$).

Using such an inductive definition we can immediately see that $E_n=\#\And_n^I=
\#\And_n^{I\!I}$, the term ${n\choose j}E_j\,E_{n+1-j}$ in $(**)$ being the number of all Andr\'e I (resp. Andr\'e II) permutations of $x_1x_2\cdots x_{n+2}$ such that $x_j=1$. Further equivalent definitions will be given in the beginning of Section~2. The first Andr\'e permutations 
from $\AndI_n$ and $\AndII_n$ are listed in Table~1.2.

\smallskip
Andr\'e permutations of the first kind:

$n=1$:\quad 1;\qquad $n=2$:\quad 12;\qquad
$n=3$:\quad 123, 213;

$n=4$:\quad 1234, 1324, 2314, 2134, 3124;

$n=5$:\quad 12345, 12435, 13425, 23415, 13245, 14235, 34125, 
24135,\hfil\break
\indent\hphantom{$n=5$:\quad}23145, 21345, 41235, 31245, 21435, 32415, 41325, 31425.

\smallskip
Andr\'e permutations of the second kind:

$n=1$:\quad 1;\qquad $n=2$:\quad 12;\qquad
$n=3$:\quad 123, 312;

$n=4$:\quad 1234, 1423, 3412, 4123, 3124;

$n=5$:\quad 12345, 12534, 14523, 34512, 15234, 14235, 34125, 
45123,\hfil\break
\indent\hphantom{$n=5$:\quad}35124, 51234, 41235, 31245, 51423, 53412, 41523, 31524.

\medskip
\centerline{Table 1.2: the first Andr\'e permutations of both kinds}

\medskip
1.3. {\it Statistics on Andr\'e permutations}.\quad 
The statistics ``$\bfF$'' and ``$\bfL$''
have been previously introduced. Two further ones are now defined: the {\it next to the last {\rm (or the penultimate)} letter} ``$\NL$'' and {\it {\bf gr}eater {\bf n}eighbor of the maximum} ``$\grn$'': for $w=x_1x_2\cdots x_n$ and $n\ge 2$ let
$\NL\,w:=x_{n-1}$; next, let $x_{i}=n$ for a certain~$i$ $(1\le i\le n)$
with the convention that $x_{0}=x_{n+1}:=0$. Then, $\grn w:=\max\{x_{i-1}, x_{i+1}\}$.

Let $({\rm Ens}_n)$ $(n\ge 1)$ be a sequence of non-empty finite sets and ``stat''
an integer-valued mapping $w\mapsto {\rm stat}(w)$
defined on each ${\rm Ens}_n$. The pair $({\rm Ens}_n,{\rm stat})$  is said to be {\it Entringerian}, if  $\#{\rm Ens}_n=E_n$ and  $\#\{w\in {\rm Ens}_n:{\rm stat}(w)=m\}=E_n(m)$ holds for each $m=0,1,\ldots, n$. We also say that ``stat'' is an {\it Entringerian statistic}. The pair $({\rm Alt}_n,\bfF)$ is Entringerian, {\it par excellence}, for all $n\ge 1$. 

\proclaim Theorem 1.1. For each $n\ge 2$ the mappings\hfil\break
\indent {\rm (i)} {\bf F} defined on $\And_n^I$,\hfil\break
\indent {\rm (ii)} $n-\NL$ defined on $\And_n^I$,\hfil\break
\indent {\rm (iii)} $(n+1)-\bfL$ defined on $\And_n^{I\!I}$,\hfil\break
\indent {\rm (iv)} $n-\grn$ defined on $\And_n^{I\!I}$,\hfil\break
are all Entringerian statistics.

\vfill\eject
\def\fleche(#1,#2)\dir(#3,#4)\long#5{%
{\leftput(#1,#2){\vector(#3,#4){#5}}}}

Statements (i) and (ii) will be proved in Section 2 by constructing two bijections $\eta$ and $\theta$ having the property
$$\matrice{
{\rm Alt}_n&\buildrel \eta\over \longrightarrow& \And_n^I 
&\buildrel \theta\over \longrightarrow &\And_n^I \cr
w&\mapsto &w'&\mapsto &w''\cr
\bfF w&=&\bfF w'&=&(n-\NL)w''\cr}\leqno(1.3)
$$
For proving (iii) and (iv) we use the properties of two new bijections $\phi:\And_n^I\rightarrow \And_n^{I\!I}$ and $g:\And_n^I\rightarrow \And_n^I$, whose constructions are described in Sections~3 and~4.  By means of those two bijections, as well as the bijection~$\theta$ mentioned in (1.3),  it will be shown in Section~5 that the following properties hold
\vskip-18pt
$$
\matrice{
\And_n^I\quad &\buildrel \strut \phi\,\circ\,g \over \longrightarrow& \And_n^{I\!I}\cr
w&\mapsto &w'\cr
\bfF w&=&(n+1-\bfL)w'\cr}
\leqno(1.4)
$$
\vskip-10pt
\noindent
and
\vskip-18pt
$$
\matrice{
\And_n^I \quad &\buildrel \strut \phi\,\circ\, \theta \over \longrightarrow& \And_n^{I\!I}\cr
w&\mapsto &w''\cr
\bfF w&=&(n-\grn)w''\cr
}\leqno(1.5)
$$
thereby completing the proof of Theorem 1.1.

\medskip
1.4. {\it The fundamental bijection $\phi$}.\quad
For proving (1.4) and (1.5) and also the next Theorem~1.2 two new statistics are to be introduced, the {\it spike} ``$\spi$'' and the {\it pit} ``$\pit$'', related to the {\it left} minimum records for the former one, and the {\it right} minimum records for the latter one. In Section~3 the bijection $ \phi$ of 
$\And_n^I$ onto $\And_n^{I\!I}$ will be shown to have the further property:
$$
(\bfF,\spi,\NL) w= (\pit,\bfL,\grn) \phi(w).\leqno(1.6)
$$
This implies that 

\smallskip
\noindent
(1.7) {\it for each pair $(m,k)$ 
the two sets
$
\{w\in\And_n^I:(\spi,\NL)w=(m,k)\}$ and $\{w\in\And_n^{I\!I}:(\bfL,\grn)w=(m,k)\}$
are equipotent.}

\smallskip

It also follows from Theorem 1.1 that  ``$(n+1)-\spi$'' on $\And_n^I$ and ``$\pit$'' on $\And_n^{I\!I}$ are two further Entringerian statistics.

\medskip

1.5. {\it The twin Seidel matrix sequence}.\quad
The next step is to say something about the joint distributions of the pairs
$(\bfF,\NL)$ on $\And_n^I$ and $(\bfL,\grn)$ on $\And_n^{I\!I}$, whose marginal distributions are Entringerian, as announced in Theorem~1.1. We shall proceed in the following way: first,
the notion of  {\it twin Seidel matrix sequence} $(A_n,B_n)$ $(n\ge 2)$ will be introduced (see Definition below), then the entry in cell $(m,k)$ of $A_n$ (resp. $B_n$) will be shown to be the number of Andr\'e I (resp. II) permutations~$w$, whose values $(\bfF,\NL)w$ (resp. $(\bfL,\grn)w$)
are equal to $(m,k)$. The definition involves the {\it partial difference operator}
$\Deltaa_{(1)}$ acting on sequences $(a_n(m,k))$ $(n\ge 2)$ of integers depending on two integral
variables $m$, $k$ as follows:
$$
\Deltaa_{(1)} a_n(m,k):= a_n(m+1,k)-a_n(m,k).
$$
The subscript (1) indicates that the difference operator is to be applied to the variable occurring at the {\it first} position, which is `$m$' in the previous equation.
\medskip

{\it Definition}.\quad
The {\it twin Seidel matrix sequence} $(A_n, B_n)$ $(n\ge 2)$ is a sequence of finite square 
matrices that obey the following five rules  (TS1)--(TS5) (see Diagram 1.3 for the values of the first  matrices, where null entries are replaced by dots):

(TS1) each matrix $A_n=(a_n(m,k))$  (resp. $B_n=b_n(m,k)$) $(1\le m,k\le n)$ is a square matrix of dimension~$n$ $(n\ge 2)$ with nonnegative entries, and zero entries along its diagonal, except for $a_2(1,1)=1$; let $a_n(m,\brullet)=\sum_ka_n(m,k)$ (resp. $a_n(\brullet,k)
=\sum_ma_n(m,k)$) be the 
$m$-th row sum (resp. $k$-th column sum) of the matrix~$A_n$ with an analogous notation for~$B_n$;

(TS2) for $n\ge 3$ the entries along the rightmost column in both $A_n$ and $B_n$ are null, as well as the entries in the bottom row of~$A_n$ and the top row of~$B_n$, i.e.,
$a_n(\brullet,n)=b_n(\brullet,n)=a_n(n,\brullet)=b_n(1,\brullet)=0$, as all the entries are supposed to be nonnegative; furthermore, $b_n(n,1)=0$;

(TS3) the first two matrices of the sequence are supposed to be:
$A_2=\matrice{1&\cdot\cr \cdot&\cdot\cr}$\quad ,\quad
$B_2=\matrice{\cdot&\cdot\cr 1&\cdot\cr }$\quad;

(TS4) for each $n\ge 3$ the matrix $B_n$ is derived from the matrix~$A_{n-1}$
by means of a transformation~$\Psi:(a_{n-1}(m,k))\rightarrow
(b_n(m,k))$ defined as follows
$$
\leqalignno{
b_n(n,k)&:=a_{n-1}(\brullet,k-1)\quad (2\le k\le n-1);&\hbox{(TS4.1)}\cr
b_n(n-1,k)&:=a_{n-1}(\brullet,k)\quad (2\le k\le n-2);&\hbox{(TS4.2)}\cr
}$$
and, by induction,
$$
\leqalignno{
\qquad \Deltaa_{(1)}b_n(m,k)&-a_{n-1}(m,k)=0\quad(2\le k+1\le m\le n-2);&\hbox{(TS4.3)}\cr
\Deltaa_{(1)}b_n(m,k)&-a_{n-1}(m,k-1) =0
\quad (3\le m+2\le k\le n-1);&\hbox{(TS4.4)}\cr
\cr}
$$

\vskip-10pt 
(TS5) for each $n\ge 3$ the matrix $A_n$ is derived from the matrix~$B_{n-1}$
by means of a transformation~$\Phi: (b_{n-1}(m,k))\rightarrow
(a_n(m,k))$ defined as follows
$$
\leqalignno{
a_n(1,k)&:=b_{n-1}(\brullet,k-1)\quad (2\le k\le n-1);&\hbox{(TS5.1)}\cr}
$$
\vfill\eject

\noindent and, by induction,
$$
\leqalignno{
\qquad \Deltaa_{(1)}a_n(m,k)&+b_{n-1}(m,k-1)=0\quad(3\le m+2\le k\le n-1);&\hbox{(TS5.2)}\cr
\Deltaa_{(1)}a_n(m,k)&+b_{n-1}(m,k)=0
\quad(2\le k+1\le m\le n-1).&\hbox{(TS5.3)}\cr
}
$$

\vskip-25pt
{\eightpoint

$$\displaylines{
\def\mapright#1{\smash{\mathop{\longrightarrow}\limits^{#1}}}
\def\esp{\hskip-5pt}
  \matrix{&&&&&&&&&&\cr
\esp  A_2\esp&=\matrice{1&\cdot\cr \cdot&\cdot\cr}&\esp\mapright\Psi
 &\esp B_3\esp&=\matrice{\cdot&\cdot&\cdot\cr
 1&\cdot&\cdot\cr
 \cdot&1&\cdot\cr}
 &\esp\mapright\Phi
 &\esp A_4\esp&=\matrice{\cdot&1&1&\cdot\cr
1&\cdot& 1&\cdot\cr
\cdot&1&\cdot&\cdot\cr
 \cdot&\cdot&\cdot&\cdot\cr}
 &\esp\mapright\Psi
  &\esp B_5\esp&=\matrice{\cdot&\cdot&\cdot&\cdot&\cdot\cr
\cdot&\cdot& 1&1&\cdot\cr
1&1&\cdot&2&\cdot\cr
1&2&2&\cdot&\cdot\cr
\cdot&1&2&2&\cdot\cr}
 &\esp\mapright\Phi\hfill
\cr}
\hfill\cr
\noalign{\vskip-5pt}
 \matrix{&&&&&&&\cr
 \esp\mapright\Phi&
\esp A_6\esp&=\matrice{\cdot&2&4&5&5&\cdot\cr
2&\cdot&4&5&5&\cdot\cr
2&4&\cdot&4&4&\cdot\cr
1&3&4&\cdot&2&\cdot\cr
\cdot&1&2&2&\cdot&\cdot\cr
\cdot&\cdot&\cdot&\cdot&\cdot&\cdot\cr}
&\esp\mapright\Psi
 &\esp B_7\esp&=\matrice{\cdot&\cdot&\cdot&\cdot&\cdot&\cdot&\cdot\cr
\cdot&\cdot& 2&4&5&5&\cdot\cr
2&2&\cdot&8&10&10&\cdot\cr
4&6&8&\cdot&14&14&\cdot\cr
5&9&12&14&\cdot&16&\cdot\cr
5&10&14&16&16&\cdot&\cdot\cr
\cdot&5&10&14&16&16&\cdot\cr
}
 &\esp\mapright\Phi
 &\esp A_8\esp&=\matrice{
   \cdot &  16 &  32 &  46 &  56 &  61 &  61 &\cdot \cr
 16 &   \cdot &  32 &  46 &  56 &  61 &  61 & \cdot\cr
16 &  32 &   \cdot &  44 &  52 &  56 &  56 &\cdot \cr
14 &  30 &  44 &   \cdot &  44 &  46 &  46 & \cdot\cr
10 &  24 &  36 &  44 &   \cdot &  32 &  32 & \cdot\cr
5 &  15 &  24 &  30 &  32 &   \cdot &  16 &\cdot \cr
\cdot  &   5 &  10 &  14 &  16 &  16 &  \cdot& \cdot\cr
\cdot&\cdot&\cdot&\cdot&\cdot&\cdot&\cdot&\cdot\cr} \cr
 }\cr}
 $$
 
\vskip-.5cm
$$\displaylines{
\def\mapright#1{\smash{\mathop{\longrightarrow}\limits^{#1}}}
\def\esp{\hskip-5pt}
\matrix{&&&&&&&&&&&\cr
\esp B_2\esp
&=\matrice{\cdot&\cdot\cr 1&\cdot\cr}
&\esp\mapright\Phi
&\esp  A_3\esp
&=\matrice{\cdot&1&\cdot\cr
 1&\cdot&\cdot\cr
 \cdot&\cdot&\cdot\cr}
&\esp\mapright\Psi
&\esp B_4\esp
&=\matrice{\cdot&\cdot&\cdot&\cdot\cr
  \cdot&\cdot& 1&\cdot\cr
  1&1&\cdot&\cdot\cr
  \cdot&1&1&\cdot\cr}
&\esp\mapright\Phi
&\esp A_5\esp
&=\matrice{\cdot&1&2&2&\cdot\cr
     1&\cdot&2&2&\cdot\cr
     1&2&\cdot&1&\cdot\cr
     \cdot&1&1&\cdot&\cdot\cr
     \cdot&\cdot&\cdot&\cdot&\cdot\cr}
&\esp\mapright\Psi\cr}\hfill\cr
 \noalign{\vskip-5pt}
\matrix{&&&&&&&\cr
\esp\mapright\Psi
&\esp B_6\esp
&=\matrice{   \cdot &   \cdot &   \cdot &   \cdot &\cdot &\cdot \cr
          \cdot &   \cdot &   1 &   2&   2 &\cdot \cr
         1 &   1 &   \cdot &   4 &   4 &\cdot \cr
          2 &   3 &  4& \cdot &   5 &  \cdot \cr
          2 &   4 &   5 &   5 &\cdot&\cdot \cr
      \cdot& 2 &   4 &   5 &   5 &\cdot\cr}
&\esp\mapright\Phi
&\esp A_7\esp
&=\matrice{
   \cdot &   5 &  10 &  14 &  16 &  16 &\cdot \cr
  5 &   \cdot &  10 &  14 &  16 &  16 &\cdot \cr
    5 &  10 &   \cdot &  13 &  14 &  14 &\cdot \cr
    4 &   9 &  13 &   \cdot &  10 &  10 & \cdot\cr
     2 &   6 &   9 &  10 &   \cdot &   5 &\cdot \cr
     \cdot &   2 &   4 &   5 &   5 &   \cdot &\cdot \cr
     \cdot&\cdot&\cdot&\cdot&\cdot&\cdot&\cdot\cr}
&\esp\mapright\Psi
&\esp B_8\esp
&=\matrice{
   \cdot &   \cdot &   \cdot &   \cdot &   \cdot &   \cdot & 
\cdot&\cdot\cr
   \cdot & \cdot&  5 &  10 &  14 &  16&    16 & \cdot\cr
   5 &  5& \cdot &  20 &  28 &  32 &    32 & \cdot\cr
  10 &  15&20 &   \cdot &  41 &    46 &46&\cdot \cr
  14 &  24 &  33&41 &   \cdot &  56 &    56 &\cdot \cr
  16 &  30 &  42 &  51 & 56&  \cdot &  61 & \cdot \cr
  16 &  32 &  46 &  56 &  61 &   61 &\cdot &\cdot\cr
  \cdot&16 &  32 &  46 &  56 &  61 &   61 &\cdot \cr}
\cr}\cr } $$

}
\centerline{Diagram 1.3: First values of the twin Seidel matrices}
\bigskip
It is worth noting that the twin Seidel matrix sequence involves two infinite subsequences:
${\rm Twin}^{(1)}=(A_2,B_3,A_4,B_5,A_6,\ldots\,)$ and
${\rm Twin}^{(2)}=(B_2,A_3,B_4,A_5,B_6,\ldots\,)$. They are independent in the sense that the matrices $A_{2n}$ (resp. $B_{2n}$) depend only on the matrices $B_{2m+1}$ and $A_{2m}$ 
(resp. $A_{2m+1}$ and $B_{2m}$) with $m<n$, with an analogous statement for the matrices $A_{2n+1}$ (resp. $B_{2n+1}$).

As easily verified, rules (TS1)--(TS5) define the twin Seidel matrix sequence by induction in a unique manner.  At each step Rules (TS1) and (TS2) furnish all the zero entries indicated by dots and rules (TS4.1), (TS4.2), (TS5.1) the initial values. It remains to use the finite difference equations (TS4.3), (TS4.4), (TS5.2), (TS5.3) to calculate the other entries.

\proclaim Theorem 1.2. The twin Seidel matrix sequence $(A_n=(a_n(m,k)),\,
B_n=(b_n(m,k)))$ $(n\ge 2,\,1\le m,k\le n)$ defined by relations {\rm (TS1)--(TS5)} provides the joint distributions of the pairs
$( \bfF,\NL)$ on $ \And_n^I$ and $(\bfL,\grn)$ on $\And_n^{I\!I}$ in the sense that
for $n\ge 2$ the following relations hold:
$$\leqalignno{
a_n(m,k)&=\#\{w\in \And_n^I:( \bfF,\NL)w=(m,k)\};&(1.8)\cr
b_n(m,k)&=\#\{w\in \And_n^{I\!I}:( \bfL,\grn)w=(m,k)\}.&(1.9)\cr}
$$

By Theorems 1.1 and 1.2 the row and column sums of the matrices $A_n$ and $B_n$ have the following interpretations
$$\leqalignno{\qquad
a_n(m,\brullet)&=E_n(m),\quad b_n(m,\brullet)=E_n(n+1-m),\quad
(1\le m\le n);&(1.10)\cr
a_n(\brullet,k)&=b_n(\brullet,k)=E_n(n-k)\quad (1\le k\le n);&(1.11)\cr
\noalign{\hbox{and furthermore the matrix-analog of the refinement of $E_n$ holds:}}
\qquad
\sum_{m,k}&\,a_n(m,k)=\sum_{m,k}b_n(m,k)=
E_n.&(1.12)\cr}
$$

1.6. {\it Tight and hooked permutations}.\quad
For proving Theorem~1.2 the crucial point is to show that the $a_n(m,k)$'s and $b_n(m,k)$'s satisfy the partial difference equations (TS4.3), (TS4.4), (TS5.2), (TS5.3),
when those numbers are defined by the right-hand sides of (1.8) and (1.9).  For each pair $(m,k)$ let
$$\leqalignno{
A_n(m,k)&:=\{w\in \And_n^I:( \bfF,\NL)w=(m,k)\};\cr
	B_n(m,k)&:=\{w\in \And_n^{I}:(\spi,\grn)w=(m,k)\}.\cr
}$$
As the latter set is equipotent with the set $\{w\in \And_n^{I\!I}:(\bfL,\grn)w=(m,k)\}$ by (1.7), 
we also have 
$a_n(m,k)=\#A_n(m,k)$ and $b_n(m,k)=\#B_n(m,k)$, by (1.8) and (1.9).
  
 For the partial difference equation (TS5.2) (resp. (TS5.3)) the plan of action may be described by the diagram
$$
\matrice{&&B_{n-1}(m,k-1)\ &&&\hskip-2.5cm{\rm (resp.}\ B_{n-1}(m,k)\, )\cr
&&\;\Big\downarrow\phi&\cr
A_n(m,k)&=&T_n(m,k)&+&N\!T_n(m,k)\cr
&&&&\Big\downarrow f\cr
&&&&A_n(m+1,k)\cr}\leqno(1.13)
$$
This means that the set $A_n(m,k)$ is to be split into two disjoint subsets
$A_n(m,k)=T_n(m,k)+N\!T_n(m,k)$ in such a way that the first component is in bijection with $B_{n-1}(m,k-1)$ (resp. $B_{n-1}(m,k)$) by using the bijection~$\phi$ defined in (6.6), and the second one with $A_n(m+1,k)$ by means of the bijection~$f$ defined in (6.5).
If this plan is realized, the above partial difference equations are satisfied, as
$\Deltaa_{(1)}a_n(m,k)=\#A_n(m+1,k)-\#A_n(m,k)
=\#N\!T_n(m,k)-\#A_n(m,k)=-\#T_n(m,k)
=-\#B_{n-1}(m,k-1)\ ({\rm resp.}\ -\#B_{n-1}(m,k)
=-b_{n-1}(m,k-1)\ ({\rm resp.}\ -b_{n-1}(m,k))$.

For the partial difference equation (TS4.3) (resp. (TS4.4)) the corresponding  diagram is the following
$$
\matrice{&&A_{n-1}(m,k)\ &&&\hskip-2.5cm{\rm (resp.}\ A_{n-1}(m,k-1\, )\cr
&&\;\Big\downarrow\Theta&\cr
B_n(m+1,k)&=&H_n(m+1,k)&+&N\!H_n(m+1,k)\cr
&&&&\Big\uparrow \beta\cr
&&&&B_n(m,k)\cr}\leqno(1.14)
$$
where $\Theta$ and $\beta$ are two explicit bijections,
defined in (6.8) and (6.13), (6.14), respectively.
The elements in $T_n(m,k)$ from (1.13) (resp. in $H_n(m+1,k)$ from (1.14)) are the so-called {\it tight} (resp. {\it hooked\/})
{\it permutations.} All details will be given in Section~6 and constitute the bulk of the proof of Theorem~1.2.
\goodbreak
\medskip
1.7. {\it Trivariate generating functions}.\quad The final step is to show that the partial difference equation systems (TS4.3), (TS4.4), (TS5.2), (TS5.3) satisfied by
the twin Seidel matrix sequence $(A_n,B_n)$ $(n\ge 2)$ make it possible to derive {\it closed expressions} for the trivariate generating functions for the sequences $(A_{2n})$, $(A_{2n+1})$, $(B_{2n})$, $(B_{2n+1})$. 
We all list them in the following theorems. See Section 8 for the  detailed proofs.

The calculations are all based on the Seidel Triangle Sequence technique developed in our previous paper [FH14]. Note that the next generating functions for the matrices~$A_n$ do not involve the entries of the rightmost columns and bottom rows, which are all zero;
they do not involve either the entries of the rightmost columns of the matrices~$B_n$,
also equal to zero, as assumed in (TS2). Finally, the generating functions for the bottom rows of the matrices~$B_n$ are calculated separately: see (1.23) and (1.24).

Also, note that the right-hand sides of identities (1.15)--(1.18) are all symmetric with respect to~$x$ and~$z$, in agreement with its combinatorial interpretation stated in Theorem~2.4. The property is less obvious for (1.16), but an easy exercise on trigonometry shows that the right-hand side is equal to the fraction
$\displaystyle{\cos x\cos z\sin(x+y+z)-\sin y\over \cos^2(x+y+z)}$.
Finally, the summations below are taken over triples $\{(m,k,n)\}$ or pairs $\{(k,n)\}$ for the last two ones; only the ranges of the summations have been written.

\vfill\eject
\proclaim Theorem 1.3 {\rm [The sequence $(A_{2n})$ $(n\ge 1)$]}. The generating function for the upper triangles is given by
$$
\displaylines{
(1.15)
\sum_{\scriptstyle 2\le m+1\le k\le 2n-1}
\kern-15pt 
a_{2n}(m,k){x^{m-1}\over (m-1)!}{y^{k-m-1}\over (k-m-1)!}
{z^{2n-k-1}\over (2n-k-1)!}
\hfill\cr
\hfill{}=
{\cos x\,\cos z\sin (x+y+z)\over \cos^2(x+y+z)}\quad\cr
\noalign{\hbox{and for the lower triangles by}}
(1.16)
\sum_{2\le k+1\le m\le 2n-1}
\kern-15pt 
a_{2n}(m,k){x^{2n-m-1}\over (2n-m-1)!}{y^{m-k-1}\over (m-k-1)!}{z^{k-1}\over (k-1)!}
\hfill\cr
\hfill{}={\cos x\,\sin z\over \cos(x+y+z)}
+{\sin x\,\cos (x+y)\over \cos^2(x+y+z)}.\quad\cr
}
$$

\proclaim Theorem 1.4 {\rm [The sequence $(A_{2n+1})$ $(n\ge 1)$]}. The generating function for the upper triangles is given by
$$\displaylines{\noalign{\vskip-5pt}
(1.17)
\sum_{2\le m+1\le k\le 2n} 
a_{2n+1}(m,k){x^{m-1}\over (m-1)!}
{y^{k-m-1}\over (k-m-1)!}{z^{2n-k}\over (2n-k)!}
\hfill\cr
\hfill{}=
{\cos x\cos z\over \cos^2(x+y+z)}\quad\cr
\noalign{\hbox{and for the lower triangles by}}
(1.18)
\sum_{2\le k+1\le m\le 2n}
\kern-20pt 
a_{2n+1}(m,k){x^{2n-m}\over (2n-m)!}{y^{m-k-1}\over (m-k-1)!}{z^{k-1}\over (k-1)!}
\hfill\cr
\hfill{}={\cos(x+y)\,\cos(y+z)\over \cos^2(x+y+z)}.\quad\cr
}$$

\goodbreak

\proclaim Theorem 1.5 {\rm [The sequence $(B_{2n})$ $(n\ge 1)$]}. The generating function for the 
upper triangles is given by
$$
\displaylines{(1.19)
\sum_{2\le m+1\le k\le 2n-1}\kern-20pt 
b_{2n}(m,k) {x^{m-1}\over (m-1)!} {y^{k-m-1}\over (k-m-1)!} {z^{2n-1-k}\over (2n-1-k)!}\hfill\cr
\hfill{}=
{\sin x\cos z\over \cos^2(x+y+z)}\quad\cr
\noalign{\hbox{and for the lower triangles by}}
(1.20)
\sum_{2\le k+1\le m\le 2n-1}
\kern-20pt 
b_{2n}(m,k){x^{2n-m-1}\over (2n-m-1)!}{y^{m-k-1}\over (m-k1)!}{z^{k-1}\over (k-1)!}
\hfill\cr
\hfill{}={\cos(x+y)\,\sin(y+z)\over \cos^2(x+y+z)}.\quad\cr
}$$


\proclaim Theorem 1.6 {\rm [The sequence $(B_{2n+1})$ $(n\ge 1)$]}. The generating function for the 
upper triangles is given by
$$\displaylines{(1.21)
\quad
\sum_{\scriptstyle 2\le m+1\le k\le 2n}
\kern-15pt 
b_{2n+1}(m,k){x^{m-1}\over (m-1)!}{y^{k-m-1}\over (k-m-1)!}{z^{2n-k}\over (2n-k)!}
\hfill\cr
\noalign{\vskip-6pt}
\hfill{}=
{\sin x\,\cos z\sin (x+y+z)\over \cos^2(x+y+z)}\quad\cr
\noalign{\hbox{and for the lower triangles by}}
(1.22)\sum_{2\le k+1\le m\le 2n}
\kern-15pt 
b_{2n+1}(m,k){x^{2n-m}\over (2n-m)!}{y^{m-k-1}\over (m-k-1)!}{z^{k-1}\over (k-1)!}
\hfill\cr
\hfill{}=-{\sin x\,\sin z\over \cos(x+y+z)}
+{\cos x\,\cos (x+y)\over \cos^2(x+y+z)}.\quad\cr
}
$$
 
 The {\it bivariate} generating functions for the bottom rows $b_n(n,k)$ $(k=1,2,\dots\,)$  are computed as follows:
$$\leqalignno{
\sum_{1\le k\le 2n-1}
b_{2n}(2n,k){x^{2n-k-1}\over (2n-k-1)!}{y^{k-1}\over (k-1)!}
	&={\cos x\over \cos(x+y)};&(1.23)\cr
\sum_{1\le k\le 2n}\kern-6pt
b_{2n+1}(2n+1,k){x^{k-1}\over (k-1)!}{y^{2n-k}\over (2n-k)!}
	&={\sin x\over \cos(x+y)}.&(1.24)\cr
}	
$$

The  previous generating functions for the matrices $A_n$, $B_n$ will be derived analytically in Section~8, from the sole definition of twin Seidel matrix sequence given in \S\thinspace 1.5, without reference to any combinatorial interpretation. It will be shown in Section~9 that, conversely, the closed expressions thereby obtained provide an analytical proof of identity (1.12), by means of the formal Laplace transform, that is, the fact that the entries of those matrices
furnish a refinement of the tangent and secant numbers.

\bigskip
\centerline{\bf 2. From alternating to Andr\'e permutations of the first kind}

\medskip
Two further equivalent definitions of Andr\'e permutations of the two kinds 
will be given (see Definitions~2.1 and~2.2). They were actually introduced in [Str74, FSt74, FSt76].
First, let $x$ be a letter of a permutation   
$w=x_1x_2\cdots x_n$ of a set of positive integers $Y=\{ y_1<y_2<\cdots<y_n\}$. The {\it $x$-factorization} of~$w$ is defined to be the sequence $(w_1,w_2,x,w_4,w_5)$, where 

(1) the juxtaposition product
$w_1w_2xw_4w_5$ is equal to~$w$;

(2) $w_2$ is the longest right factor of $x_1x_2\cdots x_{i-1}$, all letters of which are greater than~$x$;

(3) $w_4$ is the longest left factor of $x_{i+1}x_{i+2}\cdots x_n$, all letters of which are greater than~$x$.

\noindent
Next, say that~$x$ is {\it of type~I} (resp. {\it of type~II\/}) in~$w$, if whenever the juxtaposition product  $w_2w_{4}$ is non-empty, its maximum (resp. minimum) letter belongs to~$w_{4}$. Also, say that $x$ is of type~I  {\it and}~II, if $w_2$ and $w_4$ are both empty.

\medskip
{\it Definition} 2.1.\quad 
A permutation $w=x_{1}x_{2}\cdots x_{n}$ of  $Y=\{ y_1<y_2<\cdots<y_n\}$ is said to be an {\it Andr\'e permutation of the first kind} (resp. {\it of the second kind\/}) [in short, an {\it Andr\'e~I} (resp. an {\it Andr\'e~II\/})], if $x_{i}$ is of type~I (resp. of type~II) in~$w$ for every $i=1,2,\ldots,n$.

\medskip
{\it Definition} 2.2.\quad  
A permutation $w=x_1x_2\cdots x_n$ of  $Y=\{ y_1<y_2<\cdots<y_n\}$ is said to be an {\it Andr\'e permutation of the first kind} (resp. {\it of the second kind\/}), if it has no double descent (factors of the form $x_{i-1}>x_i>x_{i+1}$) and its troughs (factors $x_{i-1} x_i  x_{i+1}$  satisfying  $x_{i-1}> x_i$ and $x_i<x_{i+1}$) are all of type~I (resp. of type~II). 
By convention, $x_{n+1}:=0$.

\medskip
The following notations are being used. If $Y=\{y_{1}
<y_{2}<\cdots<y_{n}\}$ is a finite set of positive integers, let $\rho_{Y}$ be the increasing bijection of~$Y$ onto
$\{1,2,\ldots,n\}$. 
The inverse bijection of $\rho_Y$ is denoted by $\rho_Y^{-1}$. 
If $v=y_{i_1}y_{i_2}\cdots y_{i_n}$ is a permutation of~$Y$, written as a word, let $\rho_Y(v):=\rho_Y(y_{i_1})\rho_Y(y_{i_2})\cdots
\rho_Y(y_{i_n})=i_1i_2\cdots i_n$ be the {\it reduction} of the word~$v$, which is then a permutation of $1\,2\cdots n$.
When dealing with a given  word $v$, the subscript $Y$ in $\rho_Y(v)$ may be omitted, so that $\rho(v)=\rho_Y(v)$.
In the same way, the subscript $Y$ in each composition product $\rho_Y^{-1}\alpha\rho_Y(v)$ may be omitted, so that $\rho^{-1}\alpha\rho(v)=\rho^{-1}_Y\alpha\rho_Y(v)$.

Also, let
${\bf c}$ be the bijection $i\mapsto n+1-i$ of $\{1,2,\ldots,n\}$ onto itself. Furthermore, if
$v=y_{i_{1}}y_{i_{2}}\cdots y_{i_{n}}$ is a permutation of~$Y$, written as a word, let $C(v):=Y$ and the {\it length} of~$v$ be $|v|=n$. Finally, a {\it left maximum record} (resp. {\it left minimum record}\/) of~$v$ is defined to be a letter of~$v$ greater (resp. smaller) than all the letters to its left.

\medskip

\proclaim Proposition 2.1. Let $n\ge 2$ and $w=x_{1}x_{2}\cdots x_{n}$ be  Andr\'e I.\hfil\break
(1) In $w=v\min(w)v'$ both factors $v$ and $v'$ are  Andr\'e~I.\hfil\break
(2) If $w$ is from $\And^I_n$, then both permutations $1\,(x_{1}+1)\,(x_{2}+1)\,\cdots \,(x_{n}+1)$ 
and $x_{1}x_{2}\cdots x_{n}\,(n+1)$ belong to $\And_{n+1}^{I}$, and $(x_2-1)(x_3-1)\cdots (x_n-1)$ belongs to $\And_{n-1}^I$ whenever $x_1=1$.\hfil\break
(3) The last letter~$x_{n}$ is the maximum letter.\hfil\break
(4) Let $w=w' \,y\, w'' x_n$ with $y$ being the second greatest letter of~$w$. If $w''\not=e$, then $\bfF w''=\min(w'')$.\hfil\break
(5) For each left maximum record~$y$ of~$v$, less than $\max(w)$, the two factors $uy$ and $u'$ in the factorization $w=u\,y\,u'$ are themselves Andr\'e~I.\hfil\break
(6) Let $w=v\,y\,v'$ be Andr\'e I. If $y$ is a left minimum record, then $v$ is Andr\'e~I.\hfil\break
\vfill\eject

\noindent
(7) {\sl Let $w=v\,y\,v'$ be an arbitrary permutation with $y$ a letter.
If both factors $v$ and~$y\,v'$  are Andr\'e~I and~$y$ is a left minimum record, then~$w$ is Andr\'e~I.}

\proof 
(1) By the very definition given in Subsection~1.2.

(2) Clear.

(3) Write $w=v\min(w)v'$. By definition, $\max(v')=\max(vv')=\max(w)$ and by induction the last letter of~$v'$, which is also the last letter of~$w$, is equal to $\max(v')=\max(w)$.

(4) 
If $w''\not=e$, let $x:=\bfF w''$ and let $(w_1,w_2,x, w_4,w_5)$ be the $x$-factorization of~$w$. As $y$ is the maximum letter of~$w_2$, the maximum letter of~$w_4$ must be equal to $\max w$ to make~$x$ of type~I. This can be achieved only if $x$ is the minimum of~$w''$.

(5) Let $x:=\min (w)$ and $y$ be a left maximum record less than $\max(w)$, so that
$w=v\,x\,v'=u\,y\,u'$ for some factors $v$, $v'\not=e$, $u$, $u'\not=e$. Two cases are to be considered: (i) $x$ to the left of~$y$ so that $w=v\,x\,v''\,y\,u'$; (ii) $y$ to the left of~$x$ so that
$w=u\,y\,u''\,x\,v'$ for some factors $v''$, $u''$. In case (i) both factors~$v$ and $v''\,y\,u'$ are Andr\'e~I, following the definition of \S\thinspace 2.1. Now, the letter~$y$ is also a left maximum record of the word $v''\,y\,u'$. By induction on the length, both $v''\,y$ and~$u'$ are Andr\'e~I, so that the two factors  $v$ and $v''y$ of the word $u\,y=v\,x\,v''y$ are Andr\'e~I, making the latter word also Andr\'e~I. Thus, both~$uy$ and~$u'$ are Andr\'e~I. In case~(ii) the same argument applies: both factors $u\,y\,u''$ and~$v'$ are Andr\'e~I, then also $u\,y$ and~$u''$ by induction, as well as the juxtaposition product~$u''\, x\,v'$.\qed

(6) If $y=\min w$, then $v$ is Andr\'e I by definition. Otherwise, $y$ is to the left of $\min (w)$ in $w: w=v\,y\,u\,\min(w)\,u'$. But~$y$ is also a left minimum record of $v\,y\,u$. By induction on the length $v$ is Andr\'e~I.

(7) If $y=\min (w)$, then $v'$ is Andr\'e I by (2). Now, $v$ and $v'$ being both Andr\'e~I, the product  $w=v\,y\,v'$ is Andr\'e~I by definition. If $y>\min(w)$, then
$w=v\,y\,u\,\min(w)\,u'$. Nothing to prove if $v=e$. Otherwise, as $y\,v'$ is Andr\'e~I, both factors $y\,u$ and $u'$ are Andr\'e~I. As $y$ is also a left minimum record of $v\,y\,u$, the juxtaposition product $v\,y\,u$ is Andr\'e~I by induction on the length. Finally, $w$ itself is Andr\'e~I by definition, as $u'$ has been proved to be also Andr\'e~I.\qed

\medskip
In [FSch71] a bijection between $\And_{n}^I$ and ${\rm Alt}_{n}$ was constructed, but did not preserve the first letter. For proving Theorem 1.1(i) 
we need construct a bijection 
$$
\eta: w\rightarrow \eta(w),\quad\hbox{such that\quad}
{\bf F}\,w={\bf F}\,\eta(w)\leqno(2.1)
$$
of $\And_{n}^I$ onto the set ${\rm Alt}_{n}$ of all alternating permutations of length~$n$. For $n=1,2,3$ it suffices to take: $1\mapsto 1$, $12\mapsto 12$, $123\mapsto 132$, $213\mapsto 231$. When $n\ge 4$, each $w$ from $\And_{n}^{I}$ can be written $w=w'\,1\,w''$, 
where both factors $w'$, $w''$ (with $w'$ possibly empty) are Andr\'e I. 

If $w'=e$, let $v':=1$ and
$v'':=\rho^{-1}\,{\bf c}\,\eta\,
\rho(w'')$;

if $w'\not=e$, let 
\vskip-15pt
$$\leqalignno{\noalign{\vskip-12pt}
\kern2cm v'&:=\rho^{-1}\,\eta\,\rho(w');\cr
v''&:=\cases{\rho^{-1}\,\eta\,\rho(1\,w''),&if $|w'|$ even;\cr
\noalign{\smallskip}
\rho^{-1}\,{\bf c}\,\eta\,\rho(1\,w''),&if $|w'|$ odd;\cr}\cr
\noalign{\vskip-6pt}
\noalign{\hbox{and}}
\eta(w)&:=v'\,v''.&(2.2)\cr
}
$$

\goodbreak
\medskip
For instance, let $w=1234\in \And_{4}^{I}$. Then, $w'=e$, $w''=234$; then, $v'=1$,
$\rho(w'')=123$, $\eta(123)=132$, ${\bfc}(132)=312$,
$\rho^{-1}(312)=423=v''$ and
$\eta(1234)=1423$.

With $w=4361257$ we get: $w'=436$, $w''=257$; then,
$\rho(436)=213$, $\eta(213)=231$,
$\rho^{-1}(231)=463=v'$. Also,
$\rho(1257)=1234$,
$\eta(1324)=1423$, ${\bf c}(1423)=4132$,
$\rho^{-1}(4132)=7152=v''$ and
$\eta(4361257)=4637152$.

\proclaim Theorem 2.2. The mapping $\eta$ defined by $(2.2)$ is a bijection of $\And_{n}^{I}$ onto ${\rm Alt}_{n}$ such that ${\bf F}\,w={\bf F}\,\eta(w)$.

\proof
Again, factorize an Andr\'e I permutation $w$ in the form $w=x_1x_2\cdots x_n=w'1w''$.
When $w'=e$, then $\rho(w'')$
is an Andr\'e~I permutation starting with $\rho(x_{2})$.
By induction, $\eta\,
\rho(w'')$ is an {\it increasing} alternating permutation 
if $|w''|\ge 2$. Then, ${\bf c}\,\eta\,
\rho(w'')$ will be a {\it falling} alternating permutation, as well as the permutation
$v''=\rho^{-1}\,{\bf c}\,\eta\,
\rho(w'')$, which is also a permutation of $23\cdots n$. Hence, $\eta(v)=1\,v''$ will be an alternating permutation starting with~1.

When $w'\not =e$, then $v'$ is an alternating permutation of the set $C(w')$. By induction, it starts with the same letter as the first letter of~$w$, that is, $x_{1}$. If $|w'|$ is even,
$v''=\rho^{-1}\,\eta\,\rho(1\,w'')$ is an alternating permutation starting with $1$, by induction. The juxtaposition product $v'v''$ will then be an alternating permutation starting with~$x_{1}$, as the last letter of~$v'$ is necessarily greater than the first letter of~$v''$. If $|w'|$ is odd, we just have to verify that ${\bf L}\,
v'<{\bf F}\,v''$. But $w$, being an Andr\'e~I permutation, ends with its maximum letter~$n$ and so does~$w''$. 
By induction, $\eta\,\rho(1\,w'')$ starts with~1, so that ${\bf c}\,\eta\,\rho(1\,w'')$ starts with the maximum letter~$n$. Therefore, $v''$ is a falling alternating permutation starting with~$n$ and $v'v''$ is alternating permutation starting with~$x_{1}$.\qed

\goodbreak
\medskip
For each permutation $w=x_{1}x_{2}\cdots x_{n-1}x_{n}$ $(n\ge 2)$ the {\it next to the last letter} ${\bf N\!L}\,w$ of~$w$ has been defined as ${\bf N\!L}\,w:=x_{n-1}$. The construction of a bijection $\theta$ of~$\And_{n}^I$ onto itself having the property
$$\displaylines{\rlap{(2.3)}\hfill
{\bf N\!L}\,\theta(w)+{\bf F}\,w=|w|=n\hfill\cr
\noalign{\hbox{is quite simple. It suffices to define:}}
\rlap{(2.4)}\hfill
\theta(x_{1}\cdots x_{n-2} x_{n-1}\,n)
:=(n-x_{n-1})\,(n-x_{n-2})\,\cdots \,(n-x_{1})\,n.
\hfill\cr}
$$

\goodbreak\noindent
Property (2.3) is readily seen. It remains to prove that if $w$ belongs to $\And_{n}^I$, so does $\theta(w)$. This is the object of the next Proposition.

\medskip
\proclaim Proposition 2.3.
Let $Y=\{y_1<y_2<\cdots <y_n\}$ be a finite set of positive integers and $w=x_1x_2\cdots x_n$ be an Andr\'e~I permutation from the set
$\And_Y^I$. Then $\theta(w):=(x_n-x_{n-1})(x_n-x_{n-2})\cdots (x_n-x_1)x_n$ is also Andr\'e~I.

\medskip
{\it Proof}.\quad Proposition 2.3 is true for $n=2$, as $\theta(y_1y_2)=(y_2-y_1)y_2$. For $n=3$ we have $\theta(y_1y_2y_3)=(y_3-y_2)(y_3-y_1)y_3$, $\theta(y_2y_1y_3)=(y_3-y_1)(y_3-y_2)y_3$, which are two Andr\'e~I permutations. 

For $n\ge 4$ let $w=x_{1}x_{2}\ldots x_{n}\in \And_{Y}^I$ be written $w=w'y_1w''$. If $w'=e$, let $w''-y_1:=(x_{2}-y_1)\cdots (x_{n-2}-y_1)(x_{n-1}-y_1)(y_n-y_1)$. Then,  $w''$ is Andr\'e~I by Lemma 2.1~(b), as well as $w''-y_1$, since $y_1$ is the smallest element of~$Y$. By induction,
$\theta(w''-y_1)=(y_n-y_1-(x_{n-1}-y_1))\,(y_n-y_1-(x_{n-2}-y_1))\,\cdots\,
(y_n-y_1-(x_{2}-y_1))(y_n-y_1)=(y_n-x_{n-1})\,(y_n-x_{n-2})\,\cdots\, (y_n-x_{2})(y_n-y_1)$ is Andr\'e~I. Therefore, $(y_n-x_{n-1})\,(y_n-x_{n-2})\,\cdots\, (y_n-x_{2})(y_n-y_1)y_n
=(x_n-x_{n-1})\,(x_n-x_{n-2})\,\cdots\, (x_n-x_{2})(x_n-x_1)y_n$ is Andr\'e I a fortiori and is precisely the expression of $\theta(w)$ that was wanted.

\goodbreak
Let $|w'|=k\ge 1$. By induction, both $\theta(w'y_n)
=(y_n-x_{k})\cdots (y_n-x_{2})(y_n-x_{1})\,y_n$ and $\theta(y_1w'')
=(y_n-x_{n-1})\cdots (y_n-x_{k+2})\,(y_n-y_1)\,y_n$ are Andr\'e~I, and also $\check\theta(y_1w''):=(y_n-x_{n-1})\cdots (y_n-x_{k+2})\,(y_n-y_1)$. The juxtaposition product 
$\check\theta(y_1w'')\,\theta_n(w'y_n)$ reads:
$(y_n-x_{n-1})\cdots (y_n-x_{k+2})\,(y_n-y_1)
(y_n-x_{k})\cdots (y_n-x_{2})(y_n-x_{1})\,y_n$, that is, precisely 
$\theta(w)$, since $y_n-y_1=y_n-x_{k+1}$. 

Now, note that $x_k$ is the greatest letter of~$w'$ by Lemma~2.1~((3), so that $(y_n-x_k)$ is the smallest letter of the right factor $(y_n-x_{k})\cdots (y_n-x_{2})(y_n-x_{1})\,y_n$ of $\theta(w)$. On the other hand, $(y_n-y_1)>(y_n-x_k)$. Thus, $(y_n-x_k)$ is  a trough of $\theta(w)$; moreover, the $(y_n-x_k)$-factorization $(w_1,w_2,(y_n-x_k),w_4,w_5)$ of $\theta(w)$ is of type~I, since 
$w_2$ contains the letter $(y_n-y_1)$ and~$w_4$ the letter $y_n$, which is greater than
$(y_n-y_1)$. Finally, the $x$-factorizations of the other letters $x$ from 
$\check\theta(y_1w'')$ (resp. from $\theta(w'y_n)$) in each of those two factors are identical with their
$x$-factorizations in~$\theta(w)$. They are then all of type~I, and $\theta(w)$ is Andr\'e~I. \qed

\medskip
\goodbreak
By Proposition 2.3 and Identity (2.3) we have:

\proclaim Theorem 2.4. The statistics ``$\bf F$'' and ``$(n-{\bf N\!L})$'' are both Entringerian on $\And_{n}^I$. Moreover, the distribution of the bivariate statistic
$({\bf F},n-{\bf N\!L})$ on $\And_{n}^I$ is symmetric.

\vfill\eject

\centerline{\bf 3. The bijection $\phi$ between Andr\'e I and Andr\'e II permutations}

\medskip
For each permutation $w=x_1x_2\cdots x_n$ of $12\cdots n$ $(n\ge 2)$ make the convention $x_{n+1}:=0$ and
introduce the statistic {\it spike of $w$}, denoted by ``$\spi w$,"  to be equal to the letter~$x_i$ $(1\le i\le n)$ having the properties:
$$
x_1\le x_1,\ x_1\le x_2,\ \ldots \ ,\ x_1\le x_i,\quad{\rm and}\quad x_1>x_{i+1}.\leqno(3.1)
$$ 
The spike statistic may be regarded as the permutation-analog of the classical one that measures the time spent by a particle starting at the origin and wandering in  the $y>0$ part of the $xy$-plane, before crossing the $x$-axis for the first time. For instance, $\spi (253416) =4$, as all the letters to the left of~4 are greater than or equal to~2, but the letter following~4 is less than~2. Also, $\spi(425136)=4$ and $\spi(14235)=5$.

When $w$ is an Andr\'e I permutation and $\spi w=x_i$, then $x_i$ is a left maximum record, i.e., greater than all the letters to its left. Otherwise, the minimum trough between the maximum letter within $x_1x_2\cdots x_{i-1}$ would not be of type~I. Accordingly, when $w$ is an Andr\'e~I permutation, the spike~$x_i$ of~$w$ can also be defined as the {\it smallest left maximum record}
(or the {\it leftmost one}), whose successive letter~$x_{i+1}$ is less than~$x_1$.

\medskip
For introducing the statistic ``$\pit$'' we restrict the definition to all permutations $w=x_1x_2\cdots x_n$ of $12\cdots n$ such that $n\ge 2$ and $x_{n-1}<x_n$. Let $1=a_1<a_2<\cdots<x_{n-1}=a_{k-1}<x_n=a_k$ be the increasing sequence of the right minimum records of~$w$, that is to say, the letters which are  smaller than all the letters to their right. With the assumption $x_{n-1}<x_n$, 
there are always two right minimum records to the right
of each letter greater than~$x_n$. 
If $x_n=n(=\max w)$, let $\pit w:=1(=\min w)$. Otherwise, let $x_i$ be the rightmost letter greater than~$x_n$ and $a_j<a_{j+1}$ be the closest pair of right minimum records to the right of~$x_i$. Define:  $\pit w:=a_{j+1}$.

For instance, $\pit(451236) =1$, as the word ends with the maximum letter~6. With the permutation
$614235$ the letter 6  is the rightmost  letter greater than $x_n=5$ and $1<2$ is the closest pair of right minimum records to the right of~6, so that $\pit(614235)=2$.

An alternate definition for ``$\pit$" is the following: 
if $w$ ends with $\max w$, let $\pit w=\min w$. Otherwise, write
$w=w_1(\min w)w_2$ and define: $\pit w:=\pit w_2$. If $w_2$ does not end with the maximum letter, let $w_2=w_3(\min w_2)w_4$ and define $\pit w_2:=\pit w_4$, continue the process until finding a right factor $w_{2j}$ ending with its maximum letter to obtain:
$\pit w:=\pit w_2=\cdots =\pit w_{2j}=\min w_{2j}$. For instance,
$\pit(614235)=\pit(4235)=2$.

\vfill\eject
Remember that we have introduced three other statistics ``$\NL$'' (``next to the last"), ``$\bfL$" (``last'') and ``$\grn$'' (``greater neighbor of the maximum") and that $\grn w=\NL w$ whenever~$w$ is an Andr\'e~I permutation. Our goal is to prove the next theorem.

\proclaim Theorem 3.1. The triplets $(\bfF, \spi,\NL)$ on $\And_n^I$ and
$(\pit,\bfL, \grn)$ on $\And_n^{I\!I}$ are equidistributed.

Let $X=\{a_1<a_2<\cdots<a_n\}$ be a set of positive integers (or any finite totally ordered set) and $\And_X^I$ (resp. $\And_X^{I\!I}$) denote the set of all Andr\'e~I (resp. Andr\'e II) permutations of~$X$.
To prove the previous statement a bijection 
$$\leqalignno{
\phi:\And_X^I&\rightarrow \And_X^{I\!I}&(3.2)\cr
\noalign{\hbox{ will be constructed having the property:}}
(\bfF, \spi,\NL)\,w&=(\pit,\bfL,\grn)\,\phi(w).&(3.3)}$$

 When $n=0$ let $\phi(e):=e$ with $e$ the empty word. Let $\phi(a_1):=a_1$
for $n=1$; $\phi(a_1a_2):=a_1a_2$ for $n=2$. For $n\ge 3$ each permutation~$w$ from
$\And_X^I$ has one of the two forms:
$$
(i)\ w=v_0\,a_1\,v_1\,a_2\,v_2;\qquad
(ii)\ w=v_0\,a_2\,v_2\,a_1\,v_1.
$$
Note that the three factors $v_0$, $v_1$, $v_2$ and the product $v_0a_2v_2$ are all Andr\'e~I permutations and~$v_0$ is possibly empty. In case $(i)$ $v_1$ may be empty, but not $v_2$ (which ends with $a_n$ greater than $a_2$); in case $(ii)$ $v_2$ may be empty, but not~$v_1$ (which ends with $a_n$). For both cases $(i)$ and $(ii)$ define:
$$
\phi(w):=\phi(v_1)\,a_1\,\phi(v_0\,a_2\,v_2).\leqno(3.4)
$$
By induction both factors $\phi(v_1)$ and $\phi(v_0a_2v_2)$ are Andr\'e II, 
as well as $\phi(w)$, since $a_2$ is to the right of $a_1$.
\medskip

\goodbreak
\def\trait{\vrule height 10pt depth 3pt width .4pt}

{\it Example}.\quad Consider the Andr\'e I permutation
$$
\leqalignno{
w&=7\ 8\ \vtop{\hbox{5}\hbox{$v_0$}} 6\ 9\  \trait\  \vtop{\hbox{2}
\hbox{$a_2$}} \, \trait\  \vtop{\hbox{10 }\hbox{$v_2$}}\, \trait\ \vtop{\hbox{1}\hbox{$a_1$}} \trait\ 11\ 3\ \vtop{\hbox{12}\hbox{$v_1$}}\ 4\ 13;\cr
\noalign{\hbox{we successively have:}}
\phi(w)&\buildrel (ii)\over =
\phi(11\ 3\ 12\ 4\ 13)\, 1\, \phi(7\ 8\ 5\ 6\ 9\ 2\ 10);\cr
\phi(\vtop{\hbox{11}\hbox{$v_0$}}\ \vtop{\hbox{3}\hbox{$a_1$}}\ 
\vtop{\hbox{12}\hbox{$v_1$}}\ \vtop{\hbox{4}\hbox{$a_2$}}\ \vtop{\hbox{13}\hbox{$v_2$}})&\buildrel (i)\over =
12\ 3\ \phi(11\ 4\ 13);\cr
\phi(\vtop{\hbox{11}\hbox{$a_2$}}\ \vtop{\hbox{4}\hbox{$a_1$}}\ 
\vtop{\hbox{13}\hbox{$v_1$}})
&\buildrel (ii)\over =13\ 4\ 11;\cr
\phi(7\  \vtop{\hbox{8}\hbox{$v_0$}} \  \trait\  \vtop{\hbox{5}
\hbox{$a_2$}} \, \trait\  6\ \vtop{\hbox{9}\hbox{$v_2$}}\, \trait\ \vtop{\hbox{2}\hbox{$a_1$}} \trait\  \vtop{\hbox{10}\hbox{$v_1$}})
&\buildrel (ii)\over =10\ 2\ \phi(7\ 8\ 5\ 6\ 9);\cr
\phi(7\  \vtop{\hbox{8}\hbox{$v_0$}} \  \trait\  \vtop{\hbox{5}
\hbox{$a_1$}} \, \trait\  \vtop{\hbox{6}\hbox{$a_2$}}\trait  \ \vtop{\hbox{9}\hbox{$v_2$}})
&\buildrel (i)\over =
5\ \phi(7\ 9\ 6\ 9);\cr
\phi(\vtop{\hbox{7}\hbox{$a_2$}} \  \trait\  \vtop{\hbox{8}
\hbox{$v_2$}} \, \trait\  \vtop{\hbox{6}\hbox{$a_1$}}\trait  \ \vtop{\hbox{9}\hbox{$v_1$}})
&\buildrel (ii)\over =
9\ 6\ \phi(7\ 8)=9\ 6\ 7\ 8;\cr
\noalign{\hbox{so that}}
\phi(w)&= 12\ 3\ 13\ 4\ 11\ 1\ 10\ 2\ 5\ 9\ 6\ 7\ 8.\cr}
$$
We can verify that\qquad
$(\bfF,\spi,\NL)\,w=(\pit,\bfL,\grn)\,\phi(w)=(7,8,4)$.

\goodbreak
\medskip

With the previous definition of~$\phi$ we see that the maximum letter~$a_n$ of~$X$ occurs in $\phi(v_0a_2v_2)$ (resp. in $\phi(v_1$)) when $w$ is of form $(i)$ (resp. of form~$(ii)$). For constructing the inverse~$\phi^{-1}$ of~$\phi$ this suggests that we start with the factorization $$v=w_0\,a_1\,w_1$$ of each permutation~$v$ from $\And_X^{I\!I}$ with $\#X\ge 3$, after defining: 
$\phi^{-1}(e):=e$; $\phi^{-1}(a_1):=a_1$
for $n=1$; $\phi^{-1}(a_1a_2):=a_1a_2$ for $n=2$. As both $w_0$, $w_1$ are Andr\'e II permutations with fewer letters, the images $\phi^{-1}(w_0)$, $\phi^{-1}(w_1)$ are defined by induction. Let $v_1:=\phi^{-1}(w_0)$. As the minimum letter of~$w_1$ is~$a_2$, define $v_0$ and $v_2$ to be the factors in $\phi^{-1}(w_1):=v_0\,a_2\,v_2$.
Next, let
$$
\phi^{-1}(v):=\cases{v_0\,a_1\,v_1\,a_2\,v_2,&if $a_n$ is a letter of $w_1$;\cr
v_0\,a_2\,v_2\,a_1\,v_1,&if $a_n$ is a letter of $w_0$.\cr}\leqno(3.5)
$$

\proclaim Lemma 3.2.
The Andr\'e II permutation $\phi(w)$ ends with its maximum letter, if and only if~$w$ starts with its minimum letter $\min w$, and then
$\pit\phi(w)=\bfF w=\min w$. 

\proof
This is obviously true for $n=2$. For $n\ge 3$ we have: $w=a_1v_1a_2v_2$ and $\phi(w)=\phi(v_1)a_1\phi(a_2v_2)$. As $a_2v_2$ is Andr\'e~I starting with its minimum letter $a_2$, then, by induction,  $\phi(a_2v_2)$ ends with its maximum letter, which is equal to $\max w$. Hence, $\pit \phi(w)=\pit(\phi(v_1)a_1\phi(a_2v_2))=a_1=\bfF w$. For the converse take the notation 
$v=w_0a_1w_1$ of (3.5). When the maximum letter~$a_n$ occurs in~$w_1$, then
$\psi(v)= v_0a_1v_1a_2v_2$ with $v_1=\psi(w_0)$ and $\psi(w_1)=v_0a_2v_2$. By assumption, $a_n$ occurs at the end of~$v$, therefore, at the end of~$w_1$. By induction,
$\psi(w_1)=v_0a_2v_2$ starts with its minimum letter. This can be true only if $v_0=e$. Therefore, $\psi(v)=a_1v_1a_2v_2$ and starts with its minimum letter~$a_1$.\qed

\proclaim Theorem 3.3. The mapping $\phi$  is a bijection of $\And_n^I$ onto $\And_n^{I\!I}$. Moreover, relation~$(3.3)$ holds.

\proof
The bijectivity is proved by the construction of the inverse $\phi^{-1}$ (see (3.5)).
To prove identity (3.3), let $w$ be a Andr\'e I permutation, either of the form $v_0a_1v_1a_2v_2$, or of the form $v_0a_2v_2a_1v_1$. In both cases 
$$
\eqalignno{\spi w&=\spi v_0=\spi(v_0a_2v_2)\cr
&=\bfL\phi(v_0a_2v_2)&\hbox{[by induction]}\cr
&=\bfL(\phi(v_1)a_1\phi(v_0a_2v_2))=\bfL\phi(w).\cr}
$$
Next, if $w=v_0a_1v_1a_2v_2$, then
$$\eqalignno{
\NL\, w&=\NL(v_0a_2v_2)\cr
&=\grn\phi(v_0a_2v_2)&\hbox{[by induction]}\cr
&=\grn( \phi(v_1)a_1\phi(v_0a_2v_2))=\grn\phi(w),\cr}
$$
because the maximum letter $a_n$ is a letter of~$v_2$.
If $w=v_0a_2v_2a_1v_1$, then $a_1v_1$ has at least two letters
 and ends with~$a_n$, so that $a_1v_1$ is Andr\'e I of form~$(i)$.
 Consequently,
$$\eqalignno{
\NL\, w&=\NL(a_1v_1)\cr
&=\grn\phi(a_1v_1)&\hbox{[by induction]}\cr
&=\grn( \phi(v_1)a_1\phi(v_0a_2v_2))=\grn\phi(w).\cr}
$$

When  $w$ does not start with its minimum letter, then 
$\phi(w)=\phi(v_1)a_1\phi(v_0a_2v_2)$ and $\phi(v_0a_2v_2)$ does not end with $\max w$. Therefore, $\pit \phi(w)=\pit \phi(v_0a_2v_2)=\bfF v_0a_2v_2$.
As $v_0=e$ only in case~(ii), we then have:
$\pit \phi(w)=\bfF v_0a_2=\bfF w$.\qed

\bigskip
\centerline{\bf 4. The bijection $g$ of the set of Andr\'e I permutations onto itself}

\medskip
When making up the tables of the distribution of the bivariate $(\spi,\bfF)$ on $\And_n^I$
for $n=1,2,\ldots,7$, as shown in Table 4.1, it can be noticed that the matrices are symmetric with respect to their skew-diagonals. The property will hold in general if a bijection
$g$ of $\And_n^I$ onto itself can be constructed satisfying
$$
(\bfF,\spi)\,w=(n+1-\spi,n+1-\bfF)\,g(w)\leqno(4.1)
$$
for all $w$ from $\And_n^I$.
$$
\vbox{\offinterlineskip\halign{
\hfil$#$\hfil\ \vrule&&\strut\ \hfil$#$\hfil\ \cr 
\bfF  =&    1 & \cr
\noalign{\hrule}
\spi =   2 &   1 & \cr
\noalign{\medskip}
\multispan3\hfil $n=2$ \hfil\cr
}}\quad
\vbox{\offinterlineskip\halign{
\hfil$#$\  \vrule&&\strut\ \hfil$#$\hfil\ \cr 
\bfF = &    1 &   2 & \cr
\noalign{\hrule}
\spi =   2 &   . &   1 & \cr
        3 &   1 &   . & \cr
        \noalign{\medskip}
\multispan4\hfil $n=3$ \hfil\cr}}
\quad
\vbox{\offinterlineskip\halign{
\hfil$#$\ \vrule&&\strut\ \hfil$#$\hfil\ \cr 
\bfF  =&    1 &   2 &   3 & \cr
\noalign{\hrule}
\spi =   2 &   . &   1 &   . & \cr
        3 &   . &   1 &   1 & \cr
        4 &   2 &   . &   . & \cr
        \noalign{\medskip}
\multispan5\hfil $n=4$ \hfil\cr
}}\quad
\vbox{\offinterlineskip\halign{
\hfil$#$\ \vrule&&\strut\ \hfil$#$\hfil\ \cr 
\bfF  =&    1 &   2 &   3 &   4 & \cr
\noalign{\hrule}
\spi =   2 &   . &   2 &   . &   . & \cr
        3 &   . &   1 &   3 &   . & \cr
        4 &   . &   2 &   1 &   2 & \cr
        5 &   5 &   . &   . &   . & \cr
                \noalign{\medskip}
\multispan6\hfil $n=5$ \hfil\cr
}}$$
$$
\vbox{\offinterlineskip\halign{
\hfil$#$\ \vrule&&\strut\ \hfil$#$\hfil\ \cr 
\bfF=  &    1 &   2 &   3 &   4 &   5 & \cr
\noalign{\hrule}
\spi =   2 &   . &   5 &   . &   . &   . & \cr
        3 &   . &   2 &   8 &   . &   . & \cr
        4 &   . &   3 &   3 &   8 &   . & \cr
        5 &   . &   6 &   3 &   2 &   5 & \cr
        6 &  16 &   . &   . &   . &   . & \cr
                       \noalign{\medskip}
\multispan7\hfil $n=6$ \hfil\cr
}}\quad
\vbox{\offinterlineskip\halign{
\hfil$#$\ \vrule&&\strut\ \hfil$#$\hfil\ \cr 
\bfF  =&    1 &   2 &   3 &   4 &   5 &   6 & \cr
\noalign{\hrule}
\spi =   2 &   . &  16 &   . &   . &   . &   . & \cr
        3 &   . &   5 &  27 &   . &   . &   . & \cr
        4 &   . &   7 &   8 &  31 &   . &   . & \cr
        5 &   . &  11 &  10 &   8 &  27 &   . & \cr
        6 &   . &  22 &  11 &   7 &   5 &  16 & \cr
        7 &  61 &   . &   . &   . &   . &   . & \cr
                              \noalign{\medskip}
\multispan8\hfil $n=7$ \hfil\cr
}}$$
\centerline{Table 4.1: distribution of $(\spi,\bfF)$ on $\And_n^I$}

\medskip

For the construction of $g$ we proceed as follows.
An Andr\'e~I permutation $v=y_1y_2\cdots y_l$ on a set~$X$ (of cardinality $l\ge 2$) is called {\it simple}, if the first letter~$y_1$ of~$v$ is equal to $\min X$.
Consider an Andr\'e~I permutation $w=x_1x_2\cdots x_n$ from $\And_n^I$. Let $1=a_1<a_2<\cdots <a_r$ (resp. $1=b_1<b_2<\cdots<b_s=n$) be the increasing sequence of subscripts such that
$x_{a_1}>x_{a_2}>\cdots >x_{a_r}$
(resp. $x_{b_1}<x_{b_2}<\cdots <x_{b_s}$)
is the increasing (resp. decreasing) sequence of
the left minimum (resp. maximum) records of  $w=x_1x_2\cdots x_n$ from $\And_n^I$. 

\goodbreak
For the following Andr\'e I permutation the left minimum (resp. maximum) records are underlined (resp. overlined):
$$w=\underline{\strut\overline{\strut7 }}\;\overline{\strut8} \; \underline{\strut5} \; 6 \;\overline{\strut 9}\; \underline{\strut2} \;\overline{\strut 10}\; \underline{\strut1}\; \overline{\strut 11} \;3 \;\overline{\strut 12} \;4 \;\overline{\strut 13}$$

Going back to the general case let $v_1:=x_1\cdots x_{a_2-1}$,
$v_2:=x_{a_2}\cdots x_{a_3-1}$, \dots~, $v_r:=x_{a_{r}}\cdots x_{n}$, so that $w$ is the juxtaposition product $v_1v_2\cdots v_r$ and the factors $v_i$ are obtained by cutting the word~$w$ just before each left minimum record. The factorization
$(v_1,v_2,\ldots, v_r)$ is called the {\it canonical factorization} of the Andr\'e~I permutation~$w$. Furthermore, the sequence
$$(\,(\bfF v_1,\bfL v_1),(\bfF v_2,\bfL v_2),\ldots,(\bfF v_r,\bfL v_r)\,),$$
also equal to 
$(\,(x_1,x_{a_2-1}), (x_{a_2},x_{a_3-1}),\ldots,(x_{a_r,}x_n)\,)$, is called the 
{\it type} of the canonical factorization of~$w$.

With the running example the canonical factorization reads:
$$\leqalignno{
w&=\underline{\strut\overline{\strut7 }}\;\overline{\strut8} \mid \underline{\strut5} \; 6 \;\overline{\strut 9}\mid\underline{\strut2} \;\overline{\strut 10}\mid\underline{\strut1}\; \overline{\strut 11} \;3 \;\overline{\strut 12} \;4 \;\overline{\strut 13}\cr
&\hskip15pt  v_1\hskip17pt v_2\qquad v_3\qquad\quad v_4\cr}
$$
and is of type $(\,(7,8),(5,9),(2,10),(1,13)\,)$.

\goodbreak

\proclaim Proposition 4.1. Let $(v_1, v_2,\ldots, v_r)$ be the {\it canonical factorization} of the Andr\'e~I permutation~$w=x_1x_2\cdots x_n$ from $\And_n^I$.
Let~$s$ be the number of left maximum records of~$w$. Then,\hfill\break
{\rm (i)} $r\le s$;\hfill\break
{\rm (ii)} each factor $v_i$ $(i=1,2,\ldots,r)$ is a simple Andr\'e I permutation;\hfill\break
{\rm (iii)} $\bfL v_i$ is a left maximum record, so that $\bfL v_1<\bfL v_2<\cdots <\bfL v_{r-1}<\bfL v_r=n$ and, of course,
$\bfF v_1>\bfF v_2>\cdots >\bfF v_{r-1}>\bfF v_r=1$;\hfill\break
{\rm (iv)} $(\bfF v_1,\bfL v_1)=(\bfF  w,\spi w)$.

\medskip
Let $w=x_1x_2\cdots x_n$ be an Andr\'e I permutation from $\And_n^I$ and $\overline w:=\overline x_n\overline x_{n-1} \overline x_{n-2} \cdots \overline x_1 $ be the permutation defined by $\overline x_i:=N-x_{n+1-i}$ $(i=1,2,\ldots,n)$, where $N$ is some integer greater than~$n$.

\proclaim Proposition 4.2. If $w$ is a simple Andr\'e I permutation, so is $\overline w$.

\medskip
For constructing the bijection $g$ let $n\ge 3$ and $N:=n+1$. If $(v_1,v_2,\ldots,v_r)$ is the canonical factorization of a permutation $w$ from $\And_n^I$,
define $g(w)$ to be the juxtaposition product:
$$
g(w):=\overline v_1\,\overline v_2\,\ldots\,\overline v_r.\leqno(4.1)
$$
Furthermore, if $\tau=((p_1,q_1),(p_2,q_2),\ldots,(p_r,q_r))$ is the canonical factorization type of~$w$, let
$$
\overline\tau:=
((\overline q_1,\overline p_1),(\overline q_2,\overline p_2),\ldots,(\overline q_r,
\overline p_r)).
$$
We then have the fundamental property of~$g$.

\proclaim Theorem 4.3. 
The transformation $g$ is a bijection of $\And_n^I$ onto itself. Furthermore, if $\tau$
is the canonical factorization type of~$w$, then $\overline\tau$ is the canonical factorization type of $g(w)$. In particular,
$$(\bfF,\spi)\,g(w)=(n+1-\spi,\; n+1-\bfF)\,w.\leqno(4.2)
$$

\medskip
With the running example and $n+1=14$ we get:
$$g(w)=6\,7\,\mid\,5\,8\,9\,\mid\,4\,12\,\mid\,1\,10\; 2\,11\;3\,13$$
which is  of type $(\,(6,7),(5,9),(4,12),(1,13)\,)$.
\medskip

The proofs of Propositions 4.1, 4.2  and Theorem 4.3 do not present any difficulties and will be omitted. 

\vfill\eject
\centerline{\bf 5. The proof of Theorem 1.1 (iii) and (iv)}

\medskip
We reproduce the sequence (1.5) by decomposing the product $\phi\circ g$ :
$$
\matrice{
\And_n^I\quad &\buildrel  g \over \longrightarrow&\And_n^I 
&\buildrel  \phi \over \longrightarrow& \And_n^{I\!I}\cr
w&\mapsto &g(w)&\mapsto& \phi(g(w))\cr
\bfF w&=&n+1-\spi g(w)&=&n+1-\bfL\phi(g(w))\cr}
\leqno(5.1)
$$
The first (resp. second) identity $\bfF w=n+1-\spi g(w)$ (resp. $\spi g(w)=\bfL \phi (g(w))$)
is a specialization of (4.2) (resp. of (3.3)).

\medskip
Take the example of the previous section: $w=7\, 8\, 5\, 6\,9\, 2\,\twodigit10\, 1\, \twodigit11\, 3\, \twodigit12\,4\,\twodigit13$ and
$g(w)=6\,7\,5\,8\,9\,4\,\twodigit12\,1\,\twodigit10\, 2\,\twodigit11\;3\,\twodigit13$. By using the definition of~$\phi$ given in (3.4) we get:
$\phi(g(w))=\twodigit10\,1\,\twodigit11\,2\,\twodigit13\,3\,\twodigit12\,4\,8\,9\,5\,6\,7$ belonging to $\And_{13}^{I\!I}$ 
and $n+1-\bfL \phi(g(w))=14-7=7=\bfF w$.

Next, reproduce the sequence (1.4) by decomposing the product $\phi\circ \theta$:
$$
\matrice{
\And_n^I \quad &\buildrel  \theta \over \longrightarrow&\And_n^I 
&\buildrel  \phi \over \longrightarrow& \And_n^{I\!I}\cr
w&\mapsto &\theta(w)&\mapsto& \phi(\theta(w))\cr
\bfF w&=&n-\NL\, \theta(w) &=&n-\grn\phi(\theta(w))\cr}
\leqno(5.2)
$$
The first (resp. second) identity $\bfF w=n-\NL\, \theta(w)$ (resp. $\NL\, \theta(w)=\grn(\phi (g(w))$)
is a specialization of (2.3) (resp. of (3.3)). 

\medskip
For example, with
$w'=\twodigit10\,2\,\twodigit11\,3\,\twodigit12\,1\,9\,4\,5\,8\,6\,7\,\twodigit13$, we successively obtain:
$\theta(w')=6\,7\,5\,8\,9\,4\,\twodigit12\,1\,\twodigit10\,2\,\twodigit11\,3\,\twodigit13$
and
$\phi(\theta(w'))=\twodigit10\,1\,\twodigit11\,2\,\twodigit13\,3\,\twodigit12\,4\,8\,9\,5\,6\,7$.
Thus, ${\bf F}\,w'=10=n-{\rm N\!L}\,\theta(w')=13-3=n-{\bf grn}\,\phi(\theta(w'))$.

\medskip

The proofs of (iii) and (iv) of Theorem 1.1 are now completed. 
Another proof of Theorem~1.1 (iii) and (iv) makes use of the properties of a rearrangement group $G_n$, acting on the group ${\goth S}_n$ of all the permutations of $\{1,2,\ldots,n\}$, which were developed in [FSt74, FSt76] and another correspondence~$\Gamma$ on binary increasing trees, introduced in [FH13]. They constitute the main ingredients for the constructions of three bijections $\overline \Gamma$, $\Phi^I$ and $\Phi^{I\!I}$ appearing in the next diagram 
$$\displaylines{\noalign{\vskip-1pt}
\kern2.5cm
\matrix{\And_n^{I\!I}&\buildrel \textstyle\overline \Gamma\over \longrightarrow&
{\goth S}_n/G_n&\buildrel \textstyle\Phi^I\over \longrightarrow&\And_n^I\cr
&&\;\Big\downarrow\Phi^{I\!I}&\cr
&&\And_n^{I\!I}\cr}\hfill\cr
\noalign{\hbox{having the property:}}
\hfill\bfL w=1+\NL\, \Phi^I\overline\Gamma(w)
=1+\grn \Phi^{I\!I}\overline\Gamma(w).\hfill\cr}
$$

\vfill\eject

\centerline{\bf 6. Combinatorics of the twin Seidel matrix sequence}

\medskip
This Section is devoted to proving Theorem~1.2. As announced in Subsection~1.5, the question is to show that the integers $a_n(m,k)$ and $b_n(m,k)$, when taken as $a_n(m,k)=\#A_n(m,k)$, $b_n(m,k)=\#B_n(m,k)$ with
$$\leqalignno{
A_n(m,k)&:=\{w\in \And_n^I:( \bfF,\NL)w=(m,k)\};&(6.1)\cr
	B_n(m,k)&:=\{w\in \And_n^{I}:(\spi,\grn)w=(m,k)\};&(6.2)\cr
}$$
satisfy all the properties (TS1)--(TS5) 
stated in Subsection~1.5. 

The verifications of properties (TS1), (TS2), (TS3), (TS4.1), (TS4.2), (TS5.1) are easy and given in the next Subsection. The proofs of the other properties are much harder and will be developed thereafter.

\medskip
6.1. {\it The first evaluations}.\quad
By (1.7) the set $B_n(m,k)$ is equipotent with 
$$
B'_n(m,k):=\{w\in \AndII_n:(\bfL,\grn)w=(m,k)\}.\leqno{(6.3)}
$$
The evaluations in this subsection are made by using $B'_n(m,k)$ instead of $B_n(m,k)$.
\medskip
(TS1) Nothing to prove, except for the diagonals of the twin Seidel matrices $A_n$ and $B_n$. They have zero entries when $n\ge 3$, because the first and next to the last letter of each Andr\'e I permutation cannot be the same! On the other hand, the identity $\bfL w=\grn w=m$ would mean that the permutation~$w$ from $\And_n^{I\!I}$ ends with  a double descent $n>m>0$.

\medskip
(TS2) We have $a_n(k,n)=b_n(k,n)=0$, because $\grn w\le n-1$ for each~$w$ from either $\And_n^I$, or $\And_n^{I\!I}$. Also, $a_n(n,k)=0$, as each permutation from $\And_n^I$ ends with~$n$. Finally, $b_n(1,k)=0$, because each permutation from  $\And_n^{I\!I}$ cannot end with the letter~1.

\medskip
(TS3) We have: $A_2=\matrice{1&\cdot\cr \cdot&\cdot\cr}$ and
$B_1=\matrice{\cdot&\cdot\cr 1&\cdot\cr}$, because
$\And_2^I=\And_2^{I\!I}=\{12\}$ and $(\bfF,\NL,\bfL, \grn)(12)=(1,1,2,1)$.

\medskip
(TS4.1) The entry $b_n(n,k)$ counts the Andr\'e II permutations~$w$ from $\And_n^{I\!I}$ ending with the two-letter factor $k\,n$. The deletion of the ending letter~$n$ maps~$w$ onto an Andr\'e II permutation $w'$ from $\And_{n-1}^{I\!I}$ ending with~$k$ in a bijective manner.
Hence, $b_n(n,k)=b_{n-1}(k,\brullet)$, which is equal to 
$a_{n-1}(\brullet,k-1)$ by Theorem 1.1 for $1\le k\le n-1$.

\medskip
(TS4.2) The entry $b_n(n-1,k)$ counts the Andr\'e II permutations~$w$ from $\And_n^{I\!I}$ of the form $w=x_1\cdots x_{i-2}\,n\,x_i\cdots x_{n-1}\,(n-1)$ with $i\le n-1$ and~$k$ equal to $x_{i-2}$ or $x_i$.
Such a permutation can be mapped onto a permutation~$w'$ from $\And_{n-1}^{I\!I}$ defined as follows:
$$
w':=x_1\cdots x_{i-2}\,(n-1)\,x_i\cdots x_{n-1}.\leqno{(6.4)}
$$
This defines a bijection of the set of all $w$ from $\And_n^{I\!I}$ such that 
$(\bfL,\grn)w=(n-1,k)$ onto the set of all~$w'$ from $\And_{n-1}^{I\!I}$ such that
$\grn w'=k$ $(1\le k\le n-2)$. Thus, $b_n(n-1,k)= b_{n-1}(\brullet,k)$, also equal to 
$a_n(\brullet,k)$ by Theorem 1.1.

\medskip
(TS5.1)\thinspace The entry $a_n(1,\!k)$ counts the permutations $w$ from $A_{n}^I$ such~that
$(\bfF,\NL)w=(1,k)$. The bijection $1\,x_2\cdots k\,n\mapsto
(x_2-1)\cdots (k-1)\,(n-1)$ maps the set of those permutations onto the set
of all~$w'$ from $A_{n-1}^I$ such that $\grn w'=k-1$. Hence,
$a_n(1,k)=a_{n-1}(\brullet,k-1)$,

\goodbreak
\medskip
6.2. {\it Tight Andr\'e I permutations}.\quad  
As sketched in Subsection~1.6 and its  display (1.13), proving (TS5.2) and (TS5.3) 
amounts to do the following points: 

\smallskip
(a) split each set $A_n(m,k)$ into two disjoint subsets
$$A_n(m,k)=T_n(m,k)+N\!T_n(m,k),$$
in such a way that

(b) when $2\le k+1\le m\le n-2$ or $3\le m+2\le k\le n-1$ a bijection 
$$f:N\!T_n(m,k)\rightarrow A_n(m+1,k);$$

(c) and another bijection
$$
\eqalignno{\phi&:B_{n-1}(m,k)\rightarrow T_n(m,k),&\hbox{when $m>k$};\cr
\phi&:B_{n-1}(m,k-1)\rightarrow T_n(m,k),&\hbox{when $m<k$};\cr}
$$
can be duly constructed.


\medskip
Points (a) and (b).\quad Let~$f$ be the transposition of the first letter~$\bfF w=m$ within a
permutation~$w$ and the letter equal to $(m+1)$ ($1\le m\le n-2$):
$$
f: w=m\,v\,(m+1)\,v'\mapsto w'=(m+1)\,v\,m\,v'.\leqno(6.5)
$$
If $w$ is an Andr\'e I permutation,
the image $w'=f(w)$ is not always an Andr\'e~I permutation. For example, 423516 belongs to $\And^I_6$, but not $f(w)=523416$, for the trough~2 is not of type~I.  However, the reverse transposition 
$$
f^{-1}:w'=(m+1)\,v\,m\,v'\mapsto f^{-1}(w')= w=m\,v\,(m+1)\,v'
$$
whenever defined, maps each Andr\'e I permutation  onto an Andr\'e I permutation.
The Andr\'e~I permutations~$w$, whose images $f(w)$ are {\it not} Andr\'e~I  permutations
are called {\it tight}. They are characterized as follows.

\medskip
{\it Definition $6.1$}.\quad
An Andr\'e I permutation $w=m\,v\,(m+1)\,v'$ is said to be
 {\it tight}, if the following two conditions hold:

(i) either $v=e$, or $v\not=e$ and all its letters are less than~$m$;

(ii) either $v'\not=e$ and $\bfF v'$ is less than all the letters of $w$ to its left, or $v'=e$ and necessarily $m=n-1$.

\medskip
Let $T_n$ (resp. $N\!T_n$) be the subset of all Andr\'e~I permutations from $\And_n^I$, which are tight (resp. not tight), and let $T_n(m,k):=T_n\cap A_n(m,k)$, $N\!T_n(m,k):=N\!T_n\cap A_n(m,k)$.

Note that the Andr\'e I permutations from $A_n(1,k)$ are all of the form $1\,v\,2\, v'$ and, either the letters of~$v$ are all greater than~2, or $v'\not =e$ but $2<\bfF v'$, so that at least one of conditions~(i), (ii) does not hold. Accordingly, $N\!T_n(1,k)=A_n(1,k)$ for all~$k$, that is, all Andr\'e~I permutations starting with~1 are not tight. Also, note that each Andr\'e~I permutation from $A_n(n-1,k)$ is of the form $w=(n-1)\,v\,n$ and is necessarily tight, so that
$T_n(n-1,k)=A_n(n-1,k)$ for all~$k$.

\proclaim Proposition 6.1. Let $n\ge 3$ and let $w$ be a tight Andr\'e~I permutation from $\And_{n}^I$. Then, $f(w)$ (defined in $(6.5)$) cannot be an Andr\'e~I permutation.


\proof
Take the notation of (6.5) for~$w$ and $w'=f(w)$. When $m=n-1$, then $w'=n\,v\,(n-1)$ is not Andr\'e~I. When $m\le n-2$, $v=e$, and (ii) of Definition~6.1 holds, then $w'$ contains 
the double descent $(m+1)>m>\bfF v'$, therefore is not Andr\'e~I. When $m\le n-2$,
$v\not=e$ and (ii) of Definition 6.1 holds, let~$x$ be the minimum trough in~$w'$ between $(m+1)$ and $m$; then, the $x$-factorization $(w_1,w_2,x,w_4,w_5)$ of~$w'$ is such that $\max w_2w_4=m+1$ with $(m+1)$ a letter of~$w_2$. Again, $w'$ cannot be an Andr\'e~I permutation.~\qed

\proclaim Proposition 6.2.  If $2\le k+1\le m\le n-2$ or
$3\le m+2\le k\le n-1$, then $f$ maps $N\!T_n(m,k)$ onto $A_n(m+1,k)$ in a bijective manner.

\proof
To prove that~$w'$ is Andr\'e I when~$w$ is not tight, prove that (i) $w'$ has no double descent; (ii)
all the troughs of $w'$
are of type~I . 

(i) The only double descent that could be created when going from~$w$ to~$w'$ is $(m+1)>m>\bfF v'$. This could occur only if $v=e$, $v'\not=e$ and $m>\bfF v'$ and this would mean that~$w$ is tight; a  contradiction.

(ii) Let $x_i$ (resp. $x'_i$) be the $i$-th letter counted from left to right of~$w$ (resp. $w'$). Also, let $(w_1,w_2,x_i,w_4,w_5)$ (resp. $(w'_1,w'_2,x'_i,w'_4,w'_5)$) be the~$x_i$ (resp. $x'_i$)-factorization of~$w$ (resp. of~$w'$). Several cases are to be considered.

(1)  Suppose that $x_i$ is to the right of $(m+1)$ in~$w$, then $x'_i=x_i$. If~$x_i$ is a trough of~$w$, then either $(m+1)$ is a letter of~$w_2$, or not. If it is, then~$w'_2$ is derived from $w_2$ by replacing the letter $(m+1)$ by $m$. Therefore, $\max w'_2\le \max w_2<\max w_4=\max w'_4$ and the $x_i'$-factorization remains of type~I in $w'$. If it is not, then $w'_2=w_2$, $w'_4=w_4$ and the same conclusion holds.

(2) Now, suppose that $x_i=(m+1)$, so that $x'_i=m$. If $x_i$ is a trough of $w$---this is possible, as $w$ is supposed to be not tight---then, $v\not=e$ and $m$ is not a letter of~$w_2$. Furthermore, 
$(m+1)$ is a letter of~$w'_2$ only when all the letters  between~$m$ and $(m+1)$ are greater than $(m+1)$. Whatsoever, we have: $\max w'_2=\max w_2<\max w_4=\max w'_4$, so that $x'_i$ is a trough of type~I in~$w'$.

(3) Next, let $x_i$ lie between $m$ and $(m+1)$ in~$w$, so that $x'_i=x_i$ and suppose that $x_i$ is  trough of~$w$. If $x_i$ is greater than $(m+1)$,
then  $w'_2=w_2$ and $w'_4=w_4$. Moreover,  $x'_i=x_i$ will be a trough of type~I in~$w'$.
If $x_i$ is less than~$m$, the only problem arises when $m$ and $(m+1)$ are the maximum letters of~$w_2$ and $w_4$, respectively. In such a case, all the letters between~$m$ and $(m+1)$ are smaller than~$m$ and $m>\bfF v'$. Hence, $w$ would be tight. A contradiction.

Thus, the image $f(w)$ of~$w$ supposed to be not tight is Andr\'e I. If $2\le k+1\le m\le n-2$, a fortiori, $k<m+1$, so that the next to the last letter of a permutation $w$ from $N\!T_n(m,k)$, which is equal to~$k$, cannot be equal to $(m+1)$.
Thus, $f(N\!T_n(m,k))\subset A_n(m+1,k)$. In the same manner, if $3\le m+1< k\le n-1$, the inequality $m+1<k$ implies the same inclusion. As $f$ and $f^{-1}$ are inverses of each other when applied to the sets $A_n(m,k)$ and $A_n(m+1,k)$, respectively, the restriction of $f^{-1}$ to $A_n(m+1,k)$ is necessarily $N\!T_n(m,k)$ by Proposition~6.1. Thus, Proposition~6.2 is proved for $2\le k+1\le m\le n-2$ and
$3\le m+1< k\le n-1$.\qed

\smallskip
In Table 6.1 the bijection $f:N\!T_5(m,k)\rightarrow A_5(m+1,k)$ is materialized by the vertical arrows. The five tight permutations in $\And_5^I$ are reproduced in boldface.
They can only be targets of those arrows. 
This completes the program of points (a) and (b).

\medskip
\def\segment(#1,#2)\dir(#3,#4)\long#5{%
\leftput(#1,#2){\lline(#3,#4){#5}}}

{\midinsert
\def\fleche(#1,#2)\dir(#3,#4)\long#5{%
{\leftput(#1,#2){\vector(#3,#4){#5}}}}

\newbox\boxtight
\vskip 2.5cm
\def\segment(#1,#2)\dir(#3,#4)\long#5{%
\leftput(#1,#2){\lline(#3,#4){#5}}}

\setbox\boxtight=\vbox{\offinterlineskip 
\leftput(3,16){$1$}
\leftput(3,6){$2$}
\leftput(3,-5){$3$}
\leftput(3,-17){$4$}
\leftput(4,27){$\scriptstyle \NL$}
\leftput(2,23.5){$\scriptstyle \bfF$}
\leftput(16,24){1}
\leftput(33,24){2}
\leftput(49,24){3}
\leftput(66,24){4}
\segment(0,30)\dir(1,-1)\long{8}
\segment(0,30)\dir(1,0)\long {76.5}
\segment(0,22)\dir(1,0)\long {76}
\segment(0,12)\dir(1,0)\long {76}
\segment(0,0)\dir(1,0)\long {76}
\segment(0,-12)\dir(1,0)\long {76}
\segment(0,-20)\dir(1,0)\long {76}
\segment(0,30)\dir(0,-1)\long{50}
\segment(8,30)\dir(0,-1)\long{50}
\segment(25,30)\dir(0,-1)\long{50}
\segment(42,30)\dir(0,-1)\long{50}
\segment(59,30)\dir(0,-1)\long{50}
\segment(76.5,30)\dir(0,-1)\long{50}
\leftput(27,18){$13425$} 
\leftput(44,18){$12435$} 
\leftput(48,14){$14235$} 
\leftput(61,18){$12345$} 
\leftput(65,14){$13245$} 
\leftput(11,7){$23415$} 
\leftput(44,7){$21435$} 
\leftput(48,3){$24135$} 
\leftput(61,7){$21345$} 
\leftput(64,3){$\bf23145$} 
\leftput(11,-5){$\bf32415$} 
\leftput(61,-5){$31245$} 
\leftput(27,-5){$31425$} 
\leftput(30,-9){$\bf34125$} 
\leftput(27,-17){$\bf41325  $} 
\leftput(44,-17){$\bf41235$} 
\fleche(29,17)\dir(-3,-2)\long{10} 
\fleche(18,6)\dir(0,-1)\long{8} 
\fleche(28,-6)\dir(0,-1)\long{8} 
\fleche(45,17)\dir(0,-1)\long{7} 
\fleche(45,6)\dir(-2,-1)\long{16} 
\fleche(55,13.5)\dir(0,-1)\long{7} 
\fleche(55,2)\dir(-2,-1)\long{16} 
\fleche(62,17)\dir(0,-1)\long{7} 
\fleche(62,6)\dir(0,-1)\long{8} 
\fleche(72,13)\dir(0,-1)\long{7} 
\fleche(61,-6)\dir(-3,-2)\long{12} 
}

$$\box\boxtight\hskip8cm$$

\vskip2cm
\centerline{\quad Table 6.1: the bijection~$f:N\!T_5(m,k)\rightarrow A_5(m+1,k)$ $(k\not=m+1)$.\hfill}
\centerline{\hphantom{\quad Table 6.1: the bijection~}$f:A_5(m,m+1)\rightarrow A_5(m+1,m)$\hfill}

\endinsert}
{\it Remark}.\quad
When $3\le m+1=k\le n-2$ the  permutation
$w=m\,v\,(m+1)\,v'$ from $A_n(m,m+1)$ the right factor~$v'$ is equal to the one-letter word~$n$. 
This implies that $f$ maps $A_n(m,m+1)$ onto 
 $A_n(m+1,m)$ in a bijective manner. In particular, $a_n(m,m+1)=a_n(m+1,m)$.
 The fact is
illustrated in Table 6.1 by oblique arrows.

\medskip

Point (c). \quad
Let $n\ge 3$ and consider a permutation $w=x_{1}x_{2}\cdots x_{n-1}$ from $\And_{n-1}^I$. Let $x_j=\spi w$.  Define
$\phi(w):=x_{j}x'_{1}x'_{2}\cdots x'_{n-1}$, where
$$x'_{i}:=\cases {x_{i},&if $x_{i}\le x_{j}-1$;\cr
x_{i}+1,&if $x_{i}\ge x_{j}$.\cr}\leqno(6.6)
$$

The inverse bijection $\phi^{-1}$ is defined as follows: let 
$w'=x'_1x'_2\cdots x'_{n}$ belong to $T_n$; 
then, $\phi^{-1}(w'):=\rho(x'_2\cdots x'_{n})$, where $\rho$ is the reduction defined in Section 2.

\proclaim Theorem 6.3. The mapping $\phi$ is a bijection of
$\And_{n-1}^I$ onto the set~$T_n$ of all tight Andr\'e I permutations, having the following properties:\hfil\break
(i) $\spi w=\bfF \phi(w)$;\hfil\break
(ii) $\grn\phi(w)=\cases{\grn \,w,&if $\spi w>\grn w$;\cr
\grn w+1,&if $\spi w<\grn w$.\cr}
$

\vskip-12pt
{\midinsert
$$
\vbox{\offinterlineskip
\def\strat{\vrule height 12pt depth 3.5pt width 0pt}
\halign{\vrule\ \hfil$#$\quad\vrule
&\strut \quad$#$\hfil
&\quad\vrule\quad $#$\hfil
&\quad\vrule\quad $#$\hfil
&\quad\vrule\quad $#$\hfil
&\quad\vrule\quad $#$\hfil\quad\vrule\cr
\noalign{\hrule}
k={}&1&2&3&4&5\cr
\noalign{\hrule}
m=1&\cdot&\cdot&\cdot&\cdot&\cdot\cr
\noalign{\hrule}
\strat2&&&&\bf \widehat21435&\bf \widehat21345\cr
&&&&231546&231456\cr
\noalign{\hrule}
\strat3&\bf \widehat32415&\bf \widehat31425&&&\bf 2\widehat3145\cr
&342516&341526&&&324156\cr
&&&&&\bf \widehat31245\cr
&&&&&341256\cr
\noalign{\hrule}
\strat4&\bf 23\widehat415&\bf 3\widehat4125&\bf 2\widehat4135&&\cr
&423516&435126&425136&&\cr
&&\bf \widehat41325&\bf \widehat41235&&\cr
&&451326&451236&&\cr
\noalign{\hrule}
\strat5&&\bf1342\widehat5&\bf 1243\widehat5&\bf 1234\widehat5&\cr
&&513426&512436&512346&\cr
&&&\bf1423\widehat5&\bf1324\widehat5&\cr
&&&514236&513246&\cr
\noalign{\hrule}
}}
$$
\centerline{Table 6.2:  The bijection 
$\phi:B_5(m,k-1)\ {\rm (resp.\ }B_5(m,k))\rightarrow T_6(m,k)$.}
\endinsert}

\bigskip
In Table 6.2 the permutations in boldace are the  elements of $\And_5^I$. Their images under~$\phi$ are the sixteen tight permutations from~$T_6$, written in plain under them.
The box $(m,k)$ contains the permutations~$w$ from
$\And_5^I$ such that $\spi w=m$ and $\grn w=k$ (resp. $\grn w=k-1$) when $m>k$ (resp. when $m<k$).
It also contains the elements $w'$ from $T_6$
such that $\bfF w'=m$ and $\grn w'=k$. 
A hat sign  $\ \widehat{}\ $ has been put onto the spike of $w$.

\medskip
{\it Proof of Theorem $6.3$}.\quad Let $w=x_1x_2\cdots x_{n-2}x_{n-1}$ be from $\And_{n-1}^I$. Let $x_j=\spi w$. When $j=1$, then $x_1>x_2$ and
$\phi(w)=x_1\,(x_1+1)\break\,x_2 \cdots x_{n-2}'x_{n-1}'$. Accordingly, $\phi(w)$ is tight. Moreover, $\spi w=x_1=\bfF \phi(w)$, still since $x_1>x_2$. Also,
either $\grn w=x_{n-2}< x_1=\spi w$ and then $\grn \phi(w)=x'_{n-2}=x_{n-2}=\grn w$,
or $\grn w=x_{n-2}>x_1=\spi w$ and then $\grn\phi(w)=x'_{n-2}=x_{n-2}+1=\grn w+1$.

When $j\geq 2$ we have 
$$\phi(w)=x_j\,x_1\,\cdots\,x_{j-1}\,(x_j+1)\,x_{j+1}\,
x_{j+2}'\,\cdots x'_{n-1}.\leqno{(6.7)}$$ 
On the other hand, $\phi(w)$ is Andr\'e~I, because no double descent has been created; furthermore, the new trough~$x_1$ is of type~I, as the letter $(x_j+1)$ is to its right. Also, $\phi(w)$ is tight, because $x_{j+1}$ (resp. $(x_j+1)$) is less (resp. greater) than all the letters to its left. Finally, $\spi w=x_j=\bfF \phi(w)$. Moreover,
$\grn \phi(w)=x'_{n-2}$ is equal to $x_{n-2}=\grn w$ or $x_{n-2}+1=\grn w+1$, depending on whether $x_{n-2}=\grn w$ is less than or at least equal to $x_n=\spi w$.
\qed

\medskip
This  achieves the program of point (c), by definition of $B_n(m,k)$ given in (6.2).

\medskip
6.3. {\it Hooked and unhooked permutations}.\quad
Let $n\ge 3$ and consider the mapping $\Theta$, defined on $\And_{n-1}^I$ as follows.
Let $w=x_1x_2\cdots x_{n-1}$ belong to $\And_{n-1}^I$. Define:
$$
\Theta(w):=\cases{(x_1+1)x_1x'_2\cdots x'_{n-1},&if $x_1<x_2$;\cr
x_1(x_1+1)x'_2\cdots x'_{n-1},&if $x_1>x_2$;\cr}\leqno(6.8)
$$
where
$x'_i:=x_i$ (resp. $x_i+1$) if $x_i<x_1$ (resp. if $x_i>x_1$). Clearly, $\Theta$ is an {\it injection} of $\And_{n-1}^I$ into  $\And_n^I$. The permutations belonging to the subset $\Theta(\And_{n-1}^I)$ are said to be {\it hooked}. Their formal definition is next stated.

\medskip
{\it Definition $6.2$}.\quad An Andr\'e I permutation $w=x_1x_2\cdots x_n$ 
($n\geq 3$) from $\And_n^I$ is said to be {\it hooked}, if $x_1-1=x_2<x_3$  or $x_1 +1=x_2>x_3$. 

\medskip
Let $H_n$ denote the subset of all the hooked permutations from $\And_n^I$. The elements of the set-theoretic difference $N\!H_n:=\And_n^I\setminus H_n$ are said to be {\it unhooked}.
Let $H_n(m,k)$ (resp. $N\!H_n(m,k)$) denote the subset of $H_n$ (resp. of $N\!H_n$) consisting of all~$w$ such that $(\spi,\grn)w=(m,k)$. 
 
\proclaim Proposition 6.4. The injection $\Theta$ defined in $(6.8)$
of $\And_{n-1}^I$ into $\And_{n}^I$ maps $\And_{n-1}^I$ onto $H_n$.
Moreover, for each $w$ from $\And_n^I$ we have:
$$
\leqalignno{\spi \Theta(w)&=1+\bfF w;\cr
\grn \Theta(w)&=\cases{1+\grn w,&if $\bfF w<\grn w$;\cr
\grn w,&if $\bfF w>\grn w$.\cr}&(6.9)\cr}
$$

\proof
With the notation of (6.8) $\spi\Theta(w)=x_1+1$ in both cases. The identity on ``$\grn$'' follows from the very definition of~$\Theta$.\qed

\proclaim Corollary 6.5. 
The mapping~$\Theta$ is a bijection of 
$A_{n-1}(m,k)$ onto $H_n(m+1,k)$ when $1\le k< m\le n-2$, and
onto $H_n(m+1,k+1)$ when
$3\le m+2\le k\le n-1$.

In Table 6.3 have been reproduced the sixteen permutations from $\And_5^I$ in boldface and under them the hooked permutations from $H_6$, images of them under~$\Theta$.
A hat sign  $\ \widehat{}\ $ has been put onto the spike of each permutation from $H_6$.

{\midinsert

$$
\vbox{\offinterlineskip
\halign{\vrule\ \hfil$#$\quad\vrule
&\strut \quad$#$\hfil
&\quad\vrule\quad $#$\hfil
&\quad\vrule\quad $#$\hfil
&\quad\vrule\quad $#$\hfil
&\quad\vrule\quad $#$\hfil\quad\vrule\cr
\noalign{\hrule}
k={}&1&2&3&4&5\cr
\noalign{\hrule}
m=1&\cdot&\cdot&\cdot&\cdot&\cdot\cr
\noalign{\hrule}
2&&&\bf 13425&\bf 12435&\bf 12345\cr
&&&\widehat 214536&\widehat 213546&\widehat 213456\cr
&&&&\bf14235&\bf13245\cr
&&&&\widehat 215346&\widehat 214356\cr
\noalign{\hrule}
3&\bf 23415&&&\bf 21435&\bf 21345\cr
&\widehat 324516&&&2\widehat 31546&2\widehat 31456\cr
&&&&\bf 24135&\bf 23145\cr
&&&&\widehat 325146&\widehat 324156\cr
\noalign{\hrule}
4&\bf 32415&\bf 31425&&&\bf 31245\cr
&3\widehat 42516&3\widehat 41526&&&3\widehat 41256\cr
&&\bf 34125&&&\cr
&&\widehat 435126&&&\cr
\noalign{\hrule}
5&&\bf41325&\bf 41235&&\cr
&&4\widehat 51326&4\widehat 51236&&\cr
\noalign{\hrule}
}}
$$

\centerline{Table 6.3:} 
\centerline{the bijection $\Theta:A_5(m,k)\rightarrow
H_6(m+1,k)\ {\rm (resp.\ }H_6(m+1,k+1))$}

\endinsert}

We have then achieved the first two steps of the program displayed in (1.14), namely,
define the disjoint union
$B_n(m+1,k)=H_n(m+1,k)+N\!H(m+1,k)$, so that a bijection $\Theta:A_{n-1}(m,k)\rightarrow H_n(m+1,k)$ (resp. ${}\rightarrow H_n(m+1,k+1)$, for $k<m$ (resp. for $m+1<k$) can be constructed. 
The final step is devoted to the construction of the bijection $\beta$ appearing in (1.14).

\goodbreak
\medskip
6.4. {\it A bijection of $B_n(m,k)$ onto $N\!H_n(m+1,k)$}.\quad
Go back to the proofs of (TS4.1) and (TS4.2) made in \S\thinspace6.1. It was shown that
$b_n(n-1,k)=b_n(n,k+1)=b_{n-1}(\brullet,k)$ for $1\le k\le n-2$. By means of the bijections described  in \S\thinspace6.1 and Section~4,  and also the bijection~$\phi$ constructed in (5.2), we can set up a bijection of $B_n(n-1,k)$ onto $B_n(n,k+1)$. We can also proceed directly as follows.
Let $n\ge 3$ and $w=x_1\cdots x_{i-1}\,\widehat{(n-1)}\,x_{i+1}\cdots k\,n$ be an Andr\'e~I permutation such that $(\spi,\grn)w=(n-1,k)$. Then, the mapping~$\alpha$, where
$$\alpha(w):=1\,(x_1+1)\cdots (x_{i-1}+1)\,(x_{i+1}+1)\cdots (k+1)\,\widehat n,\leqno(6.10)
$$ fulfills our requirements. 

The inverse $\alpha^{-1}$ is easy to find: let $w'=x'_1\;x'_2\,\cdots \,x'_{n-1}\,n$ be a permutation from $B_n(n,k+1)$, so that $x'_1=1$, then $\alpha^{-1}(w')$ is obtained by first determining the leftmost letter $x'_{i+1}$ less than or equal to $x'_2$, and let
$$
\alpha^{-1}(w'):=(x'_2-1)\cdots (x'_{i-1}-1)\,(n-1)\,(x'_{i+1}-1)\,\cdots k\,n.\leqno(6.11)$$
The bijection $\alpha$ will be an ingredient for the next bijection~$\beta$ of
$B_n(m,k)$ onto $N\!H_n(m+1,k)$.

\medskip
First, let
$$2\le k+1\le m\le n-2\quad{\rm or}\quad 3\le m+2\le k\le n-1\leqno(6.12)
$$ 
and partition $B_n(m,k)$ into two subsets $B_n^{(1)}(m,k)$,
$B_n^{(2)}(m,k)$ as follows. Note that each permutation~$w$ from 
$B_n(m,k)$ is of the form $w=w_1mw_2(m+1)w_3$ and the factor~$w_2$ is never empty, as $m$ is the spike of $w$. Also, $w_3\not=e$ because of condition (6.12).
Say that an element of $B_n(m,k)$ belongs to $B_n^{(1)}(m,k)$ (resp. to
$B_n^{(2)}(m,k)$) if $\bfF w_3$ {\it is not} (resp. if $\bfF w_3$ {\it is}) a left minimum record, or equivalently, if 
$\min w_2<\bfF w_3$ (resp. if $\min w_2>\bfF w_3$).

Let $w=x_1x_2\cdots x_n=w_1\, m\,w_2\,(m+1)\,w_3$ be from $B_n(m,k)$ with $(m,k)$ satisfying (6.12).

(1) If $w$ belongs to $B_n^{(1)}(m,k)$,  define $w':=\beta(w)$ to be the permutation derived from~$w$ {\it by transposing} the letters~$m$ and $(m+1)$:
$$
\beta:w=w_1\, m\,w_2\,(m+1)\,w_3 \ \mapsto\ w'=w_1\, (m+1)\,w_2\,m\,w_3.\leqno(6.13)
$$

(2) If $w$ belongs to $B_n^{(2)}(m,k)$, consider the factorization $w=v_1w_3$, where $v_1=w_1 \,m\,w_2\, (m+1)$. Then, $v_1$ is Andr\'e~I by Proposition~2.1 (6). Let~$n'$ be the length of~$v_1$ and $\rho(v_1)$ 
be the {\it reduction} of~$v_1$ (by using the increasing bijection from the set  $\{x_1,\ldots,m,\ldots,m+1\}$ onto~$\{1,2,\ldots,n'\}$). Thus,  $\rho(v_1)$ is an Andr\'e~I permutation from $\And_{n'}^I$ such that $\spi\rho(v_1)=n'-1$. The bijection~$\alpha$, introduced in (6.10), can be applied to $\rho(v_1)$ and the permutation~$w':=\beta(w)$ is  defined by replacing the  
left factor~$v_1$  of~$w$ by $\rho^{-1}\alpha\rho(v_1)$:
$$
\beta:w=v_1w_3\mapsto w':=\rho^{-1}\alpha\rho(v_1)\,w_3.\leqno(6.14)
$$

{\it Example}.\quad
The permutation $w=4\,{\bf5}\, 3\,{8}\, 1\,{\bf6}\,7\,2\,9$ belongs to $B_9^{(1)}(5,2)$, as
$\min w_2=\min 381=1<7=\bfF w_3$. It then suffices to transpose~5 and~6 to get the permutation
$w'=4\,{\bf6}\ 3\,{8}\, 1\,{\bf5}\,7\,2\,9$.

Next,  $w=3\,5\,{\bf6}\, 2\,{\bf7}\,1\,8\,4\,{9}$ belongs to $B_9^{(2)}(6,4)$, as $\min w_2=2>1=\bfF w_3$. Hence, $v_1=3\,5\,6\,2\,7$, $\rho(v_1)=2\,3\,4\,1\,5$, $\alpha\rho(v_1)=
1\,3\,4\,2\,5$, $\rho^{-1}\alpha\rho(v_1)=2\,5\,6\,3\,7$ and
$w'=2\,5\,{\bf6}\,3\,{\bf7}\, 1\,8\,4\,{9}$. 
 
\smallskip
When $w$ belongs to $B_n^{(1)}(m,k)$, the letter $(m+1)$ occurs to the left of~$m$
in~$w'$. On the other hand, as $w_2$ is non-empty and $m+1>\bfF w_2$, the permutation~$w'$ is unhooked if $w_1=e$. The same conclusion also holds if $w_1\not=e$, because $\bfL w_1<m+1$ and $\bfL w_1\not=m$. Obviously, $\spi w'=m+1$ and $\grn w'=k$ by (6.12).

Let us now prove that $w'$ is Andr\'e I.
Note that the troughs remain the same in both~$w$ and~$w'$. 
Let~$x$ be a trough within~$w_2$ and $(v_1,v_2,x,v_4,v_5)$ 
(resp. $(v'_1,v'_2,x,v'_4,v'_5)$) be the $x$-factorization of~$w$ (resp. of~$w'$).
When going from~$w$ to $w'$ the type of~$x$ is not modified when at least one of the following conditions holds: $\max v_2\not=m$, $\max v_4\not=m+1$. If both were violated for a given~$x$, it would be the case for $x=\min w_2$ and all the letters of~$w_2$ would be less than~$m$. But
$\max v_4=m+1$ implies $\max v_4>\bfF w_3>\min w_2$ and $\bfF w_3$ is a trough of~$w$.  
If $(v''_1,v''_2,\bfF w_3,v''_4,v''_5)$ is the $\bfF w_3$-factorization of~$w$, the word~$\bfF w_3\,v''_4$ is necessarily a factor of~$v_4$, as all its letters are greater than $\min w_2$.
Hence, $\max v_4> m+1$, a contradiction. Thus,

\smallskip
$w'$ is an {\it unhooked permutation
from $\And^I_n$ such that $\spi w'\!=\!m\!+\!1$, $\grn w'\!=\!k$ with~$(m\!+\!1)$ to the left of $m$}. In short,
$w'\in N\!H_n^{(1)}(m+1,k)$.

\smallskip
As the transposition
$w_1\, (m+1)\,w_2\,m\,w_3 \ \mapsto\ w_1\, m\,w_2\,(m+1)\,w_3$, when applied to Andr\'e~I permutations with $(m+1)$ to the left of~$m$, always maps an Andr\'e~I onto an Andr\'e~I permutations, 

\smallskip
{\it the direct transposition~$\beta$ defined in $(6.13)$ is a bijection of $B_n^{(1)}(m,k)$ onto 
$N\!H_n^{(1)}(m+1,k)$.}

\smallskip
Next, let $w$ belong to $B_n^{(2)}(m,k)$ and consider the permutation $w'=\beta(w)$ defined in (6.14).
The left factor $\rho^{-1}\alpha\rho(v_1)$ of $\beta(w)$ is Andr\'e~I and ends with $(m+1)$. Therefore, $\beta(w)$ is of the form $w'_1\,m\,w'_2\,(m+1)\,w_3$. Again,
with the hypothesis (6.12) the letter $x_{n-1}$, equal to~$k$ in the permutation $w=x_1x_2\cdots k\,n$  remains untouched when going from~$w$ to ~$w'$. Thus, $\grn w'=\grn w=k$.
Next,  we get $\spi\rho(v_1)=n'-1$ and
$\spi\alpha\rho(v_1)=n'$; hence, $\spi\rho^{-1}\alpha\rho(v_1)=m+1$. As $w_3$ starts with a letter less than all the letters in~$v_1$, we have: $\spi w'=\spi \rho^{-1}\alpha\rho(v_1)\,w_3=m+1$. 
Moreover, $w_3$ is Andr\'e I by Proposition 2.1 (5),
so that
$\beta(w)$ is  Andr\'e~I by Proposition 2.1 (7). This shows that

\smallskip
{\it the mapping $\beta$ defined in $(6.14)$ is a bijection of $B_n^{(2)}(m,k)$ onto the set $N\!H_n^{(2)}(m+1,k)$, defined as the set of all unhooked permutations
from $\And^I_n$ such that $\spi w'=m+1$, $\grn w'=k$ with~$m$ to the left of $(m+1)$}.

\goodbreak
\smallskip
This proves the following theorem.

\proclaim Theorem 6.6. Under condition $(6.12)$ the mapping $\beta:w\mapsto w'$ defined in $(6.13)$ and $(6.14)$ is a bijection of $B_n(m,k)
=B_n^{(1)}(m,k)+B_n^{(2)}(m,k)$ onto $N\!H_n(m+1,k)=N\!H^{(1)}_n(m+1,k)
+N\!H^{(2)}_n(m+1,k)$.

{\midinsert
\def\segment(#1,#2)\dir(#3,#4)\long#5{%
\leftput(#1,#2){\lline(#3,#4){#5}}}

\def\fleche(#1,#2)\dir(#3,#4)\long#5{%
{\leftput(#1,#2){\vector(#3,#4){#5}}}}

\newbox\boxunhook
\vskip 2cm
\def\segment(#1,#2)\dir(#3,#4)\long#5{%
\leftput(#1,#2){\lline(#3,#4){#5}}}

\setbox\boxunhook=\vbox{\offinterlineskip 
\leftput(5,19){$1$}
\leftput(5,12){$2$}
\leftput(5,-2){$3$}
\leftput(5,-25){$4$}
\leftput(5,-46){$5$}
\leftput(5.5,25.5){$\scriptstyle k$}
\leftput(3.5,23.5){$\scriptstyle m$}
\leftput(16,24){1}
\leftput(34,24){2}
\leftput(54,24){3}
\leftput(78,24){4}
\leftput(94,24){5}
\leftput(10,-10){$\bf 324516$}
\leftput(27,-10){$314526$}
\segment(3,27.5)\dir(1,-1)\long{5}
\segment(3,27.5)\dir(1,0)\long {107}
\segment(3,22.5)\dir(1,0)\long {107}
\segment(3,17.5)\dir(1,0)\long {107}
\segment(3,7)\dir(1,0)\long {107}
\segment(3,-12.5)\dir(1,0)\long {107}
\segment(3,-38)\dir(1,0)\long {107}
\segment(3,-62.5)\dir(1,0)\long {107}
\segment(3,27.5)\dir(0,-1)\long{90}
\segment(8,27.5)\dir(0,-1)\long{90}
\segment(25.5,27.5)\dir(0,-1)\long{90}
\segment(45,27.5)\dir(0,-1)\long{90}
\segment(68,27.5)\dir(0,-1)\long{90}
\segment(86.5,27.5)\dir(0,-1)\long{90}
\segment(110,27.5)\dir(0,-1)\long{90}
\fleche(11,-11)\dir(0,-1)\long{5}
\leftput(10,-19){$423516$}
\leftput(12.5,-23){$\bf342516$}
\fleche(28,-11)\dir(0,-1)\long{8}
\leftput(27,-23){$413526$}
\leftput(30,-27){$\bf341526$}
\leftput(32,-31){$\bf435126$}
\leftput(47,-18){$241536$}
\leftput(50,-22){$412536$}
\leftput(53,-26){$415236$}
\leftput(56,-30){$425136$}
\fleche(11,-20)\dir(0,-1)\long{23.5}
\leftput(10,-47){$234516$}
\fleche(14,-24)\dir(0,-1)\long{16}
\leftput(13,-43){$243516$}
\leftput(33,-43){$345126$}
\fleche(34,-32)\dir(0,-1)\long{8}
\leftput(30,-47){$351426$}
\leftput(27,-51){$513426$}
\fleche(28,-24)\dir(0,-1)\long{24}
\fleche(31,-28)\dir(0,-1)\long{16}
\leftput(27,-56){$\bf451326$}
\leftput(56,-43){$245136$}
\leftput(53,-47){$514236$}
\leftput(50,-51){$512436$}
\leftput(47,-55){$251436$}
\leftput(47,-60){$\bf451236$}
\fleche(48,-19)\dir(0,-1)\long{33}
\fleche(51,-23)\dir(0,-1)\long{25}
\fleche(54,-27)\dir(0,-1)\long{17}
\fleche(57,-31)\dir(0,-1)\long{9}
\leftput(47,13){$\bf214536$}
\leftput(70,13){$\bf213546$}
\leftput(73,9){$\bf215346$}
\leftput(73,1){$315246$}
\leftput(70,-3){$312546$}
\leftput(70,-7){$\bf231546$}
\leftput(70,-11){$\bf325146$}
\fleche(71,12)\dir(0,-1)\long{11}
\fleche(74,8.5)\dir(0,-1)\long{4.5}
\leftput(88,13){$\bf213456$}
\leftput(95,9){$\bf214356$}
\fleche(106,8)\dir(0,-1)\long{4}
\fleche(106,0)\dir(0,-1)\long{27}
\fleche(103,-4)\dir(0,-1)\long{19}
\fleche(100,-7.5)\dir(0,-1)\long{11.5}
\leftput(96,1){$314256$}
\leftput(93,-3){$\bf231456$}
\leftput(90,-7){$\bf324156$}
\leftput(88,-11){$312456$}
\fleche(89,12)\dir(0,-1)\long{20}
\fleche(89,-11.5)\dir(0,-1)\long{4}
\leftput(88,-18.5){$412356$}
\leftput(91,-22.5){$234156$}
\leftput(93,-26.5){$241356$}
\leftput(96,-30.5){$413256$}
\leftput(96,-36){$\bf341256$}
}

$$\box\boxunhook\hskip11cm$$

\vskip6.1cm
\centerline{Table 6.4: the bijection~$\beta:B_6(m,k)\rightarrow N\!H_6(m+1,k)$}
\endinsert}

{\it Example}.\quad In Table 6.4 the image $\beta(w)$ of each Andr\'e~I permutation~$w$ 
from $B_6(m,k)$, with $(m,k)$ satisfying inequalities (6.12) for $n=6$,
is indicated by a downarrow. The hooked permutations are reproduced in boldface. Note that they are not bottoms of any downarrows, as $\beta$ is a bijection of $B_n(m,k)$ onto $N\!H_n(m+1,k)$.

\medskip
With the  construction of the bijection $\beta:B_n(m,k)\rightarrow N\!H_n(m+1,k)$ the program displayed in (1.14) is completed, as
$\Deltaa_{(1)}b_n(m,k)=\#B_n(m+1,k)-\#B_n(m,k)
=\#B_n(m+1,k)-\#N\!H_n(m+1,k)=\#H_n(m+1,k)
=\#A_{n-1}(m,k)$ (resp. ${}=\#A_{n-1}(m,k-1)$)
if $1\le k< m\le n-2$ (resp. if $3\le m+2\le k\le n-1$).

\bigskip
\centerline{\bf 7. The making of Seidel Triangle Sequences} 

\medskip
7.1. {\it The Seidel tangent-secant matrix}.\quad In the sequel,
three exponential generating functions will be attached to each infinite matrix $A=(a(m,k))_{m,k\geq 0}$
$$\displaylines{
	A(x,y):=\sum_{m,k\geq 0} a(m,k) {x^m\over m!}{y^k\over k!}; \cr 
	A_{m, \brullet}(y):=\sum_{k\geq 0}a(m,k) {y^k\over k!};\qquad 
	A_{\brullet, k}(x):=\sum_{m\geq 0} a(m,k) {x^m\over m!};\cr} 
$$
for $A$ itself, its $m$-th row, its $k$-th column.
Let $\overline H=(\overline h_{i,j})$ $(i,j\ge 0)$ be the infinite matrix, whose entries
are  the Entringer numbers $E_n(m)$ displayed along the skew-diagonals with the following sign:
$$\displaylines{\rlap{(7.1)}\hfill
\overline h_{i,j}=\cases{(-1)^n\,E_{i+j+1}(j+1),&if $i+j=2n$;\cr
(-1)^n\,E_{i+j+1}(i+1),&if $i+j=2n-1$;\cr}
\hfill\cr
\noalign{\hbox{or still}}
\rlap{(7.2)}\hfill
E_{2n+1}(j+1)=(-1)^n\,\overline h_{2n-j,j}
\quad(0\le j\le 2n);\hfill\cr
\rlap{(7.3)}\hfill\hskip37pt
E_{2n}(i+1)=(-1)^n\,\overline h_{i,2n-1-i}
\quad(0\le i\le 2n-1);\hfill\cr
}
$$
or still in displayed form:
{
$$\eqalignno{
\overline H\!=\!&\pmatrix{
	E_1(1)&-E_2(1)&0&E_4(1)&0&-E_6(1)&\!\!0\, \cdots\!\! &\!\!\!\cr
0&-E_3(2)&E_4(2)&E_5(4)&-E_6(2)&-E_7(6)\cr
	-E_3(1)&E_4(3)&E_5(3)&-E_6(3)&-E_7(5)\cr
0&E_5(2)&-E_6(4)&-E_7(4)\cr
E_5(1)&-E_6(5)&-E_7(3)\cr
0&-E_7(2)\cr
-E_7(1)\cr
	\vdots & \cr
}\cr
&=\pmatrix{1&-1&0&2&0&-16&0&\cdots\cr
0&-1&2&2&-16&-16\cr
-1&1&4&-14&-32\cr
0&5&-10&-46\cr
5&-5&-56\cr
0&-61\cr
-61\cr
\vdots\cr
}.\cr}
$$
}

As noted by Dumont [Du82], the definition of such a matrix~$\overline H$ goes back to Seidel himself [Se1877]. Entringer [En66] rediscovered the absolute values of the entries, when he classified the alternating permutations according to their first letters. The entries of the top row are the coefficients of the Taylor expansion of $1-\tanh y=2/(1+e^{2y})$:
$$\leqalignno{\overline H_{0,\brullet}(y)
=1-\tanh y&=1+\sum_{n\ge 1} {y^{2n-1}\over
(2n-1)!}(-1)^n
E_{2n-1}\cr
&=1-{y\over 1!}1\!+\!{y^3\over 3!}2\!-\!{y^5\over 5!}16\!+\!{y^7\over
7!}272\!-\! {y^9\over 9!}7936\!+\!\cdots\cr
}
$$
The entries of the leftmost column are the coefficients of the Taylor expansion of $1/\cosh x
=2\,e^x/(1+e^{2x})$, so that 
$$\leqalignno{
\overline H_{\brullet,0}(x)={1\over \cosh x}
&=\sum_{n\ge 0}{x^{2n}\over (2n)!}(-1)^n\,E_{2n}\cr
&=1-{x^2\over 2!}+{x^4\over 4!}5-{x^6\over 6!}61+{x^8\over 8!}1385-\cdots\cr}
$$
By means of recurrence (1.1) satisfied by the Entringer numbers and (7.1) we can verify that the entries $\overline h_{i,j}$ obey the following rule:
$\overline h_{i,j}=\overline h_{i-1,j}+\overline h_{i-1,j+1}$ for 
$j\ge 0$, $i\ge 1$, so that the entries $\overline h_{i,j}$ can be obtained by applying such a rule inductively, the entries of the top row being given. Such a matrix is called a {\it Seidel matrix} by Dumont [Du82], and its exponential generating function is directly obtained from the exponential generating function 
for its top row by the formula
$\overline H(x,y)=\overline H_{0,\brullet}(x+y)\,e^x$ (see, e.g., [DV80]). Accordingly,
  $$
 \overline H(x,y)={2\,e^x\over 1+e^{2x+2y}}.\leqno(7.4)$$

\goodbreak
\medskip
Two further matrices are derived from $\overline H$. The first one,
$\overline H_1$, is obtained  by replacing all the entries $\overline h_{i,j}$ such that $i+j$ is {\it odd} by zero, so that
$$\overline H_1\!=
\pmatrix{
	1     & \cdot & 0    & \cdot& 0 & \cdot & 0 & \cdots &\cr
	\cdot &    -1 &\cdot &2 &\cdot &-16 &  & \cr
	-1     & \cdot & 4    &\cdot& -32 & &  &\cr
	\cdot &    5 &\cdot &-46 & &&  & \cr
	5     & \cdot & -56    && & &  &\cr
	\cdot &    -61 & && &&  & \cr
	-61     & &    && & &  &\cr
	\vdots & \cr
}.
$$
As $\displaystyle
\overline H(x,y)= {2e^x\over 1+e^{2x+2y}}$, we get:
$$
\overline H_1(x,y)= {\overline H(x,y)+\overline H(-x,-y) \over 2}
=e^{x} {1+e^{2y} \over 1+e^{2x+2y}}
={\cosh y\over \cosh(x+y)}.\leqno(7.5)
$$

\goodbreak
The second one, $\overline H_2$, is derived from $\overline H$ by replacing the entries $\overline h_{i,j}$ such that $i+j$ is {\it even} by~0, so that
$$
\overline H_2=\pmatrix{\cdot&-1&\cdot&2&\cdot&-16&\cdot&272 &\cdots\cr
0&\cdot&2&\cdot&-16&\cdot&272\cr
\cdot&1&\cdot&-14&\cdot&256\cr
0&\cdot&-10&\cdot&224\cr
\cdot&-5&\cdot&178\cr
0&\cdot&122\cr
\cdot&61\cr
0\cr
\vdots \cr
}.
$$

Therefore,
$$\
\overline H_2(x,y)= {\overline H(x,y)-\overline H(-x,-y) \over 2}
=e^{x} {1-e^{2y} \over 1+e^{2x+2y}}
={-\sinh y\over \cosh(x+y)}.\leqno(7.6)
$$

In the sequel, further matrices will be derived from $\overline H_1$ and $\overline H_2$, essentially by transposing them and/or removing either their top rows, or leftmost columns. The corresponding actions on their respective exponential generating functions $\overline H_1(x,y)$ and $\overline H_2(x,y)$ are the exchange of the variables~$x$ and $y$: $TH_i(x,y):=H_i(y,x)$; then, the {\it partial derivatives} with respect to~$x$ and~$y$: $D_xH_i(x,y)$ and $D_yH_i(x,y)$
$(i=1,2)$.

\medskip
\goodbreak
7.2. {\it The generating function for the Entringer numbers}.\quad
The generating function for the Entringer numbers, already derived in [FH14],
can be obtained from relations (7.5) and (7.6). In fact, they are simply
equal to $\overline H_1(xI, yI)$ and $I\overline H_2(xI, yI)$ 
with $I=\sqrt{-1}$. Thus, 
$$\leqalignno{
\sum_{1\le k\le 2n+1}E_{2n+1}(k){x^{2n+1-k}\over (2n+1-k)!}{y^{k-1}\over (k-1)!}
&={\cos y\over \cos(x+y)};&(7.7)\cr
\sum_{1\le k\le 2n}E_{2n}(k){x^{k-1}\over (k-1)!}{y^{2n-k}\over (2n-k)!}
&={\sin y\over \cos(x+y)}.&(7.8)\cr}
$$

\medskip

7.3. {\it Seidel Triangle Sequences}.\quad For calculating the generating functions for the twin Seidel matrices we shall recourse to the techniques developed in our previous paper [FH14] for the so-called Seidel triangle sequences. Only definitions will be stated, as well as the main result.

\medskip
A sequence of square matrices $(C_{n})$ $(n\ge 1)$ 
is called a {\it Seidel triangle sequence} if the following three conditions are fulfilled:

(STS1) each matrix $C_n$ is of dimension $n$;

(STS2) each matrix $C_{n}$ has null entries along and below its diagonal; let $(c_{n}(m,k))$ ($0\le m < k\le n-1$) denote its entries strictly above its diagonal, so that
{
$$\displaylines{
C_{1}=\pmatrix{\cdot\cr };\quad
C_{2}=\pmatrix{\cdot &c_{2}(0,1)\cr
\cdot&\cdot\cr};\quad
\hfill
C_{3}=\pmatrix{\cdot&c_{3}(0,1)&c_{3}(0,2)\cr
\cdot&\cdot&c_{3}(1,2)\cr
\cdot&\cdot&\cdot\cr
};\ \ldots\ ;\hfill\cr
C_{n}=\pmatrix{
	\cdot&c_{n}(0,1)&c_{n}(0,2)&\cdots
&\cdot&c_{n}(0,n-2)&c_{n}(0,n-1)\cr
	\cdot&\cdot&c_{n}(1,2)&\cdots&\cdot&c_{n}(1,n-2)&c_{n}(1,n-1)\cr
\vdots&\vdots&\vdots&\ddots&\vdots&\vdots&\vdots\cr
	\cdot&\cdot&\cdot&\cdots&\cdot&\cdot&c_{n}(n-2,n-1)\cr
\cdot&\cdot&\cdot&\cdots&\cdot&\cdot&\cdot\cr}\!;\cr
}
$$
the dots ``$\cdot$'' along and below the diagonal referring to null entries.
}

\goodbreak
\smallskip
(STS3) for each $n\ge 3$, the following relation  holds:
$$c_{n}(m,k)-c_{n}(m,k+1)=c_{n-1}(m,k) \quad (m<k).$$

\goodbreak
Record the last columns of the triangles $C_{2}$, $C_{3}$, $C_{4}$, $C_{5}$, \dots\ , read from top to bottom, namely, 
$c_{2}(0,1)$;\quad  $c_3(0,2)$, $c_3(1,2)$;\quad   $c_{4}(0,3)$, $c_{4}(1,3)$, 
$c_{4}(2,3)$;\quad 
$c_{5}(0,4)$, $c_{5}(1,4)$, $c_{5}(2,4)$, $c_{5}(3,4)$;
\dots\ as skew-diagonals of an infinite matrix
$H=(h_{i,j})_{i,j\geq 0}$, as shown next: 
$$H:=\bordermatrix{&0&1&2&3&4 & \cdots\cr
0&c_{2}(0,1)&c_3(1,2)&c_{4}(2,3)&c_5(3,4)&c_{6}(4,5) & \cdots\cr
	1&c_3(0,2)&c_{4}(1,3)&c_5(2,4)&c_{6}(3,5)\cr
	2&c_{4}(0,3)&c_5(1,4)&c_{6}(2,5)\cr
	3&c_5(0,4)&c_{6}(1,5)\cr
4&c_{6}(0,5)\cr
\vdots & \vdots \cr
},\leqno{(7.9)}
$$
In an equivalent manner, the entries of $H$ are defined by:
$$
h_{i,j}= c_{i+j+2}(j,i+j+1).\leqno{(7.10)}
$$
The next theorem has been proved in [FH14] and will be of great use in the next sections.

\proclaim Theorem 7.1. 
The three-variable generating function for the 
Seidel triangle sequence 
$(C_{n}=(c_{n}(m,k)))_{n\geq 1}$ 
is equal to
$$
\sum_{1\le m+1\le k\le n-1}\kern-20pt  {c_{n}(m,k)} 
{x^{n-k-1}\over (n-k-1)!} 
{y^{k-m-1}\over (k-m-1)!}
{z^m\over m!}
= e^x H(x+y,z),\leqno(7.11)
$$
where $H$ is the infinite matrix defined in $(7.10)$. 

\goodbreak
With  $I:=\sqrt{-1}$
we get:
$$\displaylines{(7.12)\quad
\sum_{1\le m+1\le k\le n-1}\kern-20pt I^{n-2}
c_n(m,k){x^{n-k-1}\over (n-k-1)!}{y^{k-m-1}\over (k-m-1)!}{z^m\over m!}\hfill\cr
\hfill{}=e^{Ix}H(Ix+Iy,Iz).\quad\cr}
$$

\bigskip
\centerline{\bf 8. Trivariate generating functions}

\medskip
Each of the sequences ${\rm Twin}^{(1)}:=(A_2,B_3,A_4,B_5,A_6,\ldots\,)$, ${\rm Twin}^{(2)}:=(B_2,A_3,B_4,A_5,B_6,\ldots\,)$
(see Diagram~1.3) gives rise to two Seidel Triangle sequences, by considering the upper and lower triangles of the matrices. 

\medskip
8.1. {\it The upper triangles of ${\rm Twin}^{(1)}$}.\quad The Seidel Triangle sequence to be constructed is the following: first, $C_1:=(\cdot)$, then for $n\ge 2$
each $C_n$ will be derived from the upper triangle of $A_{n+1}$ (resp. $B_{n+1}$) by (i) dropping the rightmost column; (ii) transposing the remaining triangle with respect to its skew-diagonal; (iii) changing the signs of its entries according the following rule. More precisely,
$$
C_n:=(-1)^{(n+1)/2}\petitematrice{\cdot\ &a_{n+1}(n-1,n)&\cdots&a_{n+1}(2,n)&a_{n+1}(1,n)\cr
&&\ddots&\vdots&\vdots\cr
&&\cdot&a_{n+1}(2,3)&a_{n+1}(1,3)\cr
&&&\cdot&a_{n+1}(1,2)\cr
&&&&\cdot\cr}\ \hbox{if $n$ odd;}
$$
$$
C_n:=(-1)^{n/2}\petitematrice{\cdot\ &b_{n+1}(n-1,n)&\cdots&b_{n+1}(2,n)&b_{n+1}(1,n)\cr
&&\ddots&\vdots&\vdots\cr
&&\cdot&b_{n+1}(2,3)&b_{n+1}(1,3)\cr
&&&\cdot&b_{n+1}(1,2)\cr
&&&&\cdot\cr}\ \hbox{if $n$ even;}
$$

By referring to Diagram 1.3 we get:  $C_1=\cdot\ $;\quad  $C_2=\matrice{\cdot\ &0\cr &\cdot\ \cr}$;\quad
$C_3=\matrice{\cdot\ &1&1\cr
&\cdot&1\cr
&&\cdot\cr}$;

\noindent
$C_4=\matrice {\cdot\ &2&1&0\cr
&\cdot&1&0\cr
&&\cdot&0\cr
&&&\cdot\cr}$\ ;\quad
$C_5=\matrice{\cdot\ &-2&-4&-5&-5\cr
&\cdot&-4&-5&-5\cr
&&\cdot&-4&-4\cr
&&&\cdot&-2\cr
&&&&\cdot\cr}
$\ ;\quad
$C_6=\matrice{\cdot\ &-16&-14&-10&-5&\ 0\cr
&\cdot&-14&-10&-5&\ 0\cr
&&\cdot &-8&-4&\ 0\cr
&&&\cdot&-2&\ 0\cr
&&&&\cdot&\ 0\cr
&&&&&\ \cdot\cr}$;

$C_7=\matrice{
   \cdot &  16 &  32 &  46 &  56 &  61 &  61  \cr
 &   \cdot &  32 &  46 &  56 &  61 &  61 \cr
 &   &   \cdot &  44 &  52 &  56 &  56  \cr
&  &  &   \cdot &  44 &  46 &  46  \cr
&  &   &   &   \cdot &  32 &  32 \cr
&   &  &  &  &   \cdot &  16 \cr
&  &  &  &  &  &  \cdot\cr
}
$ ;\quad
$C_8=\matrice{\cdot &272&256&224&178&122&61&\ 0\cr
&   \cdot &  256&224&178&122&61&\ 0\cr
 &   &   \cdot &208 &  164 &  112 &  56&\ 0\cr
&  &  &   \cdot &  136 &  92 &  46&\ 0  \cr
&  &   &   &   \cdot &  64 &  32&\ 0\cr
&   &  &  &  &   \cdot &  16 &\ 0 \cr
&  &  &  &  &  &  \cdot&\ 0 \cr
  &  &  &  &  &  &&\ \cdot\cr
}$\ .

Therefore, $$c_n(m,k)=\cases{(-1)^{(n+1)/2}a_{n+1}(n-k ,n-m),& if $n$ is odd;\cr
(-1)^{n/2}b_{n+1}(n-k ,n-m),& if $n$ is even.\cr}\leqno(8.1)$$

\proclaim Proposition 8.1. The sequence $(C_n)$ $(n\ge 1)$ just defined is a Seidel Triangle sequence.

\proof
Just verify that rule (STS3) holds. If $n$ is odd and $0\le m<k\le n-2$, then
$3\le m'+2:=(n-k-1)+2\le k':=n-m\le (n+1)-1$ and
$$\displaylines{\quad c_n(m,k)-c_n(m,k+1)\hfill\cr
\kern1cm{}=(-1)^{(n+1)/2}\bigl(a_{n+1}(n-k ,n-m)-a_{n+1}( n-k-1,n-m)\bigr)\hfill\cr
\kern1cm{}=(-1)^{(n+1)/2}\Deltaa_{(1)}a_{n+1}( n-k-1,n-m)\hfill\cr
\kern1cm{}=(-1)^{(n+1)/2}\Deltaa_{(1)}a_{n+1}( m',k')\hfill\cr
\kern1cm{}=(-1)^{(n-1)/2}\,b_n(m',k'-1)\hfill\hbox{[by rule (TS5.2)]}\cr
\kern1cm{}=(-1)^{(n-1)/2}\,b_n(n-k-1,n-m-1)\hfill\cr
\kern1cm{}=(-1)^{(n-1)/2}(-1)^{(n-1)/2}\,c_{n-1}(m,k)
=c_{n-1}(m,k).\hfill\cr
}
$$
The case when $n$ is even can be proved in a similar way.\qed
\goodbreak

The next step is to determine the matrix~$H$, as defined in (7.9), whose skew-diagonals are equal to the rightmost columns of the matrices~$C_n$.
For $n\ge 2$ the skew-diagonal 
$(c_{n}(0,n-1), c_{n}(1,n-1),\ldots,c_{n}(n-2,n-1))$ of~$H$, being the rightmost column of~$C_n$, is equal to
$$\eqalignno{
&\cases{(-1)^{(n+1)/2}(a_{n+1}(1,n),a_{n+1}(1,n-1),\ldots,a_{n+1}(1,2)),& if $n$ is odd,\cr
(-1)^{n/2}(b_{n+1}(1,n),b_{n+1}(1,n-1),\ldots,b_{n+1}(1,2)),&if $n$ is even;\cr}
\cr
\noalign{\hbox{also equal to}}
&\cases{(-1)^{(n+1)/2}(b_n(\brullet,n-1),b_n(\brullet,n-2),\ldots,
b_n(\brullet,1)),& if $n$ is odd;\cr
(0, 0,\ldots,0\,),&if $n$ is even;\cr}\cr
\noalign{\hbox{by Rules (TS5.1) and (TS2); finally, equal to}}
&\cases{(-1)^{(n+1)/2}(E_n(1),E_n(2),\ldots,E_n(n-1)),& if $n$ is odd;\cr
(0, 0,\ldots,0\,),&if $n$ is even;\cr}\cr
\noalign{\hbox{by (1.11).}}
}
$$

\goodbreak
Thus, 
$$\leqalignno{H&=\pmatrix{0&E_3(2)&0&-E_5(4)&0&E_7(6)&\!\!0\ \cdots\!\! &\!\!\cr
E_3(1)&0&-E_5(3)&0&E_7(5)&0\cr
0&-E_5(2)&0&E_7(4)&0\cr
-E_5(1)&0&E_7(3)&0\cr
0&E_7(2)&0\cr
E_7(1)&0\cr
0\cr
\vdots \cr
} \cr
\noalign{\vskip -4pt}
&& {(8.2)} \cr
\noalign{\vskip -4pt}
&=\pmatrix{0&1&0&-2&0&16&0 & \cdots &\cr
1&0&-4&0&32&0\cr
0&-5&0&46&0\cr
-5&0&56&0\cr
0&61&0\cr
61&0\cr
0\cr
\vdots \cr
}.\quad\cr}
$$

This matrix is to be compared with the matrix~$\overline H_1$ (see \S7.1).
For getting~$H$ it suffices to delete the top row of~$\overline H_1$ and change the signs of all the entries. As 
$\overline H_1(x,y)=\cosh y/\cosh(x+y)$ by (7.5), we have:
$$\leqalignno{
H(x,y)&=-D_x\overline H_1(x,y)
={\cosh y\,\sinh (x+y)\over \cosh^2(x+y)}.&(8.3)\cr
	\noalign{\hbox{Hence, the right-hand side of (7.11) becomes}}
e^x H(x+y,z)&=e^x{\cosh z\,\sinh(x+y+z)\over \cosh^2(x+y+z)};
\cr
\noalign{\hbox{and the right-hand side of (7.12) is equal to}}
e^{Ix}H(Ix+Iy,Iz)&=
(\cos x+I\sin x){I\,\cos z\,\sin (x+y+z)\over \cos^2(x+y+z)}.\cr}
$$

It remains to interpret the left-hand side of identity (7.12) by using (8.1).
If $n=2l+1$, then $I^{n-2}=(-1)^{l+1}I$ and
$(-1)^{(n+1)/2}=(-1)^{l+1}$. Thus,
$I^{n-2}c_n(m,k)=I\,a_{n+1}(n-k ,n-m)$.
The imaginary part of identity (7.12) then reads:
$$\displaylines{
\sum_{\scriptstyle 1\le m+1\le k\le n-1
\atop \scriptstyle n\ {\rm odd}}
\kern-20pt 
a_{n+1}(n-k,n-m){x^{n-k-1}\over (n-k-1)!}{y^{k-m-1}\over (k-m-1)!}{z^m\over m!}
\hfill\cr
\noalign{\vskip-9pt}
\hfill{}=
{\cos x\,\cos z\sin (x+y+z)\over \cos^2(x+y+z)}.\quad\cr
}$$
With the change of variables $n\leftarrow 2n-1$, $n-k\leftarrow m$, $n-m\leftarrow k$, we get (1.15) from Theorem 1.3.
Note that the above generating function involves all the matrices $A_4$, $A_6$, \dots~
of ${\rm Twin}^{(1)}$, but not the very first term $A_2=\bigl({{1\; 0}\atop {0\;0}}\bigr)$. 
\medskip

If $n=2l$, then $I^{n-2}=(-1)^{l-1}$ and $(-1)^{n/2}=(-1)^l$, so that
$I^{n-2}c_n(m,k)=-b_{n+1}(n-k ,n-m)$.
As for the real part,
$$\displaylines{
\sum_{\scriptstyle 1\le m+1\le k\le n-1
\atop \scriptstyle n\ {\rm even}}
\kern-20pt 
b_{n+1}(n-k,n-m){x^{n-k-1}\over (n-k-1)!}{y^{k-m-1}\over (k-m-1)!}{z^m\over m!}
\hfill\cr
\noalign{\vskip-9pt}
\hfill{}=
{\sin x\,\cos z\sin (x+y+z)\over \cos^2(x+y+z)}.\quad\cr}$$
With the change of variables
$n\leftarrow 2n$, $n-k\leftarrow m$, $n-m\leftarrow k$, we get (1.21) from Theorem 1.6.

\medskip
8.2. {\it The upper triangles of ${\rm Twin}^{(2)}$}.\quad
The sequence of triangles to be considered is the following:
$C_1=\cdot\ $ and for $n\ge 2$
$$
C_n:=(-1)^{(n-1)/2}\petitematrice{\cdot\ &b_{n+1}(n-1,n)&\cdots&b_{n+1}(2,n)&b_{n+1}(1,n)\cr
&&\ddots&\vdots&\vdots\cr
&&\cdot&b_{n+1}(2,3)&b_{n+1}(1,3)\cr
&&&\cdot&b_{n+1}(1,2)\cr
&&&&\cdot\cr}\ \hbox{if $n$ odd;}
$$
$$
C_n:=(-1)^{n/2}\petitematrice{\cdot\ &a_{n+1}(n-1,n)&\cdots&a_{n+1}(2,n)&a_{n+1}(1,n)\cr
&&\ddots&\vdots&\vdots\cr
&&\cdot&a_{n+1}(2,3)&a_{n+1}(1,3)\cr
&&&\cdot&a_{n+1}(1,2)\cr
&&&&\cdot\cr}\ \hbox{if $n$ even;}
$$
that is, 
$C_1=\cdot$,\quad $C_2=\matrice{\cdot&-1\cr &\cdot\cr}$;\quad
$C_3=\matrice{\cdot&-1&0\cr  &\cdot&0\cr  &&\cdot\cr}$;\quad
$C_4=\matrice{\cdot&1&2&2\cr &\cdot&2&2\cr &&\cdot&1\cr &&&\cdot\cr}$;\quad

\indent\indent
$C_5=\matrice{\cdot&5&4&2&\ 0\cr &\cdot&4&2&\ 0\cr  
&&\cdot&1&\ 0\cr  &&&\cdot&\ 0\cr &&&&\cdot\cr}$;\quad
$C_6=\matrice{\cdot&-5&-10&-14&-16&-16\cr &\cdot &-10&-14&-16&-16\cr
&&\cdot&-13&-14&-14\cr  &&&\cdot &-10&-10\cr &&&&\cdot&-5\cr &&&&&\cdot\cr}$;

\noindent
$C_7=\kern-6pt \matrice{\cdot&-61&-56&-46&-32&-16&0\cr
&\cdot&-56&-46&-32&-16&0\cr
&&\cdot&-41&-28&-14&0\cr
&&&\cdot&-20&-10&0\cr
&&&&\cdot&-5&0\cr
&&&&&\cdot&0\cr
&&&&&&\cdot\cr}
$;\quad
$C_8=\kern-11pt\matrice{\cdot\ &61&122&178&224&256&272&272\cr
&\cdot&122&178&224&256&272&272\cr
&&\cdot&173&214&242&256&256\cr
&&&\cdot &194&214&224&224\cr
&&&&\cdot&173&178&178\cr
&&&&&\cdot&122&122\cr
&&&&&&\cdot&61\cr
&&&&&&&.\cr
}$.

Thus,
$$c_n(m,k)=\cases{(-1)^{(n-1)/2}b_{n+1}(n-k ,n-m),& if $n$ is odd;\cr
(-1)^{n/2}a_{n+1}(n-k ,n-m),& if $n$ is even.\cr}\leqno(8.6)$$
The sequence of triangles $(C_n)$ defined by (8.6) is a Seidel triangle sequence (same argument
as in the proof of Proposition 8.1).
Following the same pattern as in the preceding subsection, 
we
form the matrix~$H$, whose skew-diagonals carry the entries of the leftmost columns of the $C_n$'s:
$$
H=\pmatrix{-1&0&1&0&-5&0&61 & \cdots \cr
0&2&0&-10&0&122\cr
2&0&-14&0&178\cr
0&-16&0&224\cr
-16&0&256\cr
0&272\cr
272\cr
\vdots \cr
}.
$$
This matrix is to be compared with the matrix~$\overline H_2$ (see \S7.2). 
We see that~$H$ is obtained from $\overline H_2$ by transposition and deletion of the first row, so that
$$\leqalignno{
H(x,y)&={D_x}\overline H_2(y,x)
={D_x}\Bigl({-\sinh x\over\cosh(x+y)}\Bigr)\cr
&
={-\cosh x\,\cosh(x+y)+\sinh x\,\sinh(x+y)\over \cosh^2(x+y)}\cr
&={-\cosh y\over \cosh^2(x+y)}.\cr}
$$
Therefore,
$$\leqalignno{
e^x H(x+y,z)&=e^x
{-\cosh z\over \cosh^2(x+y+z)};\cr
\quad e^{Ix}H(Ix+Iy,Iz)&=(\cos x+I\sin x)
{-\cos z\over \cos^2(x+y+z)}.\cr
}
$$

\goodbreak
By using (8.6)
the left-hand side of identity (7.12)  can be computed as follows.
If $n=2l+1$, then $I^{n-2}=(-1)^{l+1}I$ and
$(-1)^{(n-1)/2}=(-1)^l$. Thus,
$I^{n-2}c_n(m,k)=-I\,b_{n+1}(n-k ,n-m)$.
The imaginary part of identity (7.12) reads:
$$\displaylines{
\sum_{\scriptstyle 1\le m+1\le k\le n-1
\atop \scriptstyle n\ {\rm odd}}
\kern-20pt 
b_{n+1}(n-k,n-m){x^{n-k-1}\over (n-k-1)!}{y^{k-m-1}\over (k-m-1)!}{z^m\over m!}
\hfill\cr
\noalign{\vskip-10pt}
\hfill{}=
{\sin x\cos z\over \cos^2(x+y+z)}.\quad\cr}$$
With the change of variables $n\leftarrow 2n-1$, $n-k\leftarrow m$, $n-m\leftarrow k$, we get (1.19) from Theorem 1.5.

\goodbreak
\medskip
If $n=2l$, then $I^{n-2}=(-1)^{l-1}$ and $(-1)^{n/2}=(-1)^l$, so that
$I^{n-2}c_n(m,k)=-a_{n+1}(n-k ,n-m)$.
As for the real, 
$$\displaylines{
\sum_{\scriptstyle 1\le m+1\le k\le n-1
\atop \scriptstyle n\ {\rm even}}
\kern-20pt 
a_{n+1}(n-k;n-m){x^{n-k-1}\over (n-k-1)!}{y^{k-m-1}\over (k-m-1)!}{z^m\over m!}\hfill\cr
\hfill{}=
{\cos x\cos z\over \cos^2(x+y+z)}.\quad\cr}$$
With the change of variables $n\leftarrow 2n$, $n-k\leftarrow m$, $n-m\leftarrow k$, we get (1.17) from Theorem 1.4.

\medskip
8.3. {\it The bottom rows of the matrices $B_n$'s}.\quad 
By Rule (TS4.1) and (2.6) those bottom rows, after discarding the rightmost entry which is always null, read:
$b_2(2,1)=1$; $(b_3(3,1),b_3(3,2))=(0,1)$, 
$(b_4(4,1),b_4(4,2),  \discretionary{}{}{}   b_4(4,3))\discretionary{}{}{} =(0,1,1)$,
$(b_5(5,1),b_5(5,2),b_5(5,3),b_5(5,4))=(0,1,2,2)$, \dots, which are equal to the
sequences of the Entringer numbers: $E_1(1)$, $(E_2(2),E_2(1))$,
$(E_3(3),E_3(2),E_3(1))$, $(E_4(4),E_4(3),E_4(2),E_4(1))$, \dots\ 
By (7.7) and (7.8) we recover the two identities (1.23) and (1.24) written at the end of Section~1.

\medskip
8.4. {\it The lower triangles of ${\rm Twin}^{(1)}$.}\quad
As for the upper triangles, a geometric transformation is to be made to configurate those lower triangles into Seidel triangles.
The bottom rows of the $A_n$'s and  $B_n$'s being discarded, we form the following sequence of triangles: 

\noindent
$C_1=\cdot$;\hskip5pt 
$C_2=\matrice{\cdot&1\cr&\cdot\cr}$;\hskip5pt 
$C_3=\kern-5pt \matrice{\cdot&1&0\cr
&\cdot&1\cr
&&\cdot\cr}$;\quad
$C_4=\kern-5pt \matrice{\cdot&0&-1&-1\cr
&\cdot&-1&-2\cr
&&\cdot&-2\cr
&&&\cdot\cr}$;\quad
$C_5=\kern-5pt \matrice{\cdot&-2&-2&-1&0\cr
&\cdot&-4&-3&-1\cr
&&\cdot&-4&-2\cr
&&&\cdot&-2\cr
&&&&\cdot\cr}$;

\indent\indent
$C_6=\matrice{\cdot&0&2&4&5&5\cr
&\cdot&2&6&9&10\cr
&&\cdot&8&12&14\cr
&&&\cdot&14&16\cr
&&&&\cdot&16\cr
&&&&&\cdot\cr
}$;\quad
$C_7=\matrice{\cdot&16&16&14&10&5&0\cr
&\cdot&32&30&24&15&5\cr
&&\cdot&44&36&24&10\cr
&&&\cdot&44&30&14\cr
&&&&\cdot&32&16\cr
&&&&&\cdot&16\cr
&&&&&&\cdot\cr};$

Thus, for $0\le m<k\le n-1$
$$
c_n(m,k)=\cases{(-1)^{(n+1)/2}a_{n+1}(k+1,m+1),& if $n$ is odd;\cr
(-1)^{(n+2)/2}b_{n+1}(k +1,m+1),& if $n$ is even.\cr}\leqno(8.9)
$$
The sequence of triangles $(C_n)$ defined by (8.9) is a Seidel triangle 
sequence.  The corresponding matrix~$H$ reads:
$$
\leqalignno{
	H&=\pmatrix{1&1&-2&-2&16&16 &\cdots\cr
0&-2&-2&16&16\cr
-1&-1&14&14\cr
0&10&10\cr
5&5\cr
0\cr
\vdots \cr
}\hfill\cr
&=\pmatrix{1&.&-2&.&16 \ \cdots\cr
.&-2&.&16&.\cr
-1&.&14&.\cr
.&10&.\cr
5&.\cr
\vdots\cr}\!+\!
\pmatrix{.&1&.&-2&.&16 \ \cdots\cr
0&.&-2&.&16\cr
.&-1&.&14\cr
0&.&10\cr
.&5\cr
0\cr
\vdots \cr
}\hfill\cr
&=-D_y\overline H_2-\overline H_2.\hfill\cr
}
$$
Thus,
$$\leqalignno{
H(x,y)&=D_y{\sinh y\over \cosh(x+y)}+{\sinh y\over \cosh(x+y)}\cr
&={\cosh x\over \cosh^2(x+y)}+{\sinh y\over \cosh(x+y)};&(8.10)\cr
e^{x}H(x+y,z)&=e^x\Bigl(
{\cosh (x+y)\over \cosh^2(x+y+z)}+{\sinh z\over \cosh(x+y+z)}\Bigr);\cr
e^{Ix}H(Ix+Iy,Iz)&=(\cos x+I\sin x)\Bigl(
{\cos (x+y)\over \cos^2(x+y+z)}+{I\,\sin z\over \cos(x+y+z)}\Bigr).\cr
}
 $$
If $n=2l+1$, then $I^{n-2}=(-1)^{l+1}I$ and
$(-1)^{(n+1)/2}=(-1)^{l+1}$. Thus,
$I^{n-2}c_n(m,k)=I\,a_{n+1}(k+1,m+1)$.
The imaginary part of  identity (7.12) becomes:
$$\displaylines{
\sum_{\scriptstyle 1\le m+1\le k\le n-1
\atop \scriptstyle n\ {\rm odd}}
\kern-20pt 
a_{n+1}(k+1,m+1){x^{n-k-1}\over (n-k-1)!}{y^{k-m-1}\over (k-m-1)!}{z^m\over m!}
\hfill\cr
\hfill{}={\cos x\,\sin z\over \cos(x+y+z)}
+{\sin x\,\cos (x+y)\over \cos^2(x+y+z)}.\quad\cr}$$
With the change of variables $n\leftarrow 2n-1$, $k+1\leftarrow m$, $m+1\leftarrow k$, we get (1.16) from Theorem 1.3.

If $n=2l$, then $I^{n-2}=(-1)^{l-1}$ and $(-1)^{(n+2)/2}=(-1)^{l+1}$, so that
$I^{n-2}c_n(m,k)=b_{n+1}(k+1,m+1)$. 
As for the real part
$$\displaylines{
\sum_{\scriptstyle 1\le m+1\le k\le n-1
\atop \scriptstyle n\ {\rm even}}
\kern-20pt 
b_{n+1}(k+1,m+1){x^{n-k-1}\over (n-k-1)!}{y^{k-m-1}\over (k-m-1)!}{z^m\over m!}
\hfill\cr
\hfill{}=-{\sin x\,\sin z\over \cos(x+y+z)}
+{\cos x\,\cos (x+y)\over \cos^2(x+y+z)}.\quad\cr}$$
With the change of variables $n\leftarrow 2n$, $k+1\leftarrow m$, $m+1\leftarrow k$, we get (1.22) from Theorem 1.6.

\medskip
8.5. {\it The lower triangles of ${\rm Twin}^{(2)}$.}\quad
Again, the bottom rows of the $A_n$'s and $B_n$'s having been discarded, the Seidel Triangle Sequence to be considered is the following:

\noindent
$C_1=\cdot\,$;\quad
$C_2=\matrice{\cdot&1\cr
&\cdot\cr}$;\quad
$C_3=\matrice{\cdot&0&-1\cr
&\cdot&-1\cr
&&\cdot\cr}$;\quad
$C_4=\matrice{\cdot&-1&-1&0\cr
&\cdot&-2&-1\cr
&&\cdot&-1\cr
&&&\cdot\cr}$\quad
$C_5=\matrice{\cdot&0&1&2&2\cr
&\cdot&1&3&4\cr
&&\cdot&4&5\cr
&&&\cdot&5\cr
&&&&\cdot\cr}$;

$C_6=\matrice{\cdot&5&5&4&2&0\cr
&\cdot&10&9&6&2\cr
&&\cdot&13&9&4\cr
&&&\cdot&10&5\cr
&&&&\cdot&5\cr
&&&&&\cdot\cr}
$;\qquad
$C_7=\matrice{\cdot&0&-5&-10&-14&-16&-16\cr
&\cdot&-5&-15&-24&-30&-32\cr
&&\cdot&-20&-33&-42&-46\cr
&&&\cdot&-41&-51&-56\cr
&&&&\cdot&-56&-61\cr
&&&&&\cdot&-61\cr
&&&&&&\cdot\cr}$;

\noindent
the general formula being:
$$
c_n(m,k)=\cases{(-1)^{(n-1)/2}b_{n+1}(k+1,m+1),& if $n$ is odd;\cr
(-1)^{(n-2)/2}a_{n+1}(k+1,m+1),& if $n$ is even.\cr}\leqno(8.13)
$$
Next, form the matrix~$H$, whose skew-diagonals carry the entries of the rightmost columns of the $c_n$'s, and write it as the sum of the following two matrices:
$$\leqalignno{
	&	H=\pmatrix{1&-1&-1&5&5&-61&-61 & \cdots\cr
-1&-1&5&5&-61&-61\cr
0&4&4&-56&-56\cr
2&2&-46&-46\cr
0&-32&-32\cr
-16&-16\cr
0\cr
\vdots \cr
}:=K_1+K_2\hfill\cr
&\!=\!\petitematrice{1&\cdot&-1&\cdot&5&\cdot& \! \! -61 \, \cdots\cr
\cdot&-1&\cdot&5&\ \cdot&-61\cr
0&\cdot&4&\cdot&-56\cr
\cdot&2&\cdot&-46\cr
0&\cdot&-32\cr
\cdot&-16\cr
0\cr
\vdots \cr
}
\!+\!\petitematrice{\cdot&-1&\cdot&5&\cdot&\! -61&\, \cdots\cr
-1&\cdot&5&\cdot&-61&\cdot\cr
\cdot&4&\cdot&-56&\cdot\cr
2&\cdot&-46&\cdot\cr
\cdot&-32&\cdot\cr
-16&\cdot\cr
\vdots\cr
}\!.\cr}
$$
Those matrices are to be compared with the matrix~$\overline H_1$ (see Section~7). Clearly, $K_2$ can be obtained from $\overline H_1$ by deleting the top row and then transposing the matrix, so that $K_2(x,y)=TD_x\overline H_1(x,y)$. Also, $K_1=T\overline H_1$ and then $K_1(x,y)=\overline H_1(y,x)$. As $\overline H_1(x,y)=\cosh y/\cosh(x+y)$, we get:
$$\leqalignno{
H(x,y)&=TD_x\overline H_1(x,y)+\overline H_1(y,x)\cr
&=-{\cosh x\,\sinh (x+y)\over \cosh^2(x+y)}+{\cosh x\over\cosh(x+y)};\cr
e^{x}H(x+y,z)&=e^x\Bigl(-
{\cosh(x+y)\sinh (x+y+z)\over \cosh^2(x+y+z)}+{\cosh (x+y)\over \cosh(x+y+z)}\Bigr);\cr
e^{Ix}H(Ix+Iy,Iz)&\!=\!(\cos x+I\sin x)\cr
&\ {}\times\Bigl(-
{I\,\cos(x+y)\sin (x+y+z)\over \cos^2(x+y+z)}+{\cos (x+y)\over \cos(x+y+z)}\Bigr).\cr
}
$$

If $n=2l+1$, then $I^{n-2}=(-1)^{l+1}I$ and
$(-1)^{(n-1)/2}=(-1)^{l}$. Thus,
$I^{n-2}c_n(m,k)=-I\,b_{n+1}(k+1,m+1)$.
The imaginary part of identity (7.12) becomes:
$$\displaylines{
\sum_{\scriptstyle 1\le m+1\le k\le n-1
\atop \scriptstyle n\ {\rm odd}}
\kern-20pt 
b_{n+1}(k+1,m+1){x^{n-k-1}\over (n-k-1)!}{y^{k-m-1}\over (k-m-1)!}{z^m\over m!}
\hfill\cr
\hfill{}=-{\sin x\,\cos (x+y)\over \cos(x+y+z)}
+{\cos x\,\cos(x+y)\sin(x+y+z)\over \cos^2(x+y+z)}
\quad\cr
\hfill{}={\cos(x+y)\,\sin(y+z)\over \cos^2(x+y+z)}.\kern4.8cm\cr}$$
With the change of variables $n\leftarrow 2n-1$, $m+1\leftarrow k$, $k+1\leftarrow m$, we get (1.20) from Theorem 1.5.

If $n=2l$, then $I^{n-2}=(-1)^{l-1}$ and $(-1)^{(n-2)/2}=(-1)^{l-1}$, so that
$I^{n-2}c_n(m,k)=a_{n+1}(k+1,m+1)$. 
As for the real part, 
$$\displaylines{
\sum_{\scriptstyle 1\le m+1\le k\le n-1
\atop \scriptstyle n\ {\rm even}}
\kern-20pt 
a_{n+1}(k+1,m+1){x^{n-k-1}\over (n-k-1)!}{y^{k-m-1}\over (k-m-1)!}{z^m\over m!}
\hfill\cr
\hfill{}={\cos x\,\cos (x+y)\over \cos(x+y+z)}
+{\sin x\,\cos(x+y)\sin(x+y+z)\over \cos^2(x+y+z)}
\quad\cr
\hfill{}={\cos(x+y)\,\cos(y+z)\over \cos^2(x+y+z)}.\kern4.5cm \cr}$$
With the change of variables
$n\leftarrow 2n-1$, $k+1\leftarrow m$, $m+1\leftarrow k$, we get (1.18) from Theorem 1.4.

\medskip

\centerline{\bf 9. The formal Laplace transform}

\medskip
The purpose of this Section is to show that, when the Entringer numbers $E_n(k)$ are defined by relations (1.1), without any reference to their combinatorial interpretations, they can be proved to be a refinement of the tangent/secant numbers: $\sum_kE_n(k)=E_n$ $(n\ge 1)$. In the same manner, when the twin Seidel matrix sequence $(A_n)$, $(B_n)$ is analytically defined, as it was stated in \S\thinspace 1.5, also without reference to any combinatorial interpretation, their entries $(a_n(m,k))$, $(b_n(m,k))$ make up a refinement of the Entringer numbers, by row and by column, and then $\sum_{m,k}a_n(m,k)=\sum_{m,k}b_n(m,k)=E_n$.
The proofs of those results make use of the closed expressions found for the generating functions  obtained in the preceding section, and of a well-adapted formal Laplace transform technique. 

\proclaim Theorem 9.1.
$(1)$ Let $(E_n(k))$ be the sequence of the Entringer numbers, defined by
$$\leqalignno{\noalign{\vskip-5pt}
E_1(1):=1;\quad E_n(n)&:=0\ {\rm for\ all\ }n\ge 2;\cr
\Delta E_n(m)+E_{n-1}(n-m)&=0\quad
(n\geq 2; m=n-1, \ldots, 2,1);\cr
\noalign{\vskip-5pt}}
$$
Then,
$$
\sum_{1\leq k \leq 2n-1 } E_{2n-1}(k) = E_{2n-1};\quad
\sum_{1\leq k \leq 2n} E_{2n}(k) = E_{2n};\quad (n\ge 1).\leqno(9.1)
$$
\null\hskip2.8cm $(2)$ Let $(a_n(m,k))$, $(b_n(m,k))$ be the entries
of the twin Seidel matrix sequence $(A_n)$, $(B_n)$, as they are defined in \S\thinspace $1.5$.
Then,
$$\leqalignno{\qquad
a_n(m,\brullet)&=E_n(m),\quad b_n(m,\brullet)=E_n(n+1-m),\quad
(1\le m\le n);&(9.2)\cr
a_n(\brullet,k)&=b_n(\brullet,k)=E_n(n-k)\quad (1\le k\le n).&(9.3)\cr
}
$$

The proof of (9.1) is fully given. Next, we reproduce the proof of $a_{2n}(m,\brullet)=E_{2n}(m)$, based on Theorem~1.3. The other identities in (9.2) and (9.3) can also be derived following the same method by using Theorems 1.4, 1.5, 1.6. Their proofs are omitted.

The formal Laplace transform, already used in our previous paper [FH14], maps a function~$f(x)$ onto a function ${\cal L}(f(x), x,s)$ defined by
$$\leqalignno{
{\cal L}(f(x), x,s) :&= \int_0^\infty f(x) e^{-xs}\, dx.\cr
\noalign{\hbox{In particular, ${\cal L}(\brullet,x,s)$ maps $x^k/k!$ onto $1/s^{k+1}$:}}
{\cal L}({x^k\over k!}, x,s) &= {1\over s^{k+1}}.\cr} 
$$

For proving (9.1) start with identity (7.7) involving the generating function for the numbers $E_{2n+1}(k)$ and
apply the Laplace transform twice with respect to
$(x,s)$, $(y,t)$
respectively. We get:
$$
\sum_{1\leq k \leq 2n+1 } 
{1\over  s^{2n-k+2}} {1 \over t^{k}} E_{2n+1}(k)
=  \int_0^\infty\!\!\int_0^\infty
{\cos y\over \cos(x+y)} e^{-xs-yt} dx\,dy,
$$
which becomes, with $t\leftarrow s$ and $r=x+y$
$$
\leqalignno{
\sum_{1\leq k \leq 2n+1 } 
{1\over  s^{2n+2}}  E_{2n+1}(k)
&=  \int_0^\infty\!\!\int_0^\infty {\cos y\over \cos(x+y)} e^{-xs-ys} dx\,dy \cr
&=  \int_0^\infty\!\!\int_0^r {\cos y\over \cos r} e^{-rs} dy\,dr \cr
&=  \int_0^\infty {\sin r\over \cos r} e^{-rs} dr \cr
&=  \int_0^\infty ({\tan r}) e^{-rs} dr \cr
&=  \sum_{n\geq 1} {1\over s^{2n}} E_{2n-1}. \cr
\noalign{\vskip-5pt}
\noalign{\hbox{Hence,}}
\noalign{\vskip-5pt}
\sum_{1\leq k \leq 2n-1 } E_{2n-1}(k) &= E_{2n-1}.\cr}
$$

In the same manner, 
apply the Laplace transform to identity (7.8) twice with respect to
$(x,s)$, $(y,t)$
respectively. We get
$$
\sum_{1\leq k \leq 2n } 
{1\over  s^{k}} {1 \over t^{2n-k+1}} E_{2n}(k)
=  \int_0^\infty\!\!\int_0^\infty
{\sin y\over \cos(x+y)} e^{-xs-yt} dx\,dy,
$$
which becomes with $s\leftarrow t$ and $r=x+y$:
$$
\leqalignno{
\sum_{1\leq k \leq 2n } 
{1\over  t^{2n+1}}  E_{2n}(k)
&=  \int_0^\infty\!\!\int_0^\infty {\sin y\over \cos(x+y)} e^{-xt-yt} dx\,dy \cr
&=  \int_0^\infty\!\!\int_0^r {\sin y\over \cos r} e^{-rt} dy\,dr \cr
&=  \int_0^\infty { 1-\cos r \over \cos r} e^{-rt} dr \cr
&=  \int_0^\infty (\sec r -1) e^{-rt} dr \cr 
&=  \sum_{n\geq 1} {1\over t^{2n+1}} E_{2n}.\cr
\noalign{\vskip-5pt}
\noalign{\hbox{Hence,}}
\noalign{\vskip-5pt}
\sum_{1\leq k \leq 2n} E_{2n}(k) &= E_{2n}.\cr}
$$

\goodbreak
Next, to prove $a_{2n}(m,\brullet)=E_{2n}(m)$ start with identity
(1.15) of Theorem~1.5
 and
apply the Laplace transform to its left-hand side three times with respect to
$(x,s)$, $(y,t)$, $(z,u)$,
respectively. We get
$$\displaylines{
\sum_{2\leq m+1\leq k \leq 2n-1 } {1\over s^{m}}   
{1\over  t^{k-m}} {1 \over u^{2n-k}} a_{2n}(m,k),\cr
\noalign{\hbox{which becomes}}
\rlap{(9.4)}\hfill
\sum_{2\leq m+1\leq k \leq 2n-1 } {1\over s^{m}}   
 {1 \over u^{2n}} a_{2n}(m,k),\hfill\cr}$$
when $t\leftarrow u$ and $s\leftarrow su$.
Apply the Laplace transform to the right-hand side of (1.15) three times with respect to
$(x,s)$, $(y,t)$, $(z,u)$, respectively, and let 
$t\leftarrow u,\, s\leftarrow su$. With $r=y+z$ we get:
$$
\displaylines{
\qquad \int_0^\infty\!\!\int_0^\infty\!\!\int_0^\infty
{\cos x\,\cos z\sin (x+y+z)\over \cos^2(x+y+z)}
	e^{-xsu-yu-zu}dx\,dy\,dz\hfill\cr
\qquad \qquad\qquad{}=  \int_0^\infty\!\!\int_0^\infty\!\!\int_0^r
	{\cos x \cos z \sin(x+r)\over \cos^2(x+r)} e^{-xsu-ru}dz\,dr\,dx
 \hfill\cr
\rlap{(9.5)}\qquad\qquad\qquad{}=  \int_0^\infty\!\!\int_0^\infty
 {\cos x \sin r \sin(x+r)\over \cos^2(x+r)} e^{-xsu-ru} dr\,dx. \hfill
 \cr
} 
$$

With identity (1.16) 
apply the Laplace transform to its left-hand side three times with respect to
$(x,u)$, $(y,s)$, $(z,t)$,
respectively. We get
$$\displaylines{
\sum_{2\leq k+1\leq m \leq 2n-1 } {1\over u^{2n-m}}   
{1\over  s^{m-k}} {1 \over t^{k}} a_{2n}(m,k),\cr
\noalign{\hbox{which becomes}}
\rlap{(9.6)}\hfill
\sum_{2\leq k+1\leq m \leq 2n-1 } {1\over s^{m}}   
 {1 \over u^{2n}} a_{2n}(m,k),\hfill\cr}$$
when $s\leftarrow su$ and $t\leftarrow su$.
Apply the Laplace transform to the right-hand side of (1.16) three times with respect to
$(x,u)$, $(y,s)$, $(z,t)$, respectively, and let 
$s\leftarrow su,\, t\leftarrow su$. With $r=y+z$ we get:
$$
\displaylines{
	\int_0^\infty\!\!\int_0^\infty\!\!\int_0^\infty
	\Bigl({\cos x\,\sin z\over \cos(x+y+z)}
	+{\sin x\,\cos (x+y)\over \cos^2(x+y+z)}\Bigr)
	e^{-xu-ysu-zsu}dx\,dy\,dz\hfill\cr
\quad{}=	\int_0^\infty\!\!\int_0^\infty\!\!\int_0^r
	\Bigl({\cos x\,\sin z\over \cos(x+r)}
	+{\sin x\,\cos (x+r-z)\over \cos^2(x+r)}\Bigr)
	e^{-xu-rsu}dz\,dr\,dx\hfill\cr
\quad{}=	\int_0^\infty\!\!\int_0^\infty
	\Bigl({\cos x\,(1-\cos r)\over \cos(x+r)}
	+{\sin x\,(\sin(x+r) -\sin x)\over \cos^2(x+r)}\Bigr)
	e^{-xu-rsu}dr\,dx\hfill\cr
\rlap{(9.7)}\hskip5pt \qquad{}=\!\!	\int_0^\infty\!\!\int_0^\infty\!\!
	\Bigl({\cos r\,(1\!-\!\cos x)\over \cos(x+r)}
	+{\sin r\,(\sin(x+r) \!-\!\!\sin r)\over \cos^2(x+r)}\Bigr)
	e^{-ru-xsu}dr\,dx.\hfill\cr
	} 
$$
By (9.4)---(9.7) we have
$$(9.8)
\sum_{1\leq k,m \leq 2n-1; \, k\not=m } {1\over s^{m}}   
 {1 \over u^{2n}} a_{2n}(m,k) 
=  \int_0^\infty\!\!\int_0^\infty
F(x,r) e^{-xsu-ru} dr\,dx,
$$
where
$$
\leqalignno{
	F(x,r)&=
 {\cos x \sin r \sin(x+r)\over \cos^2(x+r)} +  
	{\cos r\,(1-\cos x)\over \cos(x+r)}
	+{\sin r\,(\sin(x+r) -\sin r)\over \cos^2(x+r)}\cr
&=
 {\cos x \over \cos^2(x+r)} -1. \cr
}
$$
But from  (7.8) 
$$\leqalignno{
\sum_{1\le m\le 2n}E_{2n}(m){x^{m-1}\over (m-1)!}{r^{2n-m-1}\over (2n-m-1)!}
&={\partial \over \partial r}{\sin r\over \cos(x+r)}\cr
&={\cos x\over \cos^2(x+r)}
.&(9.9)\cr}
$$
Apply the Laplace transform to (9.9) twice with respect to
$(x,s)$, $(y,u)$,
respectively. We get:
$$\leqalignno{
\sum_{1\leq m \leq 2n } 
{1\over  s^{m}} {1 \over u^{2n-m}} E_{2n}(m)
&=  \int_0^\infty\!\!\int_0^\infty
{\cos x\over \cos^2(x+r)} e^{-xs-ru} dx\,dr,\cr
\noalign{\hbox{or still}}
\sum_{1\leq m \leq 2n } 
{1\over  s^{m}} {1 \over u^{2n}} E_{2n}(m)
&=  \int_0^\infty\!\!\int_0^\infty
{\cos x\over \cos^2(x+r)} e^{-xsu-ru} dx\,dr. &(9.10)\cr}
$$
By (9.8) and (9.10) we obtain
$$\leqalignno{
\sum_{1\leq k,m \leq 2n-1; \, k\not=m } {1\over s^{m}}   
 {1 \over u^{2n}} a_{2n}(m,k) 
&= \sum_{1\leq m \leq 2n } 
{1\over  s^{m}} {1 \over u^{2n}} E_{2n}(m)
- {1\over su^2}\cr
\noalign{\hbox{and then}}
\sum_{1\leq k,m \leq 2n-1} {1\over s^{m}}   
 {1 \over u^{2n}} a_{2n}(m,k) 
&= \sum_{1\leq m \leq 2n } 
{1\over  s^{m}} {1 \over u^{2n}} E_{2n}(m).\cr}
$$
Hence,
$$
\sum_{1\leq k\leq 2n-1 } a_{2n}(m,k) = E_{2n}(m).\qed
$$

\vskip 1cm

\centerline{\bf References}

{\eightpoint

\bigskip 
\article  An1879|D\'esir\'e Andr\'e|D\'eveloppement de $\sec x$ et
$\tan x$|C. R. Math. Acad. Sci. Paris|88|1879|965--979|

\article An1881|D\'esir\'e Andr\'e|Sur les permutations
altern\'ees|J. Math. Pures et Appl.|7|1881|167--184|

\article AF80|George Andrews; Dominique Foata|Congruences for the
$q$-secant number|Europ. J. Combin.|1|1980|283--287|

\article AG78|George Andrews; Ira Gessel|Divisibility
properties of the $q$-tangent numbers|Proc. Amer. Math.
Soc.|68|1978|380--384|

\livre Co74|Louis Comtet|\hskip-5pt Advanced 
Combinatorics|\hskip-5pt D.
Reidel/Dordrecht-Holland, Boston, {\oldstyle 1974}|

\divers Du82|Dominique Dumont|Matrices d'Euler-Seidel,
{\sl S\'eminaire Lotharin\-gien de Combinatoire}, B05c (1981), 25 pp. [Formerly: Publ. I.R.M.A. Strasbourg, 1982, 182/S-04, p. 59-78.]
\hfil\break
{\tt http://www.mat.univie.ac.at/$\sim$slc/}|

\divers DV80|Dominique Dumont; G\'erard Viennot|A
combinatorial interpretation of the Seidel generation of Genocchi
numbers, {\sl Combinatorial mathematics, optimal designs} [J.
Srivastava, ed., Fort Collins. {\oldstyle 1978}], p.~77--87\pointir
Amsterdam, North-Holland, {\oldstyle 1980} ({\sl Annals of
Discrete Math.} {\bf 6})|

\divers Di74|Filippo Disanto|Andr\'e Permutations, Right-To-Left and Left-To-Right Minima, {\sl S\'eminaire Lothar. Combin.}, B70f, {\oldstyle 2014}, 13 pp.\hfil\break
{\tt http://www.mat.univie.ac.at/$\sim$slc/}|

\article En66|R. C. Entringer|A combinatorial interpretation of the
Euler and Bernoulli numbers|Nieuw. Arch. Wisk.|14|1966|241--246|

\article FH01|Dominique Foata; Guo-Niu Han|Arbres minimax et polyn\^omes 
d'Andr\'e|Advances in Appl. Math.|27|2001|367--389|

\article FH13|Dominique Foata; Guo-Niu Han|Finite Difference Calculus
for Alternating Permutations|J.  Difference Equations and Appl.|19|2013|1952--1966|

\article FH14|Dominique Foata; Guo-Niu Han|Seidel triangle sequences and Bi-Entringer numbers|Europ. J. of Combin.|42|2014|243-260|

\divers FSch71|Dominique Foata; Marcel-Paul Sch\"utzenberger|Nombres d'Euler et
permutations alternantes. Manuscript , 71 pages, University of
Florida, Gainesville, {\oldstyle 1971},\hfill\break
{\tt http://www.mat.univie.ac.at/$\sim$slc/}|

\divers FSch73|Dominique Foata; Marcel-Paul Sch\"utzenberger|Nombres d'Euler et
permutations alternantes, in J. N. Srivastava et {\it al.} (eds.), {\it
A Survey of Combinatorial Theory}, North-Holland, Amsterdam, {\oldstyle 1973},
pp. 173-187|

\article FSt74|Dominique Foata; Volker Strehl|Rearrangements of the symmetric
group and enumerative properties of the tangent and secant
numbers|Math. Z.|137|1974|257-264|

\divers FSt76|Dominique Foata; Volker Strehl|Euler numbers and variations of
permutations, in {\it Colloquio Internazionale sulle Teorie
Combinatorie, {\oldstyle1973}}, vol.~I ({\it Atti dei Convegni Lincei}, {\bf
17}, 119-131). Accademia Nazionale dei Lincei, {\oldstyle 1976}|

\article GHZ11|Yoann Gelineau; Heesung Shin; Jiang Zeng|Bijections for
Entringer families|Europ. J. Combin.|32|2011|100--115|

\article He96|G\'abor Hetyei|On the cd-variation polynomials of Andr\'e and Simsun permutations|Discrete Comput.
Geom.|16|1996|259Ð275|

\article HR98|G\'abor Hetyei; Ethan Reiner|Permutation Trees and Variation
Statistics|Europ. J. Combin.|19|1998|847-866|

\livre Jo39|Charles Jordan|Calculus of Finite Differences|R\"ottig and
Romwalter,  Budapest, {\oldstyle 1939}|

\article KPP94|A. G. Kuznetsov; I. M. Pak;  A. E.
Postnikov|Increasing trees and alternating permutations|Uspekhi Mat.
Nauk|49|1994|79--110|

\article MSY96|J. Millar; N. J. A. Sloane; N. E. Young|A new operation on sequences: the Boustrophedon transform|J. Combin. Theory ser.~A|17|1996|44-54|

\livre Ni23|Niels Nielsen|Trait\'e \'el\'ementaire des nombres
de Bernoulli|Paris, Gauthier-Villars, {\oldstyle 1923}|

\divers OEIS|OEIS Foundation|Sequence A008282, {\it The On-Line Encyclopedia of Integer Sequences}, {\oldstyle 2015}, {\tt http://oeis.org}| 

\article Po89|Christiane Poupard|Deux propri\'et\'es des arbres
binaires ordonn\'es stricts|Europ. J. Combin.|10|1989|369--374|

\article Pu93|Mark Purtill|Andr\'e permutations, lexicographic shellability
and the {\it cd}-index of a convex polytope|Trans. Amer.
Math. Soc.|338|1993|77-104|

\divers Se1877|L. Seidel|\"Uber eine einfache Enstehungsweise der Bernoullischen Zahlen und einiger verwandten Reihen, {\sl Sitzungberichte der M\"unch. Akad. Math. Phys. Classe}, {\oldstyle 1877}, p.~157--187|

\article St76|Richard P. Stanley|Binomial posets, M\"obius inversion, and permutation enumeration|J. Combin. Theory Ser. A|20|1976|336-356|

\article St94|Richard P. Stanley|Flag f-vectors and the cd-index|Math. Z.|216|1994|483--499|

\divers St10|Richard P.
Stanley|A Survey of Alternating Permutations, in {\sl Combinatorics and graphs}, 165--196, {\sl Contemp. Math.}, {\bf 531}, Amer. Math. Soc. Providence, RI, {\oldstyle 2010}|

\divers Str74|Volker Strehl|Geometrische und arithmetische Eigenschaften von Andr\'e-Poly\-nomen, Ph. D. Thesis, Friedrich-Alexander-Universit\"at Erlangen-N\"urnberg, {\oldstyle 1974}|

\article Vi80|G\'erard Viennot|Une interpr\'etation
combinatoire des coefficients des d\'evelop\-pements en
s\'erie enti\`ere des fonctions elliptiques de Jacobi|J.
Combin. Theory, Series
A|29|1980|121--133|

\divers Vi88|Xavier G. Viennot|S\'eries g\'en\'eratrices
\'enum\'eratives, chap.~3, Lecture Notes, 160~p., 1988, notes de
cours donn\'es
\`a l'\smash{\'E}cole Normale Sup\'erieure Ulm (Paris), UQAM (Montr\'eal,
Qu\'ebec) et Universit\'e de Wuhan (Chine)\hfil\break
{\tt
http://www.xavierviennot.org/xavier/cours.html}|

\bigskip
\hbox{\vtop{\halign{#\hfil\cr
Dominique Foata \cr
Institut Lothaire\cr
1, rue Murner\cr
F-67000 Strasbourg, France\cr
\noalign{\smallskip}
{\tt foata@unistra.fr}\cr}}
\qquad
\vtop{\halign{#\hfil\cr
Guo-Niu Han\cr
I.R.M.A. UMR 7501\cr
Universit\'e de Strasbourg et CNRS\cr
7, rue Ren\'e-Descartes\cr
F-67084 Strasbourg, France\cr
\noalign{\smallskip}
{\tt guoniu.han@unistra.fr}\cr}}}

}

\bye